\documentclass[preprint,12pt]{elsarticle}

\usepackage{natbib}
\usepackage{multirow}
\usepackage{amssymb}
\usepackage{amsmath}
\usepackage[table,xcdraw]{xcolor}
\usepackage{amsthm}
\usepackage{ragged2e}
\usepackage{bm}
\usepackage{cancel} 
\newtheorem*{remark}{Remark}
\newtheorem*{defA*}{Subproblem A}
\newtheorem*{defB*}{Subproblem B}
\newtheorem*{MP*}{Material Parameters}
\usepackage{cleveref}
\def \vec#1{{\bf{#1}}}
\crefrangeformat{equation}{(#3#1#4--#5#2#6)}
\usepackage[ruled]{algorithm2e}
\usepackage{comment}
\usepackage{longtable}
\usepackage{anyfontsize}
\usepackage{pgfplots} 
\pgfplotsset{width=6.5cm, compat=1.6}
\usepackage{lineno}
\usepackage{graphicx}
\usepackage{wrapfig}
\DeclareMathOperator{\tr}{tr}
\usepackage{subcaption}
\DeclareMathOperator*{\argmin}{arg\,min}
\newcommand{\mc}{\mathcal}

\newcommand{\be}{\begin{equation}}
\newcommand{\ee}{\end{equation}}

\newcommand{\bi}{\begin{itemize}}
\newcommand{\ei}{\end{itemize}}

\usepackage{tikz-3dplot}
\usepackage{pgfplots}
\usetikzlibrary{3d} 

\journal{Computers $\&$ Mathematics with Applications}

\begin{document}
\begin{frontmatter}
\title{Nonlinear Methods for Shape Optimization Problems in Liquid Crystal Tactoids}
\author[label1]{J. H. Adler\corref{cor1}}\ead{James.Adler@tufts.edu}
\author[label1]{A. S. Andrei}\ead{Anca.Andrei@tufts.edu}
\address[label1]{Department of Mathematics, Tufts University, Medford, MA 02155}
\author[label2]{T. J. Atherton} \ead{Timothy.Atherton@tufts.edu}
\address[label2]{Department of Physics and Astronomy, Tufts University, Medford, MA 02155}
\cortext[cor1]{Corresponding author}

\begin{abstract}
Anisotropic fluids, such as nematic liquid crystals, can form non-spherical equilibrium shapes known as tactoids. Predicting the shape of these structures as a function of material parameters is challenging and paradigmatic of a broader class of problems that combine shape and order. Here, we consider a discrete shape optimization approach with finite elements to find the configuration of two-dimensional and three-dimensional tactoids using the Landau--de Genne framework and a Q-tensor representation. Efficient solution of the resulting constrained energy minimization problem is achieved using a quasi-Newton and nested iteration algorithm. Numerical validation is performed with benchmark solutions and compared against experimental data and earlier work. We explore physically motivated subproblems, whereby the shape and order are separately held fixed, respectively, to explore the role of both and examine material parameter dependence of the convergence. Nested iteration significantly improves both the computational cost and convergence of numerical solutions of these highly deformable materials.
\end{abstract}


\end{frontmatter}

\section{Introduction}\label{Introduction}

Liquid crystals (LCs) are intermediate phases of matter that exhibit long-range order like a crystal but retain fluid properties \cite{DeGennes1993}. The \emph{nematic} LC phase, in particular, lacks translational order but possesses orientational order characterized by a locally-preferred axis of molecular or particulate alignment; this direction may vary spatially at the cost of elastic energy. Due to the presence of orientational order, nematic liquid crystals possess anisotropic physical properties, such as surface tension, dielectric response, and elasticity. In contrast to an isotropic fluid, which only exhibits surface tension and no elastic effects, LCs may form non-spherical droplets known as \emph{tactoids} \cite{Prinsen2003} when suspended in a surrounding host isotropic fluid. Tactoids can assume various shapes and director field configurations depending on their size \cite{Prinsen2003}, elastic properties, and anisotropic surface tension strength \cite{Chandrasekhar1966}. Due to their potential to change shape and ability to conform to complex geometries \cite{Dzubiella2000, Sitta2018}, tactoids are an exciting geometry for emerging technologies \cite{Lagerwall2012}. This includes enhancing LC displays’ performance \cite{Asdonk2017}, serving as carriers for pharmaceuticals \cite{Stealey2023, Mascarenhas2023}, and developing materials with adaptive stiffness as seen in soft robots \cite{Wehner2016, Shah2021, Schwarzendahl2021}.

In this paper, we develop efficient and robust numerical methods based on nonlinear optimization techniques to predict the solution to tactoid shape-order problems. Using this approach, we compare our results with earlier work, Bates \cite{Bates2010} and Prinsen and van der Schoot \cite{Prinsen2003}, all while improving upon previous research on modeling tactoids' various morphologies. 

Numerical efforts to investigate the configuration of a nematic tactoid droplet have received extensive attention. Monte Carlo simulations have shown that a tactoid's aspect ratio can be temperature dependent \cite{Bates2003} and that their morphologies can depend on the LC’s orientational ordering \cite{Xing2012}. For instance, various tactoid formations are determined by competition between the bending and surface tension energies as LCs exhibit phase transitions \cite{Ding2021}. Monte Carlo methods have also been used to model tactoid defects \cite{Dzubiella2000, Ludwig2020, Sitta2018}. 
While Monte Carlo methods are versatile, they are computationally expensive and require many simplifying assumptions on the model to achieve convergence. Phase field methods use an auxiliary scalar field to interpolate between the interior and exterior of a shape and have been used for modeling tactoids \cite{Ludwig2020}. Such methods are powerful, but challenges arise when dealing with cusps in the manifold \cite{Chen2002}. Level set methods, which represent the free boundary of the system as a contour or a level set of a scalar function defined in a higher-dimensional space, have been used to model interfaces of tactoids and depict their defects \cite{Cermelli2007}. While level set methods can model materials that change shape and topology, they require sophisticated numerical techniques, and including constraints can be challenging \cite{Gibou2018}.

In contrast, finite--element discretizations with gradient descent (GD) based algorithms have been designed to iteratively adjust the position and orientational degrees of freedom in the tactoid to find stable nematic tactoid solutions under different conditions \cite{Nitschke2020, Schimming2021}. DeBenedictis \emph{et al.} \cite{DeBenedictis2016} use a director formulation and Frank--Oseen energy to develop an alternating optimization scheme that takes GD iterations to solve for the director configuration, then solves for the optimal shape and repeats until convergence. While GD methods are computationally cheap per iteration using gradient-only calculations and are easy to implement, they possess linear convergence and hence need a high number of iterations to converge.

This paper aims to improve upon \cite{DeBenedictis2016} by developing an integrated optimization method that simultaneously determines the shape and director configuration while ensuring physical validity. Predicting the optimal shape and physical fields involves solving a nonlinearly constrained optimization problem where we minimize the sum of bulk terms defined on a manifold, $\mathcal{M}$, and surface terms defined on the boundary $\partial\mathcal{M}$ while satisfying a nonlinear volumetric or surface area constraint. To expedite convergence, we use Newton’s method with Lagrange multipliers \cite{Adler2015, Nocedal1999}, which offers local quadratic convergence and fewer iterations for faster solutions. However, Newton’s method requires calculation of the Hessian matrix of all functionals concerning shape and field and is sensitive to initial guesses, potentially hindering convergence to the correct solution. To mitigate computational expense, we approximate the Hessian using the well-established BFGS quasi-Newton method, sacrificing quadratic for local superlinear convergence. The method is implemented and run for physically-relevant parameters, reproducing several expected physical phenomena. We use nested iteration and Newton damping with line search \cite{Broyden1973, Nocedal1999} to handle relatively poor initial guesses for efficient iterative convergence. Nested iteration techniques are a hierarchical approach to solving complex numerical problems, where coarse-grid solutions are used to accelerate convergence on finer grids, improving computational efficiency \cite{Starke2000, Briggs2000, Trottenberg2000}. Nested iteration has been successfully applied to a variety of problems, including LC optimization problems \cite{Adler2016}. Here, we demonstrate its efficiency in resolving the above-mentioned tactoid shapes.  Many preconditioning techniques \cite{Wathen_2015} would also improve the efficiency of the gradient descent and quasi-Newton algorithms we present.  For example, using an appropriate Riesz map, either on the gradient or Hessian, may yield mesh-independent results.  However, the use of nested iteration reduces the need for such methods, as most computation is done on coarser grids, resulting in only a few iterations needed per grid level at high resolution.

Finally, we note that a great amount of work has been done on solution methods for solving LC problems in fixed boundaries. This includes deflation and parameter continuation for finding multiple stable LC configurations \cite{Farrell2015, Adler2017, Emerson2018, Xia2021} and multigrid preconditioners for the resulting linear systems \cite{Ramage2013, Nochetto2017, MacDonald2020, Xia2023} to improve performance. As our focus is on the nonlinear methods with nested iteration, we use direct solvers when needed and consider multigrid methods as future work.

The remainder of the paper is organized as follows. We first pose the shape optimization problem and construct its discrete version in Section \ref{TactoidEnergyModel}. In Section \ref{discreteTactoidSection}, we derive the GD method from \cite{DeBenedictis2016}, the all-at-once quasi-Newton-based method, and describe the nested iteration approach used to improve computational efficiency in finding the optimal tactoid shape. Numerical experiments are reported in Section \ref{results}. Together with the full problem, we study two physically motivated subproblems, whereby either the shape or liquid crystal itself is held fixed, to understand the role of each better. We also study the formation of three-dimensional nematic tactoids and explain their similarities and differences to the two-dimensional case. We give concluding remarks in Section \ref{conclude} and consider opportunities for future work.
	

\section{Q-tensor Model for Tactoids}\label{TactoidEnergyModel}
The configuration of a nematic liquid crystal is described by a vector or tensor field that encodes information about the local orientational ordering of the constituent molecules or particles. Several choices of representation are commonly used. The first possibility is to represent the local average molecular orientation by a unit vector field known as the \emph{director}, $\vec{n}$. The director fully describes the nematic in the absence of \emph{disclinations}, special points where $\vec{n}$ is not uniquely defined and $\nabla \vec{n}$ diverges. If the director formulation is used, $\vec{n}$ is a minimizer of the Frank--Oseen free energy model \cite{Frank1958, Adler2015} subject to the constraint that $\vec{n}\cdot\vec{n} = 1$ everywhere.

An alternative formulation, known as the Q-tensor approach \cite{Mottram2014, DeGennes1993}, encodes both orientational information as well as the local degree of alignment, denoted by a scalar field $S$, into a single tensor order parameter. In $d-$dimensions, and assuming that the alignment is uniaxial, the Q-tensor has the form,
\begin{equation}\label{Q-tensor}
    \mathcal{Q} = S\left(\vec{n}\otimes\vec{n} - \frac{\mathcal{I}}{d}\right),
\end{equation}
where $\mathcal{I}$ is the identity matrix.

In this paper, we model nematic tactoids using a three-dimensional $\mathcal{Q}$ where, by choice of parameterization, $\mathcal{Q}$ is symmetric and traceless.
In the isotropic phase, where there is no orientational order, $\mathcal{Q}=0$, i.e., a zero-tensor. Given an instance of $\mathcal{Q}$, both $\vec{n}$ and $S$ can be reconstructed by eigenanalysis: the largest eigenvalue of $\mathcal{Q}$ is $\frac{2}{3}S$ and $\vec{n}$ is the associated normalized eigenvector. If $\mathcal{Q}$ is uniaxial, the eigenvalues lie on the interval $[-1/3,2/3]$ \cite{Mottram2014}.  While we do not enforce uniaxiality explicitly, we will verify that most of the solutions we find possess eigenvalues of $\mathcal{Q}$ on this interval.  Some three-dimensional tests, however, show the emergence of biaxial configurations.

The Q-tensor approach has a number of advantages over the director formulation. For one, it permits defects since when large gradients of $\vec{n}$ arise, $S$ compensates by tending to zero; it can accommodate defects with a non-integer winding number naturally. Moreover, it obviates the need to impose any unit-length constraints. Generally, these advantages are at the expense of requiring more degrees of freedom overall. 

The local values of $\mathcal{Q}$ are obtained by minimizing a free energy that includes bulk terms defined on a manifold, $\mathcal{M}$, and surface terms defined on the boundary $\partial\mathcal{M}$. The free energy has the form,
\begin{equation} \label{qtensorEnergy}
    \mathcal{F}(\vec{x}, {\mathcal{Q}},\nabla{\mathcal{Q}}, ...)=\int_{\mathcal{M}}f(\vec{x}, {\mathcal{Q}},\nabla {\mathcal{Q}}, ...)\ d\vec{x} + \int_{\partial \mathcal{M}}g(\vec{x}, {\mathcal{Q}}, \nabla {\mathcal{Q}},...)\ dS, 
\end{equation}
where $\vec{x} = (x_1,x_2,x_3) \in \mathcal{M}$ and $f$ and $g$ are linear or nonlinear energy densities that depend on ${\mathcal{Q}}$ and its derivatives.
 The function $g$ is defined on $\partial\mathcal{M}$ and may impose a preferred orientation of the LC relative to the boundary tangent plane, a phenomenon referred to as \emph{anchoring}. Furthermore, the minimization may be subject to nonlinear \emph{global} (integral) equality constraints, such as one that fixes the manifold's volume,
\begin{align} \label{areaConstraint}
    \mathcal{C}({\mathcal{Q}},\nabla {\mathcal{Q}}) = \int_{\mathcal{M}}c({\mathcal{Q}},\nabla {\mathcal{Q}})\ d \vec{x} = 0.
\end{align}
\subsection{Landau--de Gennes Energy Model}
Absent any external forces, we specifically consider the one constant Landau--de Gennes model, where the free energy sufficient to capture the physics of nematic tactoids, $\mathcal{F}$, depends on the state variables of the system: the shape $\mathcal{M}$ and the tensor $\mathcal{Q}(\vec{x})$. 
The free energy is,
 \begin{align}\label{QEnergy}
		\mathcal{F}(\vec{x}, \mathcal{Q}) &= \int_{\mathcal{M}} a\tr{({\mathcal{Q}}^2)} + \frac{2b}{3}\tr{({\mathcal{Q}}^3)} + \frac{c}{2}\tr{({\mathcal{Q}}^2)}^2 + \frac{L_1}{2}|\nabla{{\mathcal{Q}}}|^2\ d \vec{x} \nonumber \\
  &+ \int_{\partial\mathcal{M}} \sigma - \frac{W}{2}\tr{\big(({\mathcal{Q}}-{\mathcal{Q}}_s)^2\big)}\ dS,
\end{align}
where the first three bulk terms represent a Landau expansion of the free energy of ${\mathcal{Q}}$ with Landau coefficients $a$, $b$, and $c$ with units $Nm^{-2}$. These parameters, in effect, select a particular uniform value of $S$ in the bulk; by convention, $a$ is chosen to be temperature dependent, $a = a_0(T-T_0)$, where $T_0$ is the temperature below which the isotropic phase is no longer stable. The fourth term involving $\nabla{\mathcal{Q}}$ represents elasticity, and here, we adopt the commonly used one-constant approximation with a single elastic constant $L_1$ with units $N$. 
On the boundary, $\partial\mathcal{M}$, the constant term with prefactor $\sigma$ represents the isotropic surface tension. The second term is the anisotropic surface tension with anchoring coefficient $W$ and is constructed to favor the alignment of the LC in the tangent plane of the boundary surface. The surface parameters have units $Nm^{-1}$. Here, $\mathcal{Q}_s$ is a $Q-$tensor with alignment in the normal direction of the surface,
\begin{align}
    {\mathcal{Q}}_s = S_0\left(\bm{\nu}\otimes\bm{\nu} - \frac{\mathcal{I}}{d}\right),
\end{align}
where $\bm{\nu}$ is the local normal vector at the surface, $S_0$ is the degree of order induced by the surface, which may differ from the value set by the Landau coefficients in the bulk depending on the chemistry of the liquid crystal-host interface.  The prefactor $W>0$ controls the orientation such that, combined with the negative sign in front, higher values penalize the director's orientation toward the tangent plane.
As discussed above, we impose a volume constraint, 
\begin{align}\label{eq:globalVolumeConstraint}
    \mathcal{C}(\vec{x}) = \int_{\mathcal{M}} d\vec{x} - V_0 = 0, 
\end{align}
where, $V_0$ is the target volume of the tactoid. 
 
 In three dimensions, $\mathcal{Q}$ can be parameterized to be symmetric and traceless, 
 \begin{align}\label{generalQ3D}
		\mathcal{Q} = \begin{bmatrix}
			q_{xx} & q_{xy} & q_{xz}\\
			q_{xy} & q_{yy} & q_{yz}\\
			q_{xz} & q_{yz} & -q_{xx}-q_{yy}
		\end{bmatrix},
	\end{align}
and hence includes five independent degrees of freedom $\{q_{xx}, q_{xy}, q_{xz}, q_{yy}, q_{yz}\}$. 
We nondimensionalize $\mathcal{F}$ and $\mathcal{C}$ by introducing a length scale $\xi$, so that $\vec{x}\to \xi \bar{\vec{x}}$ where $\bar{\vec{x}}$ is non-dimensional, and divide \eqref{QEnergy} by $L_1$ to get,
 \begin{align}\label{eq:LandauEnergyNonDim}
		\overline{\mathcal{F}}(\overline{\vec{x}}, {\mathcal{Q}}) &= \int_{\mathcal{M}} \overline{a}\tr{({\mathcal{Q}}^2)} + \frac{2\overline{b}}{3}\tr{({\mathcal{Q}}^3)} + \frac{\overline{c}}{2}\tr{({\mathcal{Q}}^2)}^2 + \frac{1}{2}|\nabla{{\mathcal{Q}}}|^2\ d \overline{\vec{x}} \nonumber \\
  &+ \int_{\partial\mathcal{M}} \overline{\tau} - \frac{\overline{\tau}\overline{\omega}}{2}\tr{\big(({\mathcal{Q}}-{\mathcal{Q}}_s)^2\big)}\ d \overline{S},
\end{align}
with dimensionless parameters,
\begin{align}\label{eq:nondimparams}
    \quad
    \overline{a}= \frac{a\xi^2}{L_1},\quad \overline{b}= \frac{b\xi^2}{L_1},\quad \overline{c}= \frac{c\xi^2}{L_1},\quad \overline{\tau} = \frac{\sigma\xi}{L_1},\quad \overline{\omega} = \frac{W}{\sigma}.
\end{align}
The volume constraint, \eqref{eq:globalVolumeConstraint}, is trivially non-dimensionalized, and hence we simply choose $V_0$ to be dimensionless. Note that $\mathcal{Q}$ and $S$ are also non-dimensional.

Hence, to find the equilibrium order and mesh configuration, we minimize \eqref{eq:LandauEnergyNonDim} subject to the global nonlinear constraint, \eqref{eq:globalVolumeConstraint}. For the rest of the paper, we assume all quantities are suitably nondimensionalized and drop the bar notation. 

\subsection{Discrete Energy Model}
To set up the discrete optimization problem, we represent $\mathcal{M}\subset\mathbb{R}^d$ as a simplicial complex, $M$, consisting of $p$ spatial coordinate points, $\vec{x}_i\in M$, $i=1\hdots p$, and $N$ triangular ($d=2$) or tetrahedral ($d=3$) elements. We denote by $\vec{X}\in\mathbb{R}^{pd}$ the collection of all $\vec{x}_i\in M$.  We then use finite elements to discretize the Q-tensor, \eqref{generalQ3D}, over $M$, representing each of the five components as a piecewise polynomial function over the elements of $M$.  Note that this enforces the symmetry and tracelessness properties of the discrete Q-tensor, which we denote as $Q$. 

Finite-element approaches for Q-tensor models have been considered in \cite{Xia2023,Nestler2019, Farrell2023}, as well as for other LC frameworks \cite{Diegel2019, Borthagaray2020, Hirsch2023}. The algorithms we develop in this paper are independent of the specific finite-element spaces chosen as long as the systems remain well-posed.  However, we choose a linear piecewise-polynomial representation for simplicity, noting 
that higher-order representations are possible when needed. Using such linear finite elements, the degrees of freedom for the Q-tensor are defined as $Q_i:=\mathcal{Q}(\vec{x}_i)$ at each vertex point $\vec{x}_i\in M$.  Thus, in three dimensions, this discretization leads to $8p$ degrees of freedom: three for each coordinate of $\vec{x}_i$, and five for the $Q_i$, one for each $\{q_{xx}, q_{xy}, q_{xz}, q_{yy}, q_{yz}\}$ at $\vec{x}_i$.

Next, $\mathcal{F}(\vec{x}, \mathcal{Q})$ and $\mathcal{C}(\vec{x})$ are defined over $M$ to obtain $F(\vec{X}, Q)$ and $C(\vec{X})$,
\begin{align}
\label{eq:discreteLandauEnergyNonDim}
F(\vec{X}, Q) &= \int_{M} a\tr{({Q}^2)} + \frac{2b}{3}\tr{({Q}^3)} + \frac{c}{2}\tr{({Q}^2)}^2 + \frac{1}{2}|\nabla{{Q}}|^2\ d\vec{x} \nonumber 
    \\
&+ \int_{\partial M} \tau - \frac{\tau\omega}{2}\tr{\big(({{Q}}-{{Q}}_s)}^2\big)\ dS, \\
    \label{eq:discreteVolume}
    C(\vec{X}) &= \sum_{i = 1}^{N}\text{Vol}(T_i)-V_0,
\end{align}
where $\text{Vol}(T_i)$ is the volume of the $i^{th}$ tetrahedron, $T$, in the simplicial complex. The corresponding discrete shape optimization problem is then defined as,  
\begin{align}\label{eq:discretetactoidOptProblem}
    (\vec{X}^*,Q^*) &= \argmin_{\vec{X}, Q} \ F(\vec{X},Q),\\
    \text{s.t.} \quad 0 &= C(\vec{X}),  \nonumber 
\end{align} 
where again $\vec{X}\in M$ and $Q$ is a symmetric and traceless 2nd-order tensor whose components are piecewise linear functions defined on $M$.  Here, $\vec{X}^*$ and $Q^*$ represent the equilibrium configurations of the mesh points in the manifold and the Q-tensor approximation on $M$ for the minimized configuration.


\section{Energy Minimization of the Discrete Tactoid Energy}\label{discreteTactoidSection}
Considering the discrete $(\vec{X}, Q)$ shape optimization problem, \eqref{eq:discretetactoidOptProblem}, we now describe both a GD-based and a Newton-based method. We note that the GD-based approach in \cite{DeBenedictis2016} was formulated for a director formulation with a Frank--Oseen energy \cite{Frank1958}; we, therefore, describe in the first subsection how the method is modified for the Q-tensor formulation used in this paper. Next, we discuss the quasi-Newton approach, which is the main target of this work, and nested iteration, which is used to improve the performance of both methods.

\subsection{Gradient Descent}\label{PGDDescription}
In \cite{DeBenedictis2016}, the authors note that the combined $(\vec{X},\vec{n})$ optimization problem exhibited stiffness due to intrinsic differences in length scale between the vertices, $\vec{X}$, and directors, $\vec{n}$, degree of freedom. To alleviate this, they use an alternating gradient descent scheme, first taking descent steps in $\vec{n}$ followed by $\vec{X}$ and repeating until convergence. Given the connection between $\vec{n}$ and $Q$, we use a similar approach.

Consider a vector of the spatial vertices at the $k^{th}$ iteration, $\vec{X}^k\in M^k$, where $M^k$ is the current triangulation of the grid, $M$. We first compute the gradients of $F$ and $C$ with respect to $\vec{X}$ on the given grid $M$, denoted by $F^k_{M},$ and $C^k_M,$ respectively, and evaluate each gradient at $\vec{x}_i \in M^k.$ This trio of data, $\{\vec{X}^k, F^k_{M}, C^k_{M}\}$ is then projected onto the constraint's tangent space, 
\begin{align}\label{xHalfUpdate}
    \widetilde{\vec{X}}^{k} \leftarrow \vec{X}^k + \alpha_k\bigg(F^k_{M} - \frac{F^k_{M} \cdot C^k_M}{C^k_M \cdot C^k_M}C^k_M \bigg).
\end{align}
Here, $\alpha_k$ is the $k^{th}-$iteration step-size found from performing a one-dimensional line search. The update defined in \eqref{xHalfUpdate} is an intermediate step towards finding $\vec{X}^{k+1}\in M^{k+1}$ since $\widetilde{\vec{X}}^{k}$ only satisfies the constraint to linear order. Following  \cite{DeBenedictis2016}, we perform additional reprojection steps,
\begin{align}\label{xFullUpdate}
    \widetilde{\vec{X}}^k \leftarrow \widetilde{\vec{X}}^k + \frac{C(\widetilde{\vec{X}}^k)}{C^k_M \cdot C^k_M}C^k_M.
\end{align} 
In each reprojection step, $C^k_M$ is evaluated at $\widetilde{\vec{X}}_i \in \widetilde{M}^k.$ The prefactor, $\frac{C(\widetilde{\vec{X}}^k)}{C^k_M \cdot C^k_M},$ represents the difference between the constraint value and its target; reprojection steps are repeated until the constraint is satisfied to a given tolerance. This projected value is used as the next iterate $\vec{X}^{k+1} \leftarrow \widetilde{\vec{X}}^k$. 

To find $Q$ we follow a similar procedure but with $\{Q^k, F^k_Q\}$,
\begin{align}\label{QUpdate}
    Q^{k+1} \leftarrow Q^k + \beta_kF^k_Q,
\end{align}
where $\beta_k$ is obtained with a separate one-dimensional line search. A second projection is unnecessary as the global constraint only depends on $\vec{X}\in M.$ In practice,  \cite{DeBenedictis2016} starts by optimizing for the field and then for the spatial coordinates. To control mesh quality, a procedure called \textit{equiangulation} is used after a few iterations of each optimization routine.  This ensures that no irregular triangles (i.e., long and skinny) appear in the discrete manifold, and is crucial when the amount of spatial or order deformation increases.

	
    \subsection{Quasi-Newton}\label{QN}
	The GD algorithm introduced in the previous section converges quite slowly and can stagnate under certain conditions. It also incorporates a number of metaparameters, such as the mesh quality control and the number of alternating iterations to be taken for shape and field degrees of freedom that must be hand-tuned for each problem. To address these issues, we develop a quasi-Newton (QN) based method to solve the full $(\vec{X},Q)$ minimization problem all at once. Given the presence of the nonlinear constraint, we introduce the Lagrangian $\mathcal{L}$ of the system based on \eqref{eq:discretetactoidOptProblem},
	\begin{align}\label{eq:tactoidLagrangian}
		\mathcal{L}(\vec{X},Q,\lambda) = F(\vec{X},Q) - \lambda C(\vec{X}),
	\end{align}
	where $\lambda \in \mathbb{R}$. The necessary first-order optimality conditions for minimizing the Lagrangian are 
	\begin{align}
		\vec{\mathcal{L}}_{M} &:= F_{M} - \lambda C_{M} = \vec{0}, \label{eq:LM} \\
		\vec{\mathcal{L}}_{Q} &:= F_{Q} = \vec{0}, \label{eq:LQ} \\
		\mathcal{L}_{\lambda} &:= -C(\vec{X})  = 0, \label{eq:Llambda}
	\end{align}
	where $F_{M}$ and $F_{Q}$ are the first-order G$\hat{\text{a}}$teaux derivatives with respect to every $\vec{x}_i\in M$ and $Q_i$ on $M$, respectively (i.e., the gradients as computed for GD), and $C_{M}$ is the gradient of the constraint with respect to each $\vec{x}_i$.
	
	These equations are nonlinear, so we linearize \eqref{eq:LM}-\eqref{eq:Llambda} and set up the corresponding iterative scheme. Let $(\vec{X}^k, Q^k,\lambda^k)$ be the current approximations for $(\vec{X}, Q,\lambda)$, respectively. The update, $(\vec{d}_{M}, \vec{d}_{Q}, d_{\lambda})$, to the approximation is the solution to the QN iteration,
	\begin{align}\label{eq:QNIterations}
		\begin{bmatrix}
			\vec{d}_{M,Q}\\
			d_{\lambda}
		\end{bmatrix} = \begin{bmatrix}
			\vec{B}^k & -(\vec{A}^k)^T \\
			-{{\vec{A}^k}} & 0
		\end{bmatrix}^{-1}
		\begin{bmatrix}
			{\vec{R}^k}\\
			-\mathcal{L}_{\lambda}
		\end{bmatrix} =: \mathcal{A}^{-1}\mathcal{R},
	\end{align} 
	with
	\begin{align*}
		\vec{d}_{M,Q} &= \begin{bmatrix}
			\vec{d}_{M}\\ \vec{d}_{Q}\end{bmatrix}, &
		\vec{A}^k &= \begin{bmatrix}
			C_{M}&\vec{0}\\ \end{bmatrix}, &
		\vec{R}^k &= \begin{bmatrix}
			-\vec{\mathcal{L}}_{M}\\ -\vec{\mathcal{L}}_{Q}\end{bmatrix}.
	\end{align*}
	Here, $\vec{B}^k$ is an approximation to the Hessian of the $(\vec{X},Q)$ portion of the Lagrangian, denoted by $\nabla^2\mc{L}^k := \begin{bmatrix}
                    \mc{L}_{MM}^k & \mc{L}_{MQ}^k \\
                    \mc{L}_{QM}^k & \mc{L}_{QQ}^k
                \end{bmatrix} $, where the entries are second-order G$\hat{\text{a}}$teaux derivatives with respect to the mesh and the $Q-$tensor at $(\vec{X}^k,Q^k)$. 
 Additionally, $\vec{A}^k$ contains the gradient of the constraint function at $\vec{X}^k$, and  $\vec{R}^k$ is the nonlinear residual for the $(\vec{X},Q)$ portion of the system. 
 Since $\mathcal{A}$ is a saddle-point matrix, which poses challenges for building efficient solvers, we choose a BFGS approximation and use explicit formulae to compute $\vec{H}^k = {(\vec{B}^k)}^{-1}$ and $\mathcal{A}^{-1}$ \cite{Benzi2005}. This allows us to find a closed form for computing the updates $(\vec{d}_{M}, \vec{d}_{Q}, d_{\lambda})$ which we use to set up a recursive matrix-free approach of implicitly updating the approximation of $\vec{H}^k.$ 

	Next, we find the step size for the field, $\alpha_Q$, and the step size for the spatial coordinates $\alpha_M.$ The step size $\alpha_Q$ is found with backtracking line search such that it satisfies,
	\begin{align}
		F(\vec{X}^k,Q^{k} + \alpha_Q\vec{d}_{Q}) &\leq F(\vec{X}^k,Q^k)+ \eta\alpha_Q(F^k_Q)^T\vec{d}_{Q},
	\end{align}
	where $\eta \in (0,1).$ 
	Then, we update the field accordingly,
	\begin{align}
		Q^{k+1} &\leftarrow Q^k +\alpha_Q\vec{d}_Q.
	\end{align}
	The step size $\alpha_M$ for the spatial coordinates is found using an $\ell^1$ monitor function \cite{Bach2011} defined as 
	\begin{align}
		\phi(\vec{X},Q^k,\mu) &:= F(\vec{X},Q^k) + \mu\|C(\vec{X})\|_1,
	\end{align} where $\mu \in \mathbb{R}$ is a penalty term for the constraint. This monitor function decides what factor of $\vec{d}_{M}$ and $d_\lambda$ should be accepted in order to decrease $F(\vec{X}^{k+1},Q^{k})$. 
	Similar to finding $\alpha_Q$, we use backtracking until the following inequality is satisfied,
	\begin{align}\label{tactoidarmijo}
		\phi(\vec{X}^k + \alpha_M\vec{d}_{M},Q^{k}, \mu_k) &\leq \phi(\vec{X}^k,Q^k,\mu_k) + \eta\alpha_M(\vec{\phi}_M^k)^T\vec{d}_{M},\nonumber
	\end{align}
	where $\eta \in (0,1)$, $\mu_k > \|\lambda_k\|_{\infty}+\rho$ ($\rho$ is a positive constant),  and 
	\begin{align}
		(\vec{\phi}_M^k)^T\vec{d}_{M} &:= \alpha_MF^k_M\vec{d}_{M} - \mu_k\|C(\vec{X}^k)\|_1. 
	\end{align}
	Finally, we update the spatial coordinates and Lagrange multiplier,
	\begin{align}
		\vec{X}^{k+1} &\leftarrow \vec{X}^k +\alpha_M\vec{d}_M\\
		\lambda^{k+1} &\leftarrow \lambda^k +\alpha_Md_{\lambda}.
	\end{align}
    Given that the constraint is prescribed only to the mesh, we allow for its respective Lagrange multiplier, $\lambda^{k}$, to have the same step size $\alpha_M$ as the mesh update. 
    
    Finally, we note that quasi-Newton methods have been extensively studied, and several convergence results have been established. For the BFGS approach we consider, a well-known result \cite[Thm. 18.6]{Nocedal1999} guarantees superlinear convergence of the method under certain assumptions.  Though this result does not consider line search, for all of the tests presented here, the step sizes approach one as the method converges.  The assumptions of the theorem include conditions on the continuity and differentiability of the objective and constraint function that can be shown by direct calculation.  Moreover, there are conditions on the linear independence of the gradients of the constraints in order to guarantee that the first-order optimality conditions, \eqref{eq:LM}-\eqref{eq:Llambda}, are satisfied.  Since we only have one constraint here, this is trivially true.  We also need to assume that the Hessian, $\nabla^2\mc{L}$, at the local solution, $(\vec{X}^*,Q^*)$ and its initial approximation, $\vec{B}^0$, are symmetric and positive definite.  The latter is straightforward since $\vec{B}^0 = \delta\mathcal{I}$ with $\delta>0$, while the former requires several assumptions on the physical parameters.  The numerical experiments we present indicate that we are within the regime where this is satisfied.  Finally, as with all Newton-based methods, we need the initial guess to be ``good enough" in order for the iterates to remain in the basin of attraction.  To guarantee this, we introduce the notion of nested iteration.

	\subsection{Nested Iteration}\label{NI}
	Both QN and GD, as described above, are sensitive to initial guesses. In particular, as the anisotropic surface tension parameters increase, we expect to see the nematic LC tactoids elongating towards the characteristic eye shape \cite{Bates2003, Prinsen2003}. 
	Numerically, this means that by starting from the same initial guess for every parameter value, we may not be close to the basin of attraction for some regions of parameter space. This can lead the method to stagnate at locally optimal solutions, not the expected global minima. To remedy this issue, we wrap each method with nested iteration. 
	
	Nested iteration \cite{Starke2000, Briggs2000} begins with an initial coarse grid, denoted by $M_i$ with only a few vertex points, $p$, that represents the problem domain. This coarse grid may not capture all the details of the solution, but solving the nonlinear system on this grid is computationally inexpensive. Thus, \eqref{eq:discretetactoidOptProblem} is solved on the coarse grid with either QN or GD until a preferred convergence criterion is reached. The coarse grid, $M_i$, is then subdivided into smaller elements. In this work, we consider uniform refinement such that each element is divided into four (for triangular elements) or eight (for tetrahedral elements) smaller elements. 
	The coarse solution from $M_i$ is then linearly interpolated onto the finer grid $M_{i+1}$ and used as an initial guess for solving the problem now represented on $M_{i+1}$ of size $8p$ (in three-dimensions) and $4p$ (in two-dimensions). This initial solution on $M_{i+1}$ should then be a more accurate guess as it came from solving a similar problem on $M_i$, and thus fewer iterations to converge are expected.
 \begin{remark}
     Through the use of nested iteration, we progressively improve the initial guess, $(\vec{X}^0,Q^0)$, thereby maintaining that $\|(\vec{X}^0,Q^0)-(\vec{X}^*,\vec{Q}^*)\|$ is sufficiently small. Consequently, this helps guarantee convergence results, such as those in \cite{Nocedal1999}.
 \end{remark}

\section{Numerical Results}\label{results}

We demonstrate the robustness of QN combined with Nested Iteration (NI) on challenging two-dimensional and three-dimensional problems involving the formation of a nematic tactoid. We compare the approach with the GD-based techniques described in Section \ref{PGDDescription}. 
Numerical methods and benchmark tests are implemented and executed in \emph{Morpho}, an open-source programmable environment for shape optimization \cite{joshi2024}. \emph{Morpho} is able to evaluate the objective function of interest as well as its gradients with respect to the configuration's and field's degrees of freedom. Furthermore, grid quality control in \emph{Morpho} does not require user intervention and provides the user the option of automatic domain refinement and easy object-oriented programming. All timed numerical results are done using a workstation with an $8$-core $3$-GHz Intel Xeon Sandy Bridge CPU and $256$ GB of RAM. Force and energy calculations in \emph{Morpho} are parallelized using a symmetric multiprocessing model with a user-controllable number of worker threads. In the numerical experiments we report below, we use 16 worker threads.
For all timings reported, we do not include the time spent during refinement and equiangulation \textit{between} nested iteration grids, because these components are not yet parallelized and are common to all methods. 

Throughout this section, all test problems  use material constants from \cite{Mottram2014} with $a_0 = 0.042\times10^{6}~Nm^{-2}K^{-1}$ and $T-T_0 = -0.1~K$,  Thus, the prefactors in the \emph{dimensional} free energy, \eqref{QEnergy}, are: 
\begin{align}
    a&=-0.042\times10^{5}~Nm^{-2},&
    b&=-0.64\times10^{6}~Nm^{-2},\\ c&=0.35\times10^{6}~Nm^{-2},& 
    L_{1}&=1\times 10^{-11}~N.
\end{align}
Finally, we note that in all two-dimensional experiments the resulting equilibrium configurations are found to be uniaxial, such that the spectrum of the Q-tensor lies in the interval $[-1/3,2/3]$.  However, some of the three-dimensional results yield biaxial configurations and, as such, the spectrum deviates slightly from this interval.

\subsection{Two-Dimensional Nematic Tactoids}\label{2dSection}

We begin by modeling two-dimensional nematic tactoids, where the computational domain represents a slice of an object that is infinitely extended in the perpendicular direction. We investigate solutions where the nematic director lies in-plane on the computational domain, and hence $\mathcal{Q}$ can be described by a reduced parameterization,
\begin{align}\label{generalQ}
    \mathcal{Q} = \begin{bmatrix}
        q_{xx} & q_{xy} & 0\\
        q_{xy} & \frac{1}{3}-q_{xx} & 0\\
        0 & 0 & -\frac{1}{3}
    \end{bmatrix},
\end{align}
that includes only two independent degrees of freedom $q_{xx}$ and $q_{xy}$. 
We construct the preferred surface value of the Q-tensor ${\mathcal{Q}}_s$ in \eqref{eq:LandauEnergyNonDim} from the local tangent vector $\vec{t}$ on the boundary,
\begin{align}
    {\mathcal{Q}}_s = S_0\left(\vec{t}\otimes\vec{t} - \frac{\mathcal{I}}{3}\right).
\end{align}
To favor alignment with the local tangent vector, we must prefactor the anchoring term in \eqref{eq:discreteLandauEnergyNonDim} with a positive sign,
\begin{equation*}\int_{\partial M} \tau + \frac{\tau\omega}{2}\tr{\big(({{Q}}-{{Q}}_s)}^2\big)\ dS.
\end{equation*}
The volume constraint becomes a cross-sectional area constraint for a target area $A_0$, 
\begin{align}
    \mathcal{C}(\vec{x}) = \int_{\mathcal{M}} d\vec{x} - A_0 = 0. \nonumber
\end{align}
Consequently, the constraint has the explicit discrete form,
\begin{align}
    C(\vec{X}) &= \sum_{i = 1}^{p-2}\frac{1}{2}\|(\vec{x}_{i+1}-\vec{x}_{i}) \times (\vec{x}_{i+2}-\vec{x}_{i+1})\| - A_0. \nonumber
\end{align}
We use the same parameters as those listed in \eqref{eq:nondimparams} for the two-dimensional geometry. In this two-dimensional setting, we envision our domain as a $2D$ slab embedded in a three-dimensional space. Therefore, the energy functional \eqref{eq:LandauEnergyNonDim} now has dimensions of energy per unit length or $Jm^{-1} = N$.
Setting $\xi = 1\times10^{-7}~m$, then, the  constants corresponding to those in the \emph{non-dimensionalized} free energy, \eqref{eq:discreteLandauEnergyNonDim}, are 
\begin{align}\label{params}
    \overline{a}&=-4.2,&
    \overline{b}&=-640,& \overline{c}&=350.
\end{align}
For the numerical experiments below, we let the surface tension, $\overline{\tau}$, vary from $1$ to $100$ and the surface anchoring, $\overline{\omega}$ vary from $0.01$ to $1.$ This corresponds to a dimensionalized surface tension, $\sigma$, ranging from $10^{-4}-10^{-2}\,N$ and anchoring values, $W$, ranging from $10^{-6}-10^{-2}\,N$. Moreover, we artificially scale the Landau coefficients so that the defect size is $\sim10$ times its true value. The scalar order parameter $S_0$ for the tangential anchoring is initialized as, 
\[S_0 = \frac{-\overline{b}+\sqrt{\overline{b}^2-24\overline{a}~\overline{c}}}{4\overline{c}} = 0.933567.\] For the rest of the section, we refer to nondimensionalized quantities only and drop the bar notation. 

For most test problems, the initial guess on the coarsest grid is defined on a circle of area $1$ with equally distributed vertices representing the degrees of freedom, $\vec{X}\in M$. The initial director field is aligned parallel to the $x-$axis.
We start by solving the discrete problem on the coarsest grid of size $|M_i| = 16$, where the notation  $|M_i|$ denotes the number of vertices on grid $M_i$, and propagate the solution across nine levels of refinement until the finest grid of size $|M_i| = 591,361$. For both methods described in this work, after each level of refinement, we perform the equiangulation procedure mentioned above, that improves the quality of the mesh, ensuring there are no irregular elements. The preferred convergence criterion is the relative change in the energy $F^k$, with the tolerance set at $10^{-6}$. We compare the results from the NI algorithm to those obtained running solely on the finest grid, $M_9$. We compare the final energy, $F^k$, the runtime in seconds, and the number of iterations on each grid. 

\subsubsection{Fixed Shape and Fixed Field Subproblems}
Before performing the full $(\vec{X},Q)$ optimization, we aim to understand the separate role of shape and $Q$ degrees of freedom by defining two subproblems derived from \eqref{eq:discretetactoidOptProblem}:\\[1em]
\textbf{Subproblem A} (Fixed Shape)
    \begin{align}\label{eq:fsmf}
        {Q}^*=\argmin_{Q} F(\vec{X}^*,Q);
    \end{align}
\textbf{Subproblem B} (Fixed Field)
    \begin{align}\label{eq:msff}
        \vec{X}^* &= \argmin_{\vec{X}} \ F(\vec{X},Q^*), \\
        \text{s.t.} \quad 0 &= C(\vec{X}). \nonumber
    \end{align}
These subproblems are used to motivate the use of our QN-based method. Therefore, we do not compare it with the gradient descent method. In Subproblem A, \eqref{eq:fsmf}, which we call the Fixed Shape model, we assume the grid, $M^*$, with vertices $\vec{X}^*$ is fixed and at equilibrium. Similarly, for Subproblem B, \eqref{eq:msff}, which we call the Fixed Field model, the field, or Q-tensor, $Q^*$ defined on $M$, is fixed and at equilibrium. 
 
\begin{figure}[h!]
		\begin{subfigure}{0.49\textwidth}
			\centering
			\includegraphics[width=0.7\linewidth]{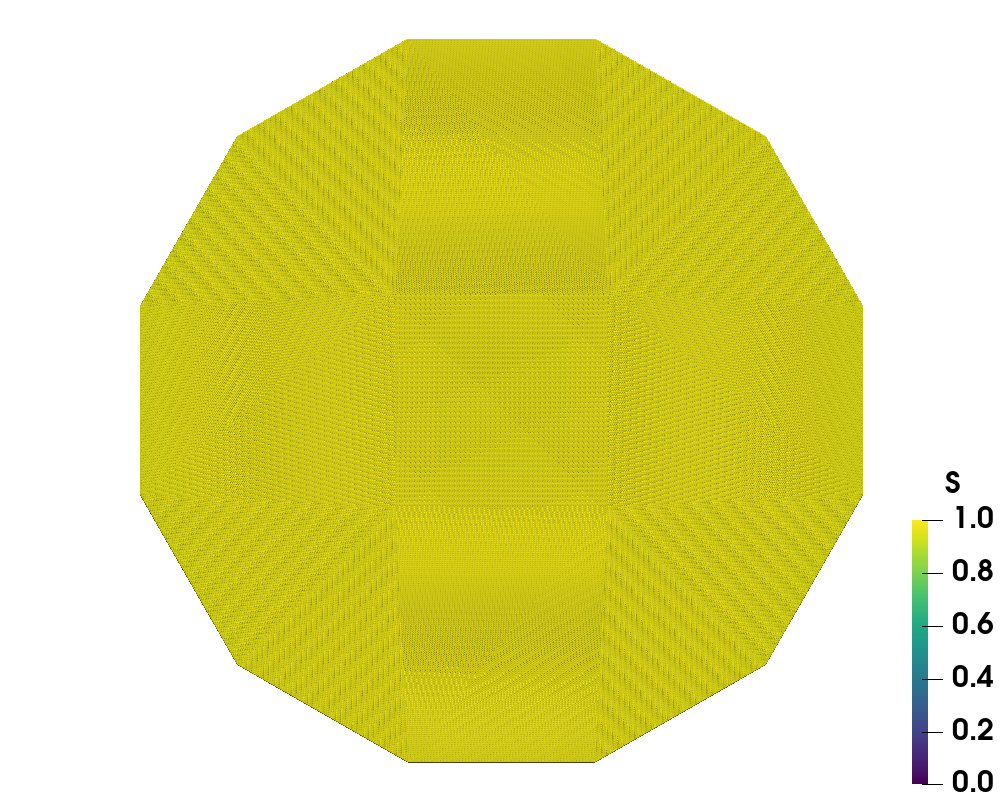}
			\caption*{(a) $\tau = 1$}
		\end{subfigure}%
		\begin{subfigure}{0.49\textwidth}
			\centering
			\includegraphics[width=0.7\linewidth]{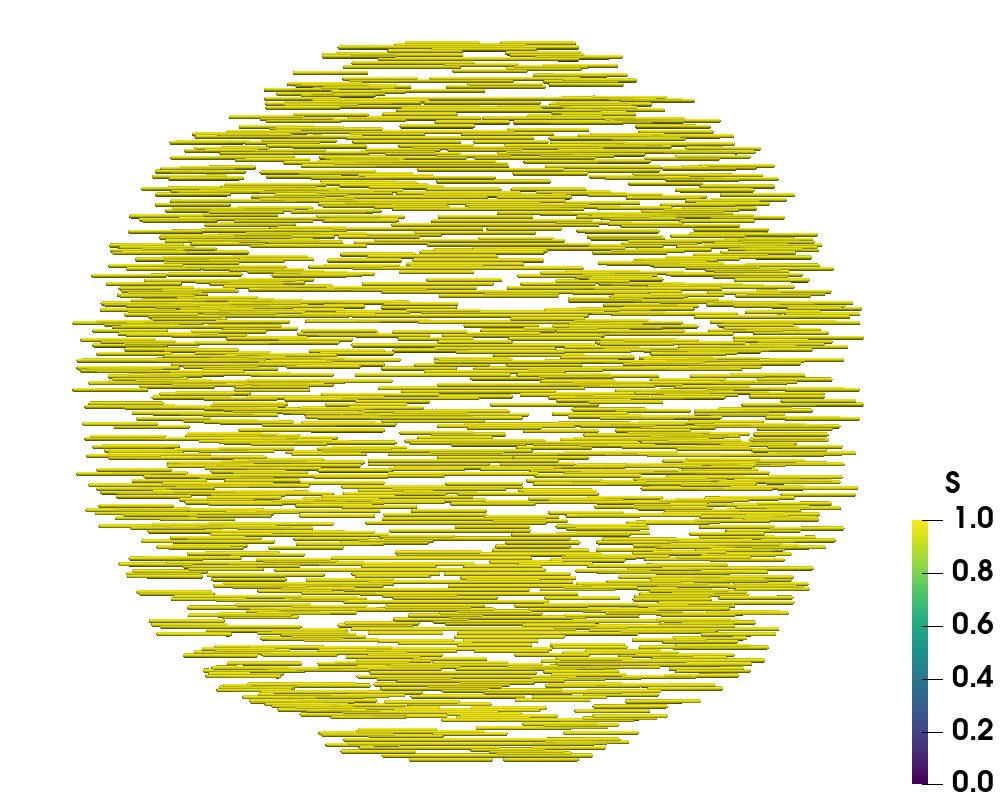}
			\caption*{(b)}
		\end{subfigure}\\
		\begin{subfigure}{0.49\textwidth}
			\centering
			\includegraphics[width=0.7\linewidth]{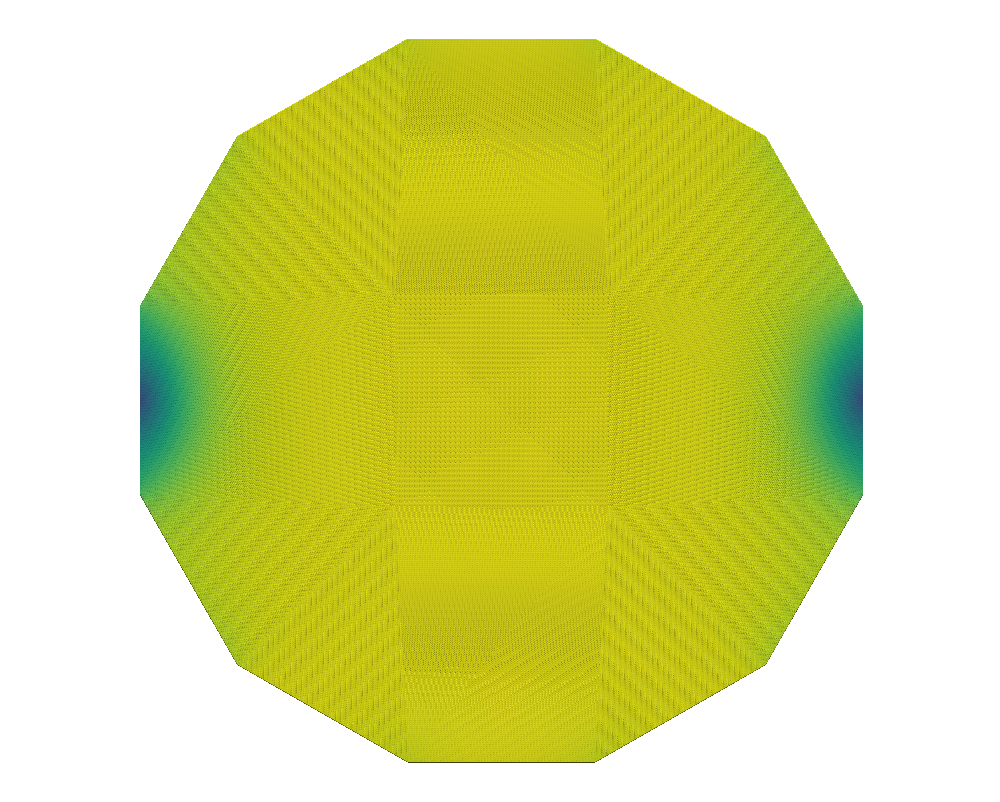}
			\caption*{(c) $\tau = 50$}
		\end{subfigure}%
		\begin{subfigure}{0.49\textwidth}
			\centering
			\includegraphics[width=0.7\linewidth]{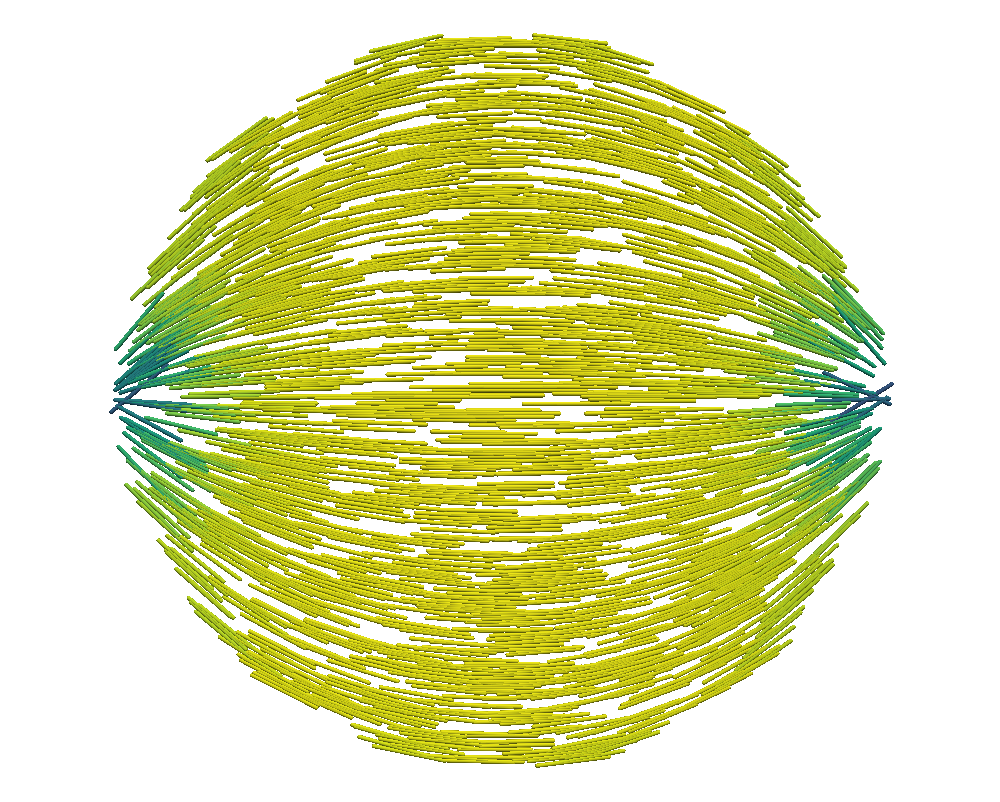}
			\caption*{(d)}
		\end{subfigure}\\
        \begin{subfigure}{0.49\textwidth}
			\centering
			\includegraphics[width=0.7\linewidth]{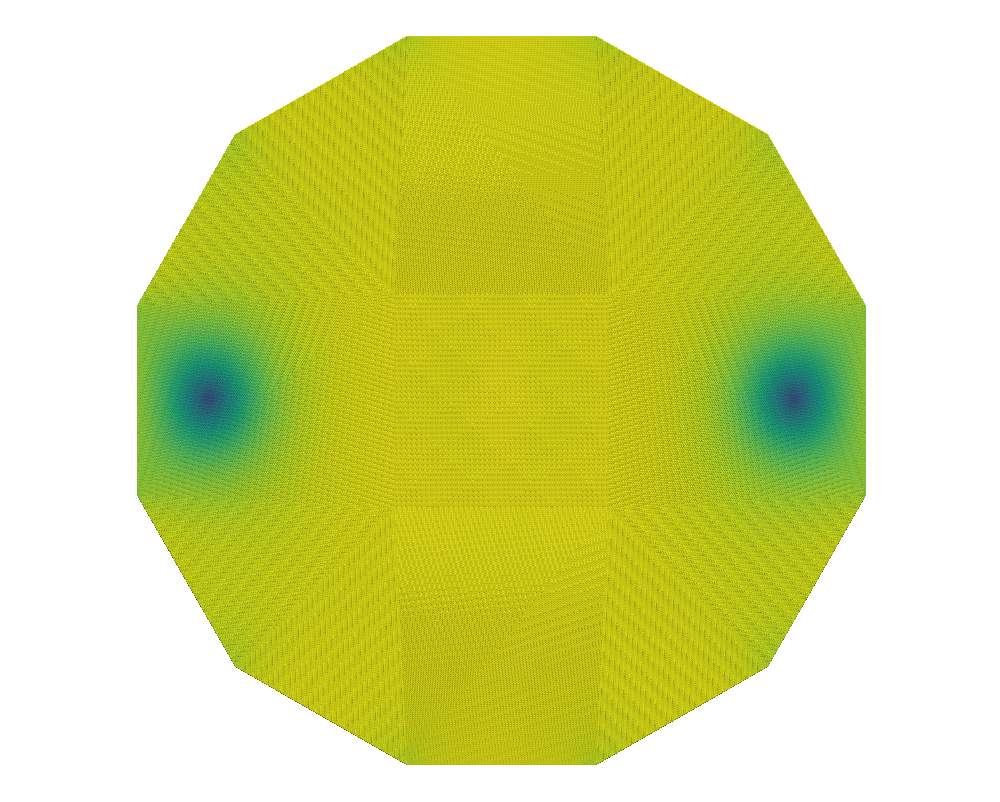}
			\caption*{(e) $\tau = 100$}
		\end{subfigure}%
		\begin{subfigure}{0.49\textwidth}
			\centering
			\includegraphics[width=0.7\linewidth]{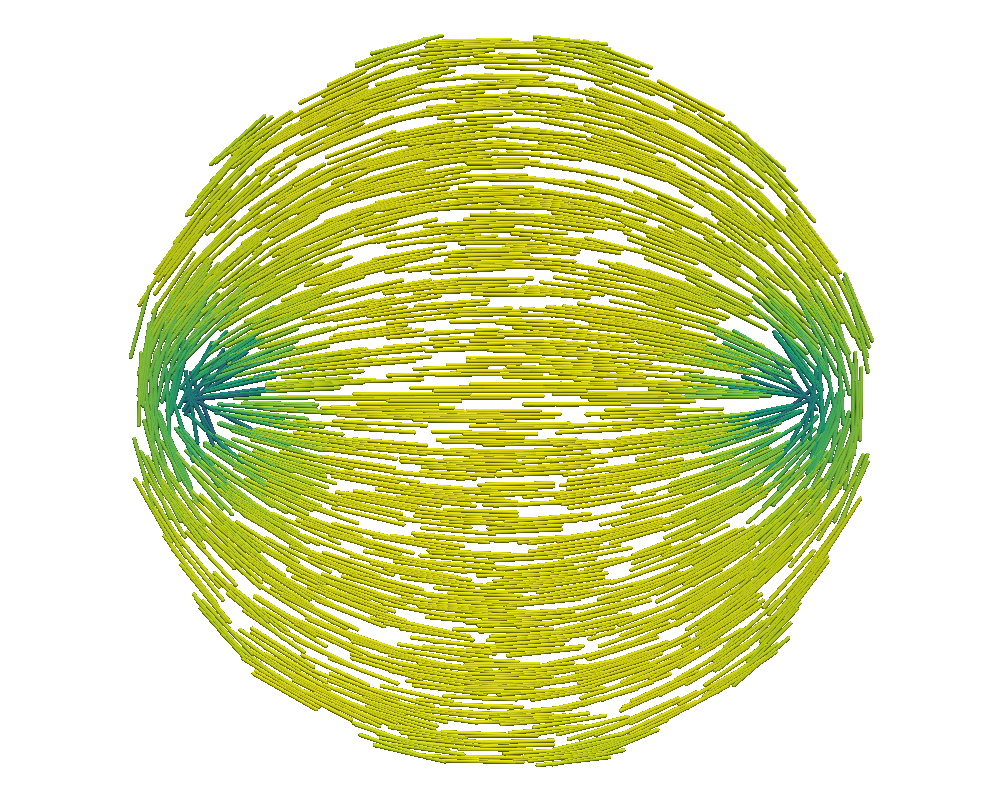}
			\caption*{(f)}
		\end{subfigure}
		\caption{\justifying Applying QN \textit{without} NI for Subproblem A (Fixed Shape). Results shown on grid $M_9$. The color bar indicates the value of $S$, i.e., the order of the director field in the domain. Left plots, (a), (c), (e), depict the order's distribution, illustrating two defects for higher $\tau$. Right plots, (b), (d), (f), show the directors' anchoring to the boundary as $\tau \rightarrow 100$ and also demonstrate the disorder around the defects.}
		\label{fig:QNfsmf}
	\end{figure}

\begin{figure}[h!]
		\begin{subfigure}{0.49\textwidth}
			\centering
			\includegraphics[width=0.7\linewidth]{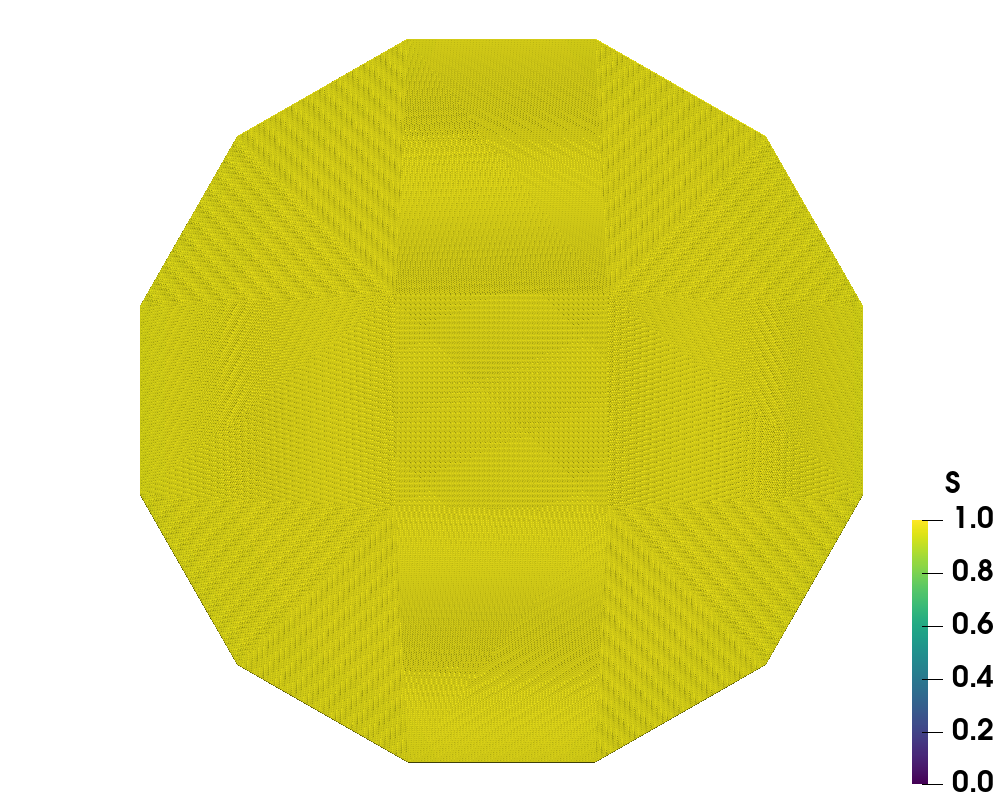}
			\caption*{(a) $\tau = 1$}
		\end{subfigure}%
		\begin{subfigure}{0.49\textwidth}
			\centering
			\includegraphics[width=0.7\linewidth]{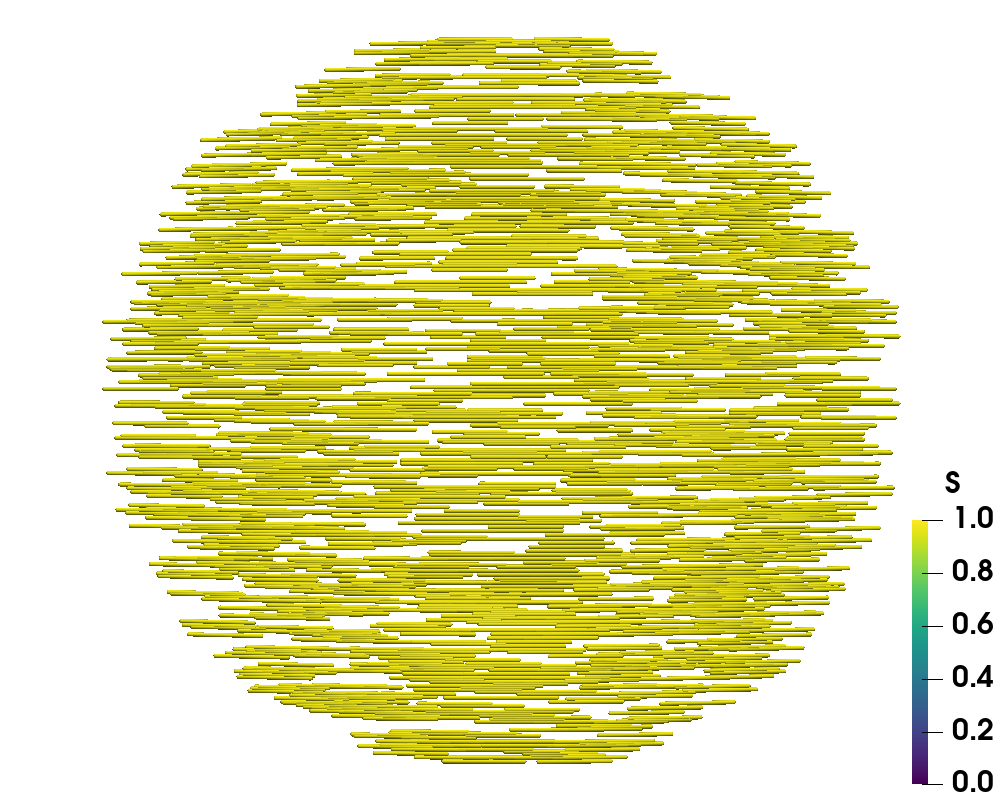}
			\caption*{(b)}
		\end{subfigure}\\
		\begin{subfigure}{0.49\textwidth}
			\centering
			\includegraphics[width=0.7\linewidth]{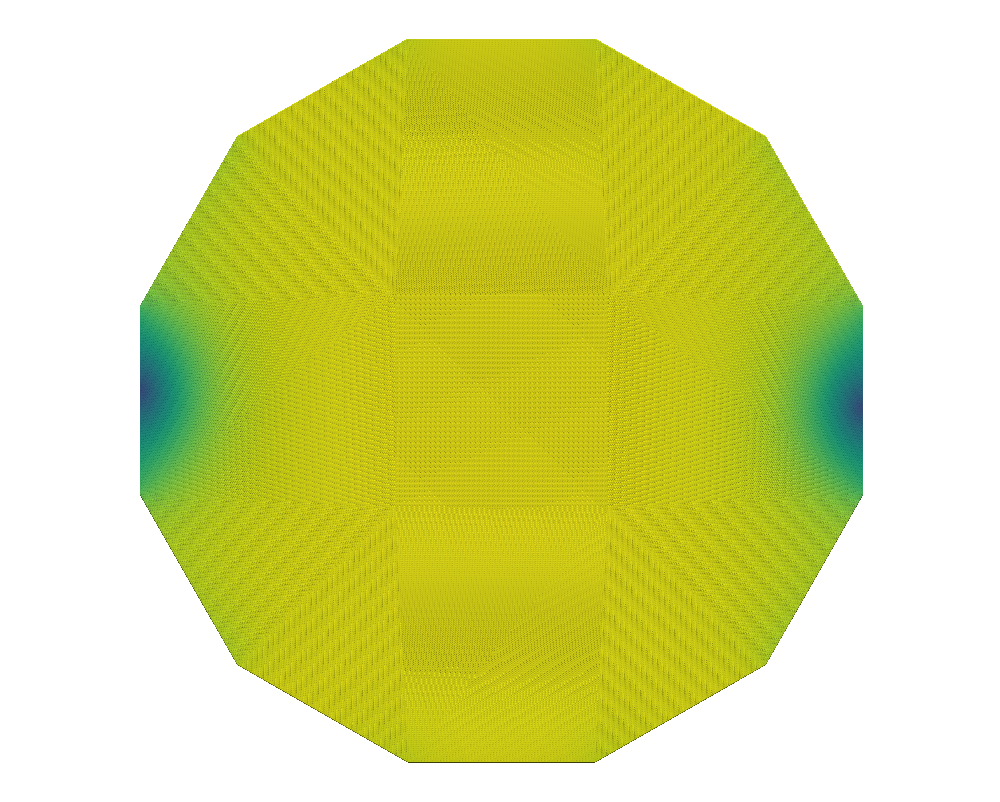}
			\caption*{(c) $\tau = 50$}
		\end{subfigure}%
		\begin{subfigure}{0.49\textwidth}
			\centering
			\includegraphics[width=0.7\linewidth]{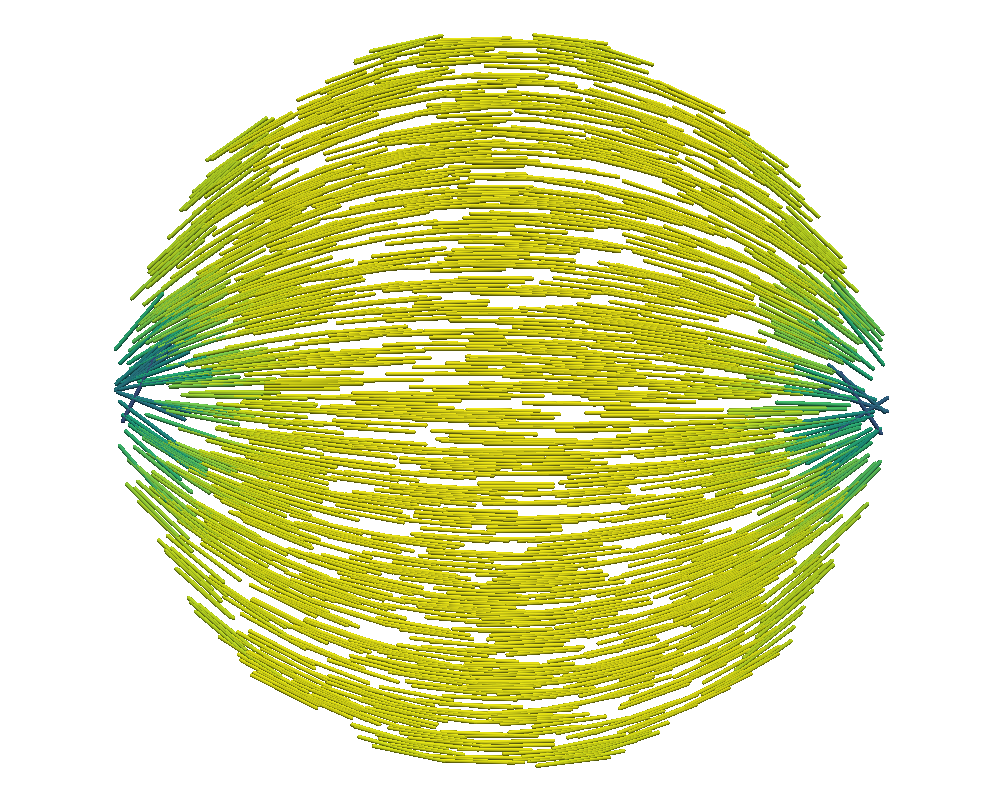}
			\caption*{(d)}
		\end{subfigure}\\
        \begin{subfigure}{0.49\textwidth}
			\centering
			\includegraphics[width=0.7\linewidth]{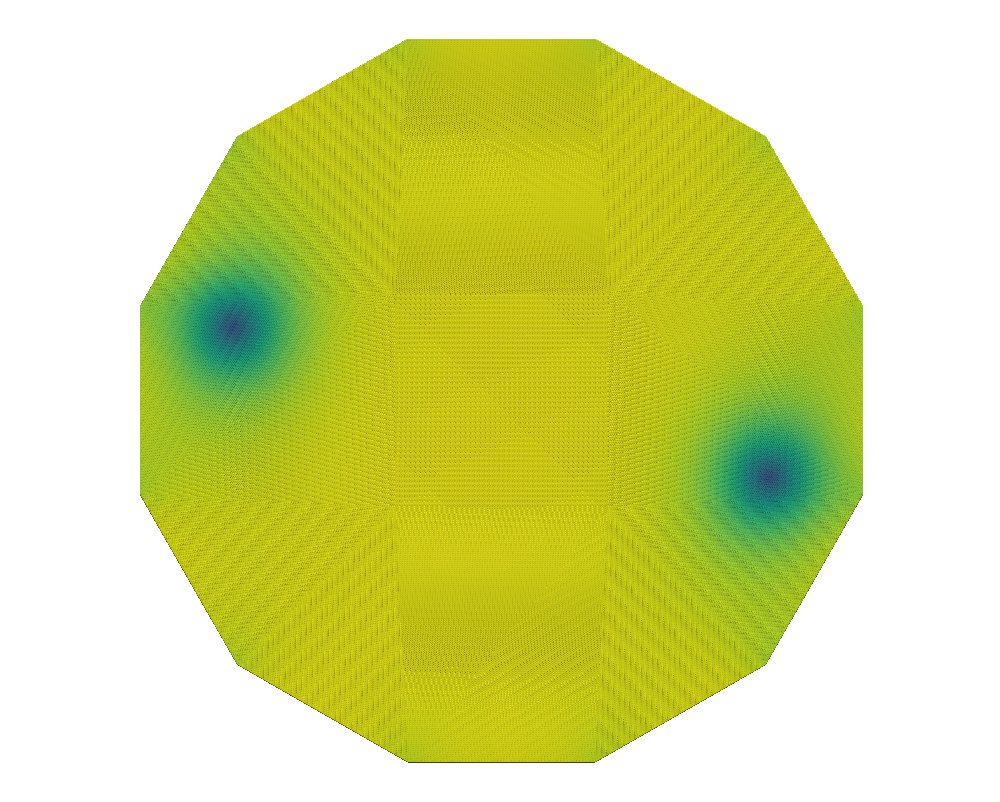}
			\caption*{(e) $\tau = 100$}
		\end{subfigure}%
		\begin{subfigure}{0.49\textwidth}
			\centering
			\includegraphics[width=0.7\linewidth]{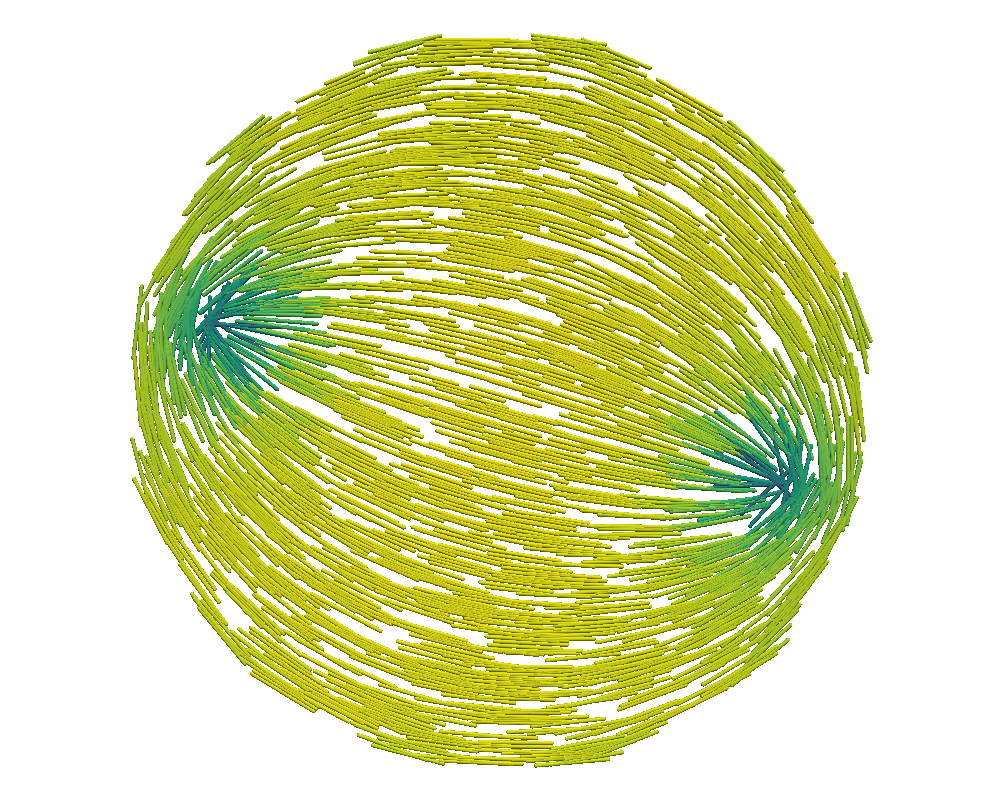}
			\caption*{(f)}
		\end{subfigure}
		\caption{\justifying Applying QN \textit{with} NI for Subproblem A (Fixed Shape). Results shown on grid $M_9$. The color bar indicates the value of $S$, i.e., the order of the director field in the domain. Left plots, (a), (c), (e), depict the order's distribution illustrating two defects for higher $\tau$. Right plots, (b), (d), (f), show the directors' anchoring to the boundary as $\tau \rightarrow 100$ and also demonstrate the disorder around the defects.}
		\label{fig:QNNIfsmf}
	\end{figure}

 Solutions to the Fixed Shape model are director fields with different alignments strongly affected by the isotropic surface tension $\tau$ while the anchoring parameter remains small at $\omega = 0.2$. The boundary terms from \eqref{eq:fsmf} can be split into two parts: the line tension integral,
	\begin{align}\label{eq:lineTensionIntegral}
		\tau\int_{\partial M}\, dl,
	\end{align}
	and the surface anchoring integral,
	\begin{align}\label{eq:surfaceAnchoringIntegral}
		\frac{\tau\omega}{2} \int_{\partial M} \tr{({Q}-{Q}_s)}^2\ dl.
	\end{align}
	With the grid fixed, as we increase $\tau$, we only see an effect from the surface anchoring integral, \eqref{eq:surfaceAnchoringIntegral}, as it depends on the Q-tensor. Here, the prefactor $\Gamma := \frac{\tau\omega}{2}$ promotes stronger tangential anchoring as $\tau$ increases. Following Bates' conclusion \cite{Bates2010}, we expect to see stronger tangential alignment on the domain's boundary as $\tau$ increases as well as the emergence of two defects \cite{Dzubiella2000, Sitta2018} inside the tactoid domain. This is confirmed in Figures \ref{fig:QNfsmf} and \ref{fig:QNNIfsmf}.

	We compare the results of QN with and without nested iteration for various values of $\tau$ on grids up to $M_9$. We retain the overall shape of the domain, a dodecagon, from the coarsest grid. The left graphics of Figures \ref{fig:QNfsmf} and \ref{fig:QNNIfsmf} show the distribution of the scalar order parameter $S$, and the corresponding director field is given on the right. Both standalone QN and QN with NI depict two defects appearing inside the tactoid domain with stronger surface anchoring, as expected. Both methods also demonstrate that the director field anchors tangentially to the boundary as $\tau$ tends to $100$. We note that while the solutions from nested iteration locate the defect pair near two opposing vertices of the polygonal boundary, we also see a rotational invariance in the energy. Nevertheless, the two methods have converged to similar energies. 

     \begin{table}[!ht]
                    		\centering
                    		\resizebox{\columnwidth}{!}{%
                    			\begin{tabular}{cc|ccccc}
                    				\hline
                    	     &   & $\tau = 1$  & $5$         & $10$               & $50$                & $100$      \\ \hline
                    	NI Grid& $|M_i|$ &\multicolumn{5}{c}{Iterations}\\ \hline
                        $M_1$& $16$   & $19$        & $26$        & $25$               & $37$                & $32$       \\ 
                        $M_2$& $49$  & $5$         & $17$        & $21$               & $22$                & $83$     \\ 
                    	$M_3$& $169$  & $3$         & $10$        & $10$                & $20$                & $35$   \\ 
                    	$M_4$& $625$  & $5$         & $10$        & $13$            & $16$                & $17$   \\ 
                    	$M_5$& $2,401$  & $2$         & $9$         & $10$             & $12$                & $16$   \\ 
                    	$M_6$& $9,409$  & $2$         & $8$         & $8$                & $11$                & $12$        \\
                        $M_7$& $37,249$ & $1$         & $1$         & $1$                & $8$                  & $10$        \\
                        $M_8$& $148,225$  & $1$         & $4$         & $5$                  & $6$                  & $8$        \\ 
                        $M_9$& $591,361$  & $1\ (35)$   & $1\ (79)$   & $1\ (549)$       & $3\ (666)$        & $4\ (851)$        \\ \hline
$F^k$ &  & $-19.87\ (-19.85)$ & $-18.73\ (-18.65)$ & $-17.53\ (-17.53)$  & $-12.99\ (-12.98)$  & $-12.04\ (-12.01)$ \\ \hline
\textbf{Runtime [sec]}& & $\vec{37.32\ (925.66)}$& $\vec{67.79\ (2,122.25)}$& $\vec{75.53\ (18,553.0)}$  & $\vec{154.34\ (22,792.3) }$ & $\vec{200.84\ (31,390.7)}$    \\ \hline
                    			\end{tabular}
                    		}
                \caption{\footnotesize \justifying Subproblem A (Fixed Shape): Iteration count for QN with NI on each grid level. Iteration count for QN without NI is given on level $M_9$ in parenthesis. The grid size, $|M_i|$, final energy, $F^k$, and runtime in seconds, in bold, for the full simulation of QN with NI (standalone QN in parenthesis) are also given. NI improves efficiency in terms of iteration and runtime by taking minimal steps on the finest grid.}
        		\label{tab:fsmfData}
    \end{table}
	
	Table \ref{tab:fsmfData} depicts the number of QN iterations both on grid $M_9$ alone (in parenthesis) and using NI to build to $M_9$, showing the iteration count on each level of refinement. For small values of $\tau$, the initial guess improves along the grids in the nested iteration process, thereby requiring fewer iterations to find the equilibrium solution on the finer grids. For increasing values of $\tau$, while each level's iterations increase, we still see it decreasing down the levels, as expected. The table also shows the converged energy, $F^k$, 
 and the runtime in seconds for the simulations. We see that both QN alone and with NI yield the very similar converged solutions that minimize the energy. However, the addition of nested iteration allows the method to converge to a slightly lower energy which is physically expected. 
NI significantly improves the timing compared to standalone QN by a factor of $24$ to $245$ times.
	
	Conversely, solutions to the Fixed Field model, \eqref{eq:msff}, are the various shapes that result at a fixed $\tau = 10$, but varying the anchoring strength, $\omega$. With the molecular alignment fixed across the shape, we see an effect from both line tension \eqref{eq:lineTensionIntegral} and anchoring \eqref{eq:surfaceAnchoringIntegral} integrals, as $\Gamma$ now indicates stronger \emph{anisotropic} surface tension on the spatial coordinates of the nematic LC molecules. The increasing elongation of the shape illustrates this. As Prinsen and van der Schoot show \cite{Prinsen2003}, the tactoid morphology depends on the balance of surface and bulk forces and on the ratio of the anisotropic to isotropic surface tension, $\frac{\Gamma}{\tau}$. With the molecular alignment fixed, their results illustrate that for $\omega \ll 1$, the aspect ratio is expected to scale as $1+\omega.$ Even though we fix $\tau = 10$ and bulk constants, $a,b,c$ in the Fixed Field model, we mimic their results numerically by increasing $\omega$ from $0.01$ to $1$ indicating that the shape's aspect ratio is $\omega-$dependent (see Figure \ref{fig:QNNImsmfRho}).

 \begin{figure}[!ht]
		\centering
		\begin{subfigure}{0.26\textwidth} 
			\includegraphics[width=0.9\linewidth]{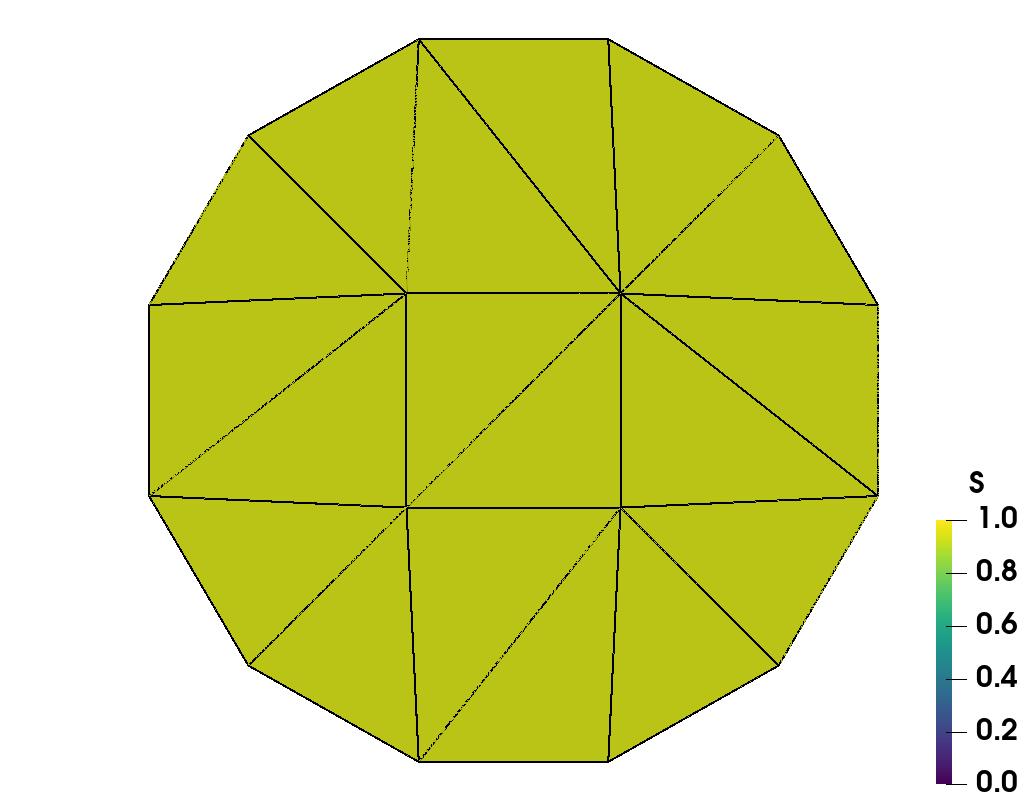}
			\caption{\footnotesize $|M_1|=16$, $\omega=0.01$}
		\end{subfigure}
		\hfill
		\begin{subfigure}{0.22\textwidth} 
			\includegraphics[width=1.1\linewidth]{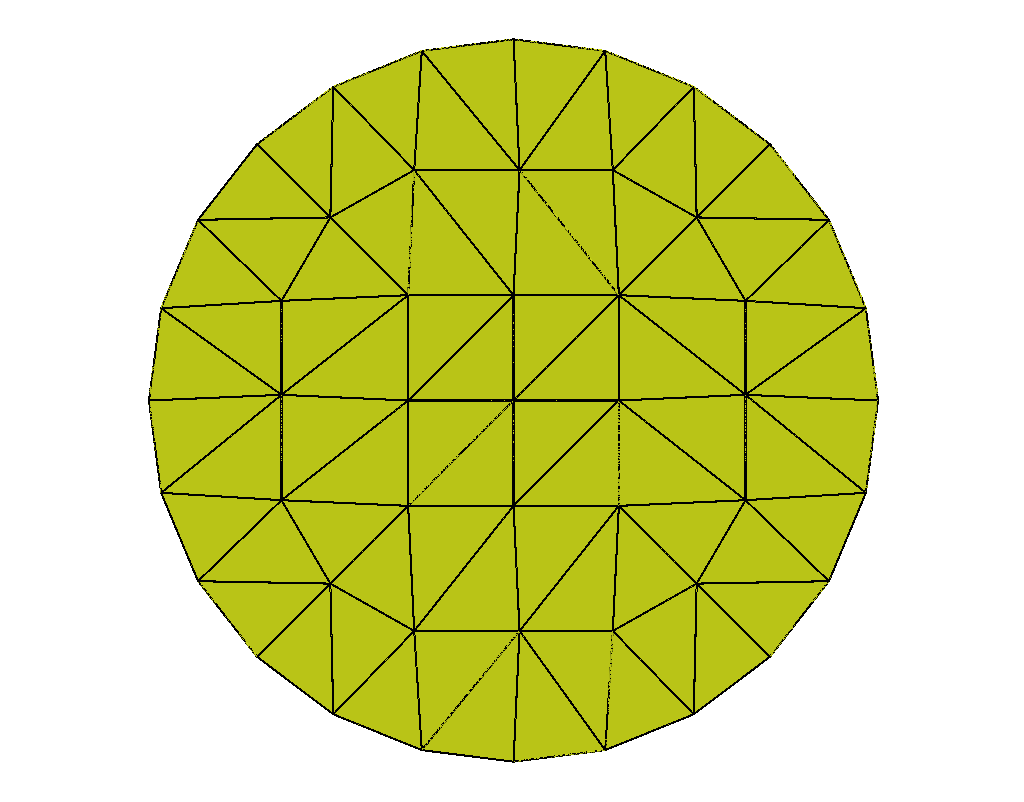}
			\caption{\footnotesize $|M_2|=49$}
		\end{subfigure}
		\hfill
		\begin{subfigure}{0.24\textwidth}
			\includegraphics[width=\linewidth]{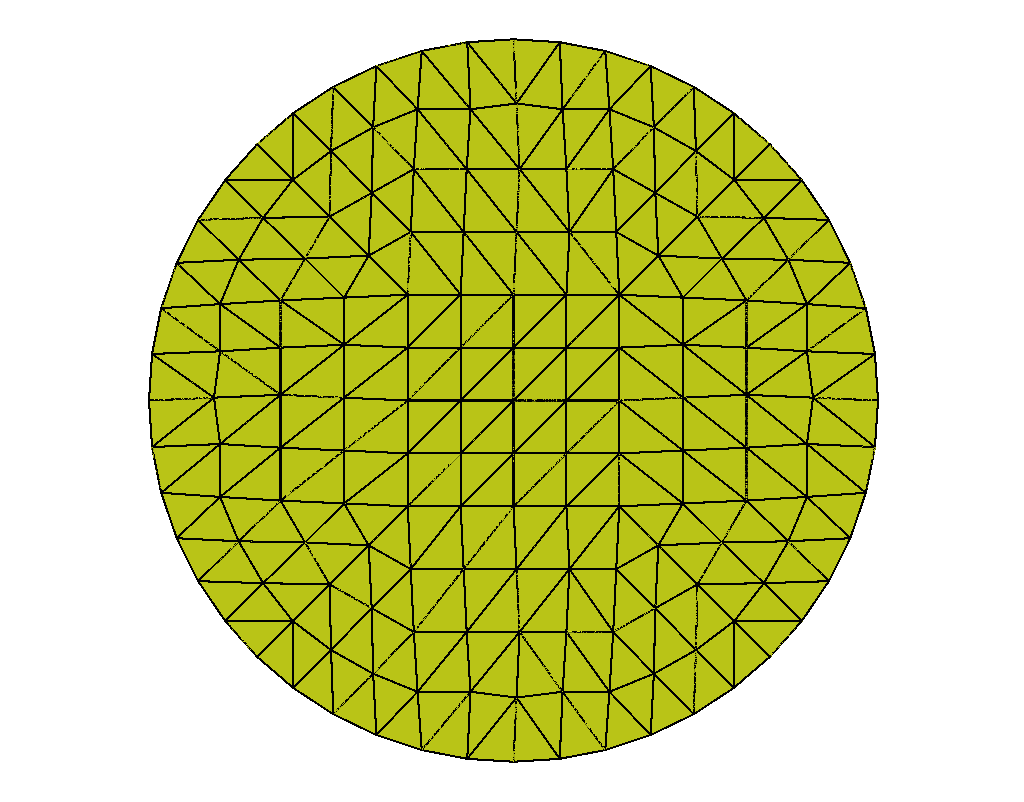}
			\caption{\footnotesize $|M_3|=169$}
		\end{subfigure}
        \begin{subfigure}{0.24\textwidth} 
			\includegraphics[width=\linewidth]{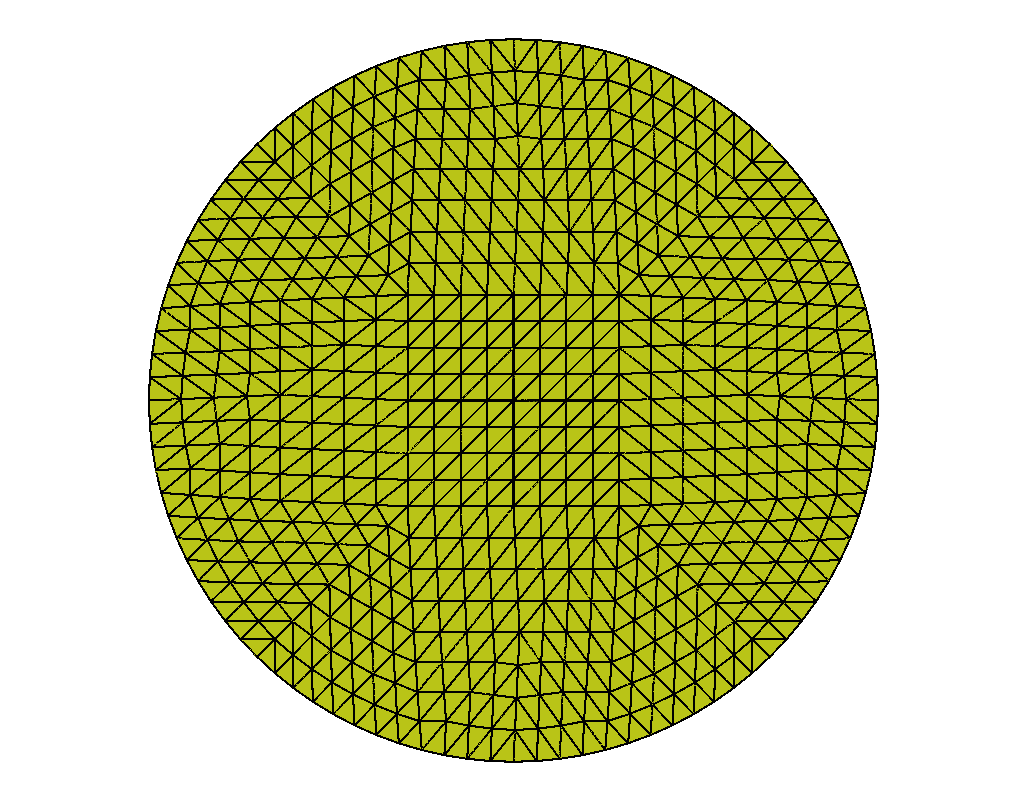}
			\caption{\footnotesize $|M_4|=625$}
		\end{subfigure}
		\hfill
		
		
		 \vspace{-0.35cm} 
		
		\begin{subfigure}{0.24\textwidth} 
			\includegraphics[width=\linewidth]{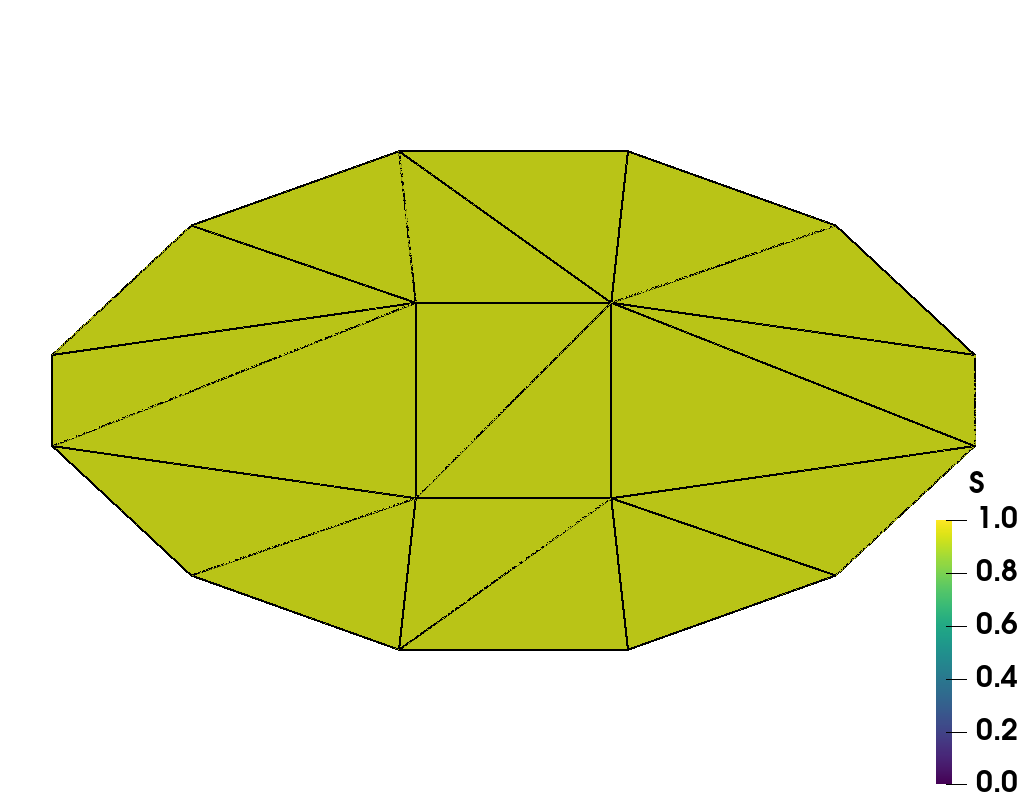}
			\caption{\footnotesize $|M_1|=16$, $\omega=1$}
		\end{subfigure}
		\hfill
		\begin{subfigure}{0.24\textwidth} 
			\includegraphics[width=\linewidth]{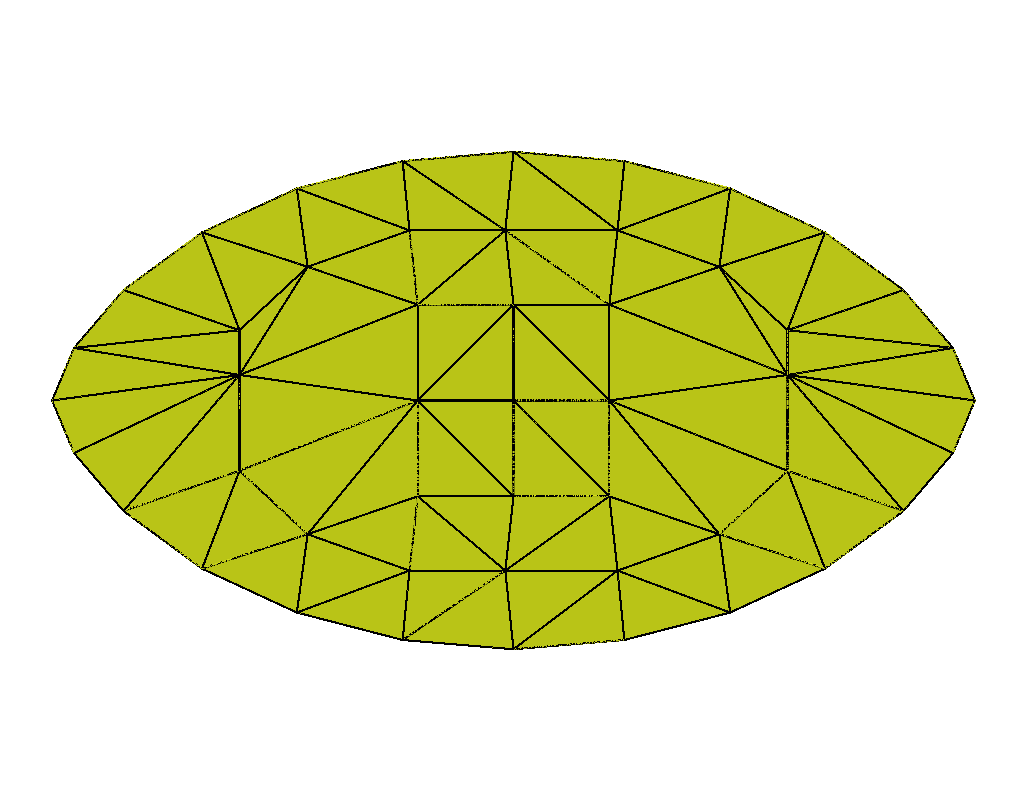}
			\caption{\footnotesize $|M_2|=49$}
		\end{subfigure}
		\hfill
		\begin{subfigure}{0.24\textwidth} 
			\includegraphics[width=\linewidth]{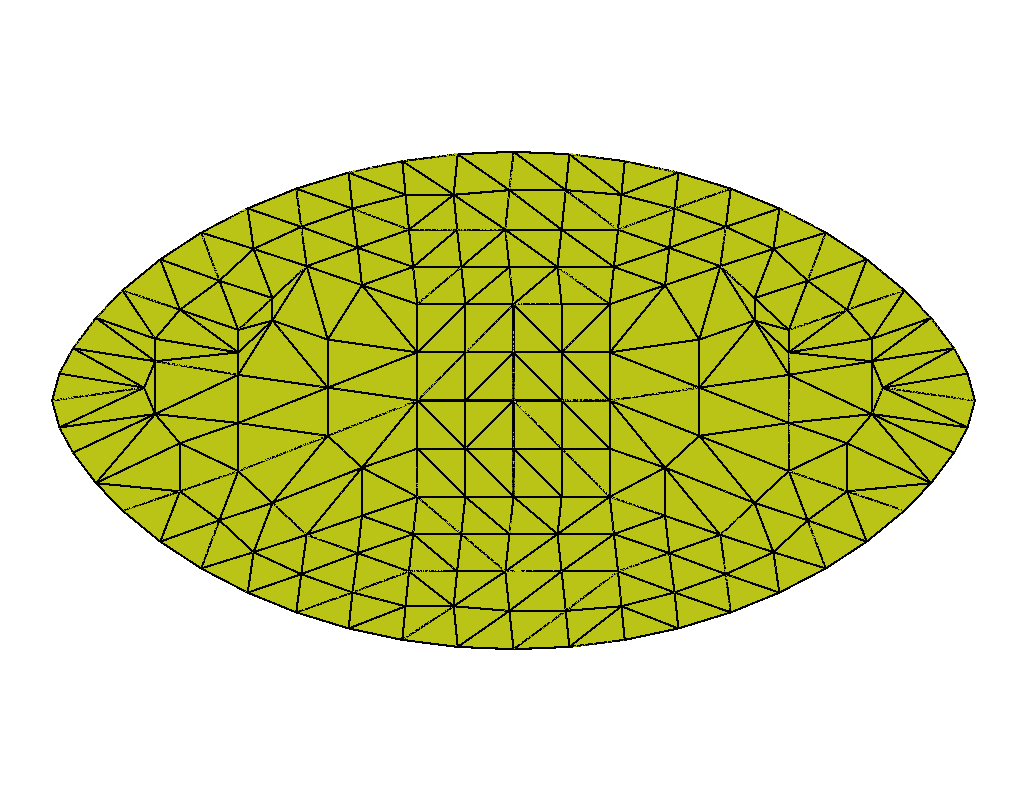}
			\caption{\footnotesize $|M_3|=169$}
		\end{subfigure}
        \hfill
		\begin{subfigure}{0.24\textwidth} 
			\includegraphics[width=\linewidth]{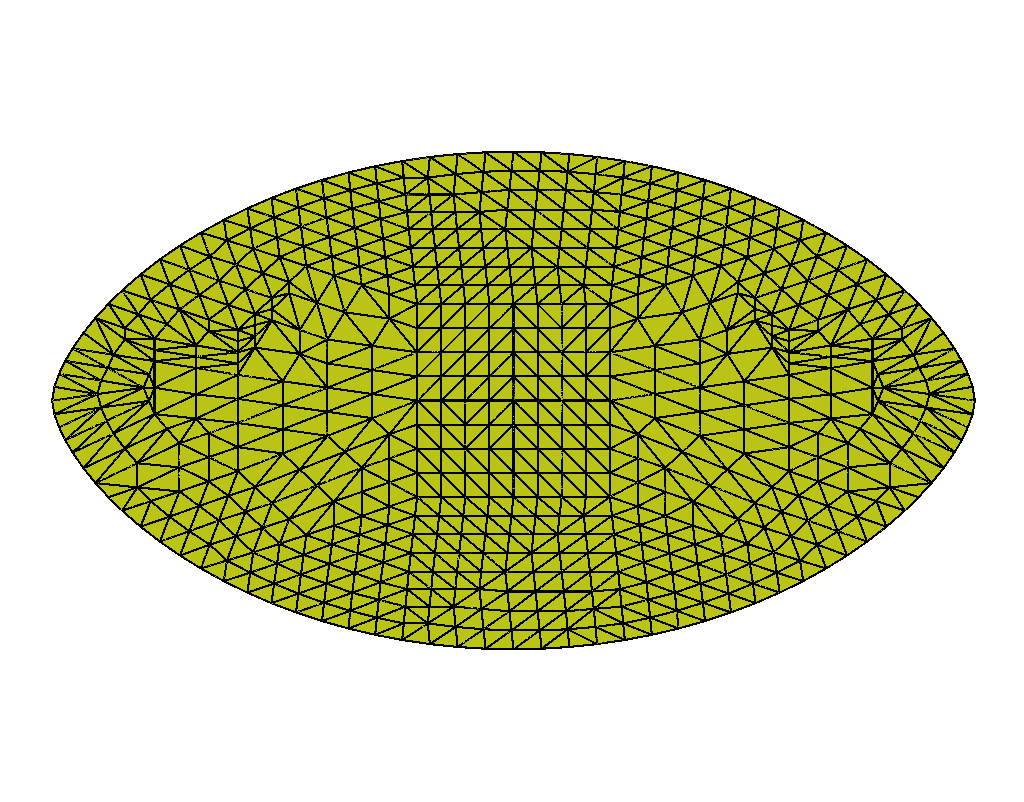}
			\caption{\footnotesize $|M_4|=625$}
		\end{subfigure}
  \hfill
		
		 \vspace{-0.1cm} 
		
		\caption{\justifying Applying QN with NI for Subproblem B (Fixed Field): Grids for each NI level for $\omega=0.01$ (top row) and $\omega=1$ (bottom row). 
  }
		\label{fig:msffgrids}
	\end{figure}
 \begin{figure}[!ht]
		\begin{subfigure}{0.49\textwidth}
			\centering
			\includegraphics[width=0.7\linewidth]{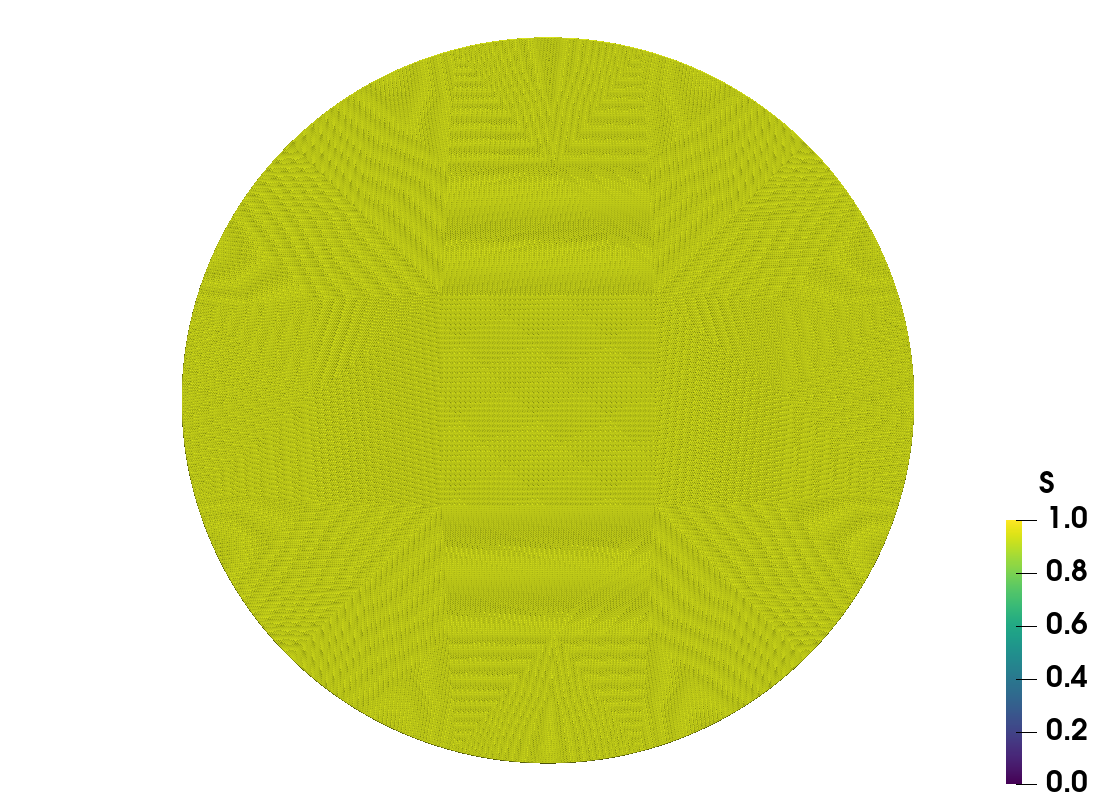}
			\caption*{(a) $\omega = 0.01$}
		\end{subfigure}%
		\begin{subfigure}{0.49\textwidth}
			\centering
			\includegraphics[width=0.7\linewidth]{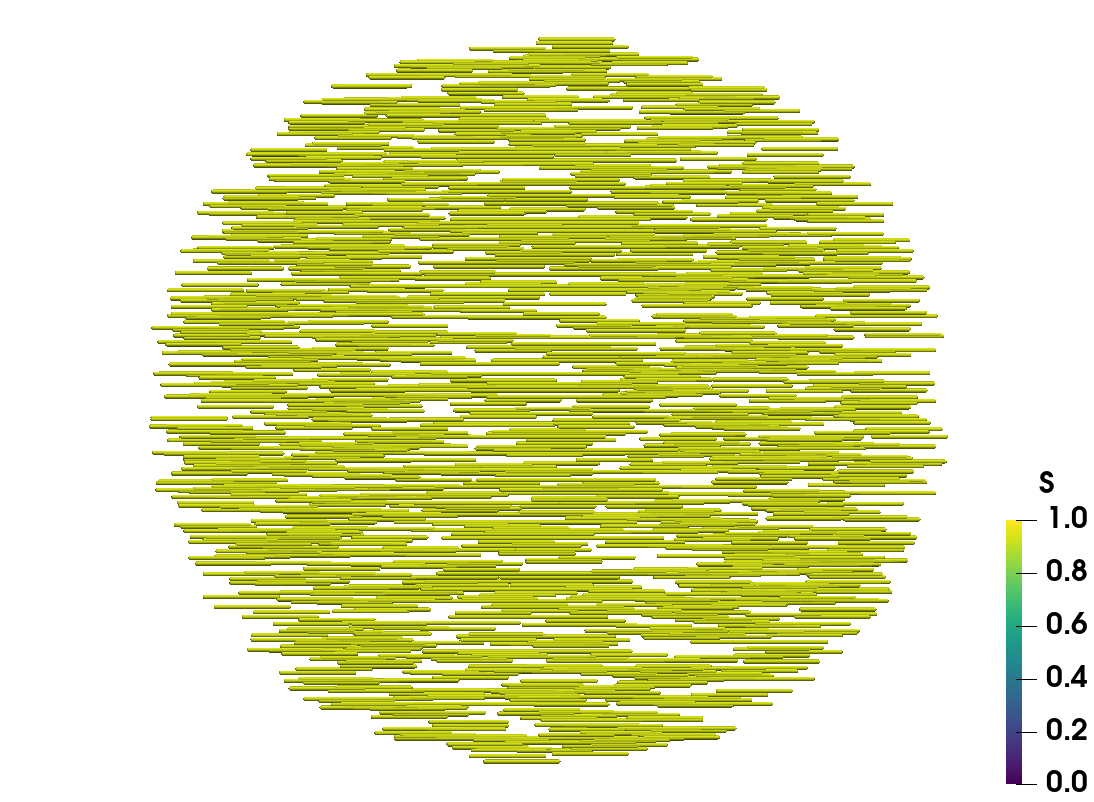}
			\caption*{(b)}
		\end{subfigure}\\
		\begin{subfigure}{0.49\textwidth}
			\centering
			\includegraphics[width=0.7\linewidth]{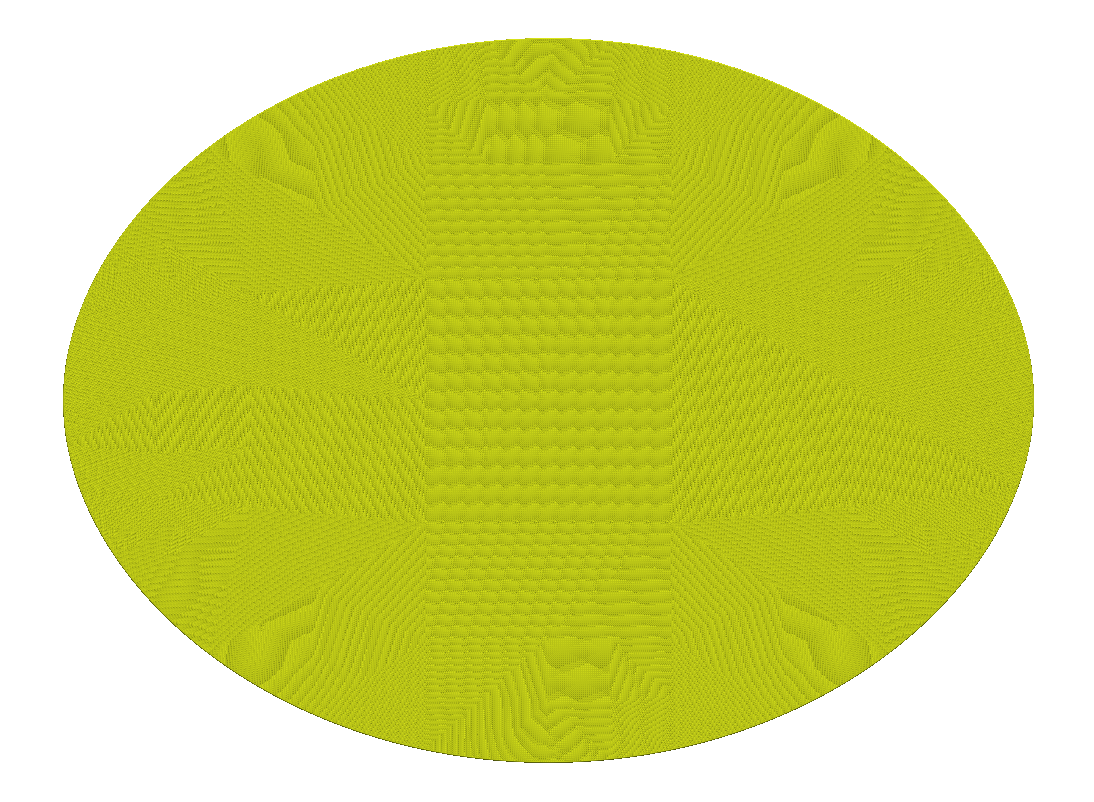}
			\caption*{(c) $\omega = 0.4$}
		\end{subfigure}%
		\begin{subfigure}{0.49\textwidth}
			\centering
			\includegraphics[width=0.7\linewidth]{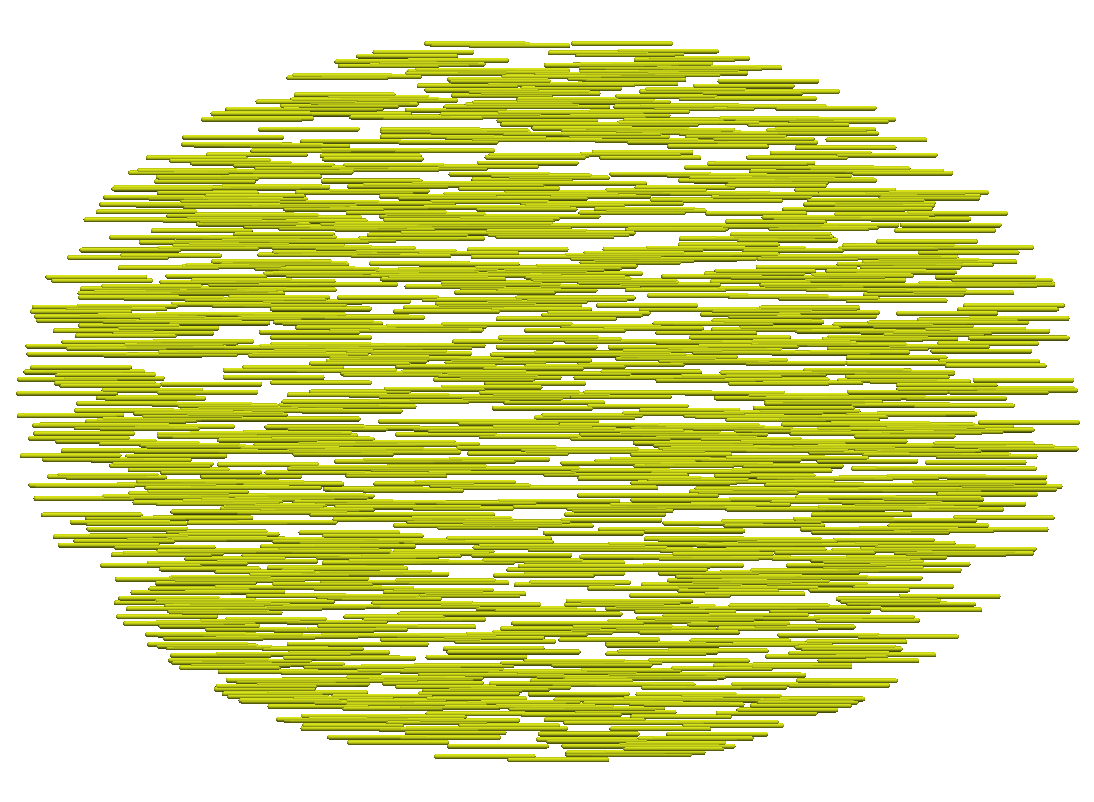}
			\caption*{(d)}
		\end{subfigure}\\
		\begin{subfigure}{0.49\textwidth}
			\centering
			\includegraphics[width=0.7\linewidth]{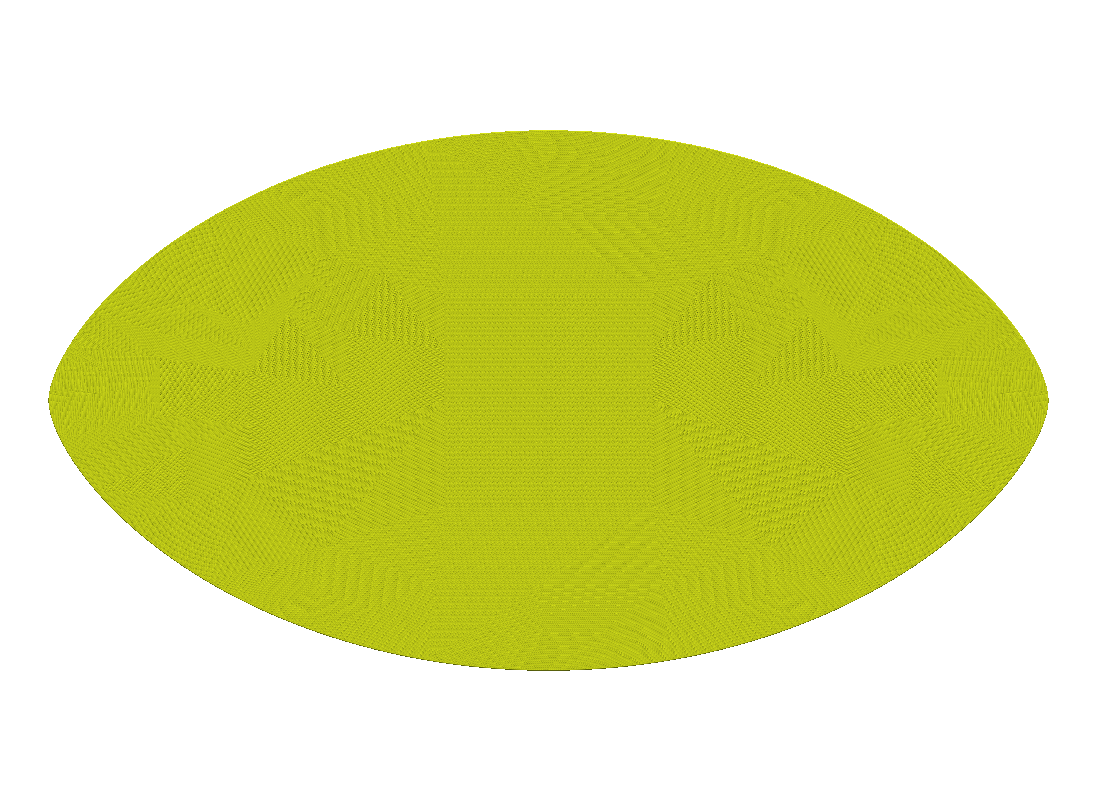}
			\caption*{(e) $\omega = 1$}
		\end{subfigure}%
		\begin{subfigure}{0.49\textwidth}
			\centering
			\includegraphics[width=0.7\linewidth]{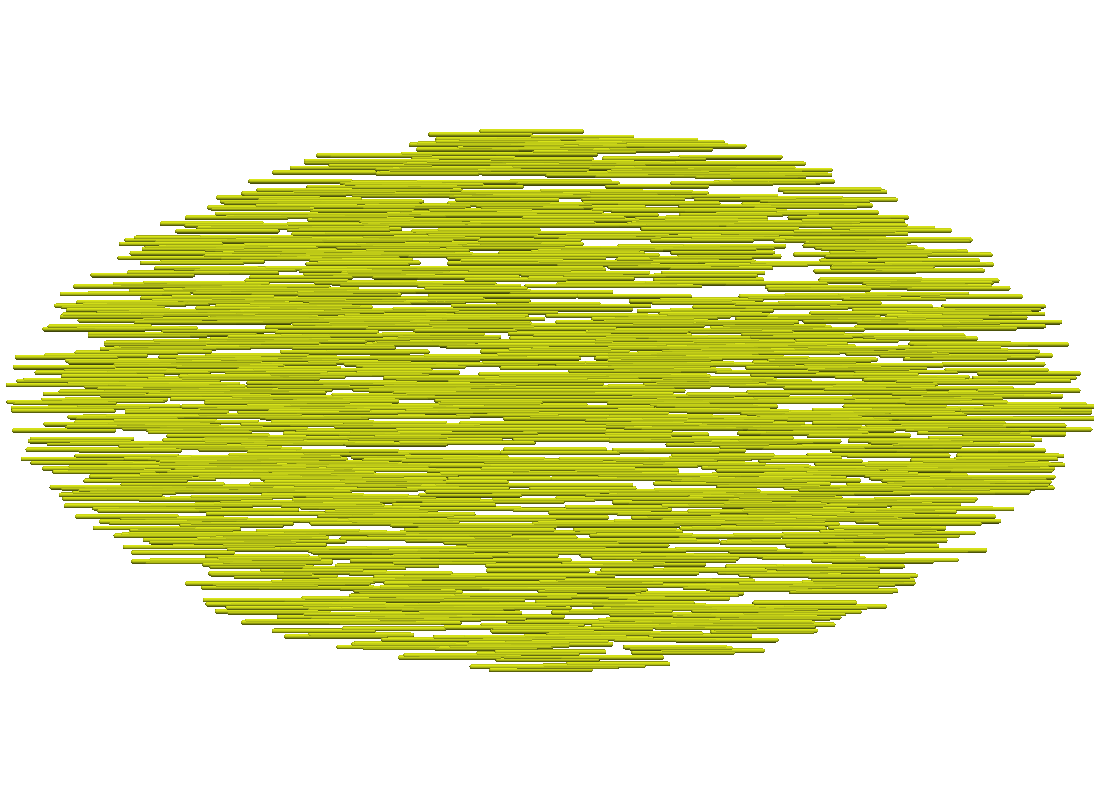}
			\caption*{(f)}
		\end{subfigure}
		\caption{\justifying Applying QN with NI for Subproblem B (Fixed Field). Results shown on grid $M_9$. The color bar indicates the value of $S$, i.e., the order of the director field in the domain. Left plots, (a), (c), (e), depict the order's distribution. Right plots, (b), (d), (f), show the fixed horizontally-aligned directors. Dramatic shape change is shown as $\omega \rightarrow 1$.}
		\label{fig:QNNImsff}
	\end{figure}
	
	We do not present standalone QN results as these did not converge for every value of $\omega$. As the shape changes, the initial circular guess is further from the optimal configuration, causing the method to stagnate. As expected, NI improves the initial guesses on each successive grid level, providing convergent results. 
 Figure \ref{fig:msffgrids} displays the hierarchy of adaptive grids, before any mesh regularization, used in the simulation for the extreme ends of shape change (i.e., $\omega = 0.01$ and $\omega = 1$). Each refined grid, $M_i$ is produced by quadrisecting the coarser elements from $M_{i-1}$. After each grid refinement, we regularize with equiangulation as a way to maintain a regular triangulation.  In Figure \ref{fig:QNNImsff}, we again show the distribution of the scalar order parameter $S$ (left plots), and the corresponding director field (right plots). We see the aspect ratio of the shape increasing as $\omega$ increases, validating Prinsen and van der Schoot's theory. Note that with the $Q$ held constant, the order of the director field is fixed at $S_0 = 0.933567$, and each director is horizontal to the $x-$axis. We do not expect to see defects forming, only the tactoid elongating as the anisotropic surface tension strength increases.  
 
	\subsubsection{Full $(\vec{X},Q)$ Optimization}\label{sec:Full2DOptimization}
The previous subproblems allow us to understand the individual effect on the presence of defects with $\tau$ and on the shape's trend of deformation driven by $\omega$. Numerically, we demonstrated that nested iteration significantly improved timings for examples with high orientational deformations and a notable improvement in convergence for regions of high spatial deformations. In this section, we allow both spatial and orientational movements, with the intention of comparing QN, derived in Section \ref{QN}, with the GD method in Section \ref{PGDDescription}, and seeing the effect of NI on both approaches. 

We validate the same results from Subproblems A \eqref{eq:fsmf} and B \eqref{eq:msff}, while tracking the distribution of the scalar order parameter $S$, the corresponding director field, the converged energy $F^k$, the number of iterations for each grid level, and the runtime in seconds. In these numerical experiments, we mimic the Fixed Field model \eqref{eq:msff} tests by holding $\tau$ fixed at $10$ and varying $\omega$ from $0.01$ to $1$. Since we allow for spatial and orientational displacements at once, we expect to see an effect from both isotropic surface tension \eqref{eq:lineTensionIntegral} and anchoring  \eqref{eq:surfaceAnchoringIntegral} effects. Visually, in this case, we see the colloidal particles elongating and two defects (emerging possibly outside the tactoid) on the opposite ends of the shape as $\omega$ tends to $1$. This is the combined effect of the dynamic interplay of the isotropic and anisotropic strengths on tactoids.
	
\begin{figure}[!ht]
		\centering
		\begin{subfigure}{0.26\textwidth} 
			\includegraphics[width=0.9\linewidth]{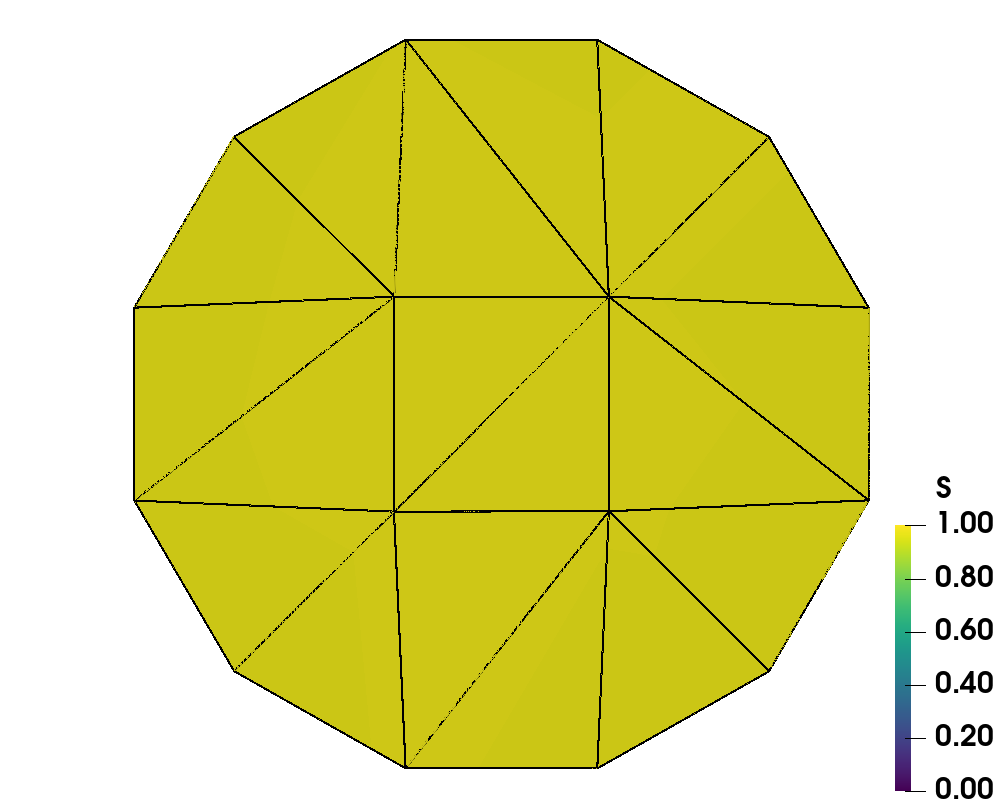} 
			\caption{\footnotesize $|M_1|=16$, $\omega=0.01$}
		\end{subfigure}
		\hfill
		\begin{subfigure}{0.22\textwidth} 
			\includegraphics[width=1.1\linewidth]{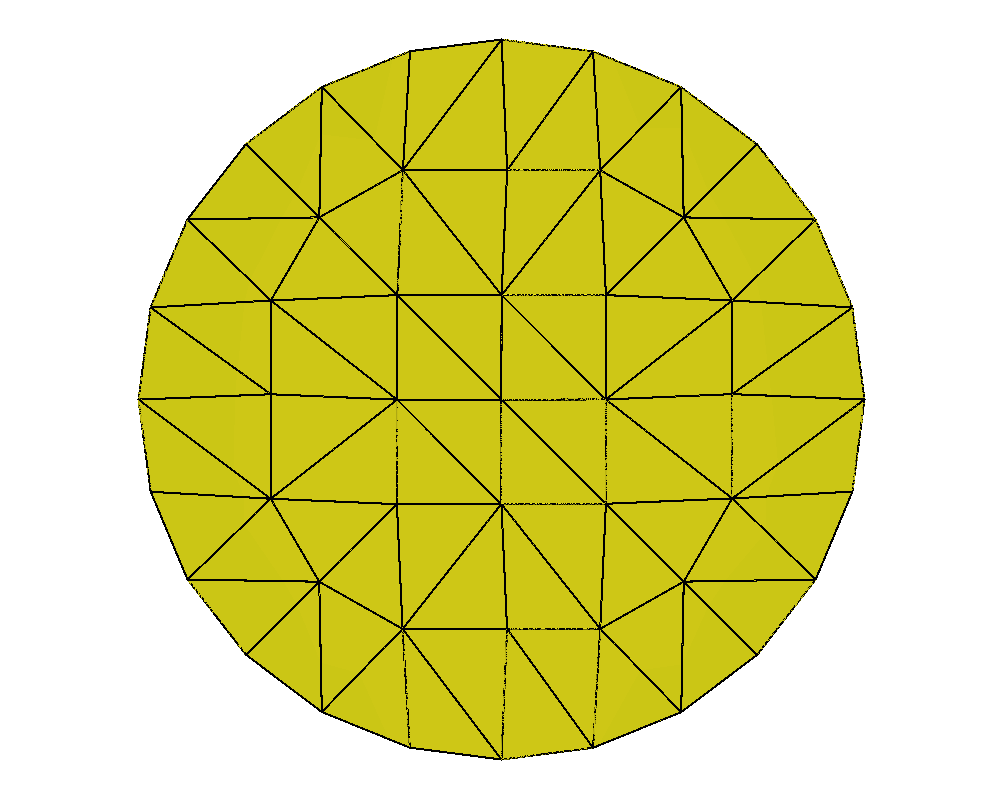} 
			\caption{\footnotesize $|M_2|=49$}
		\end{subfigure}
		\hfill
		\begin{subfigure}{0.24\textwidth} 
			\includegraphics[width=1\linewidth]{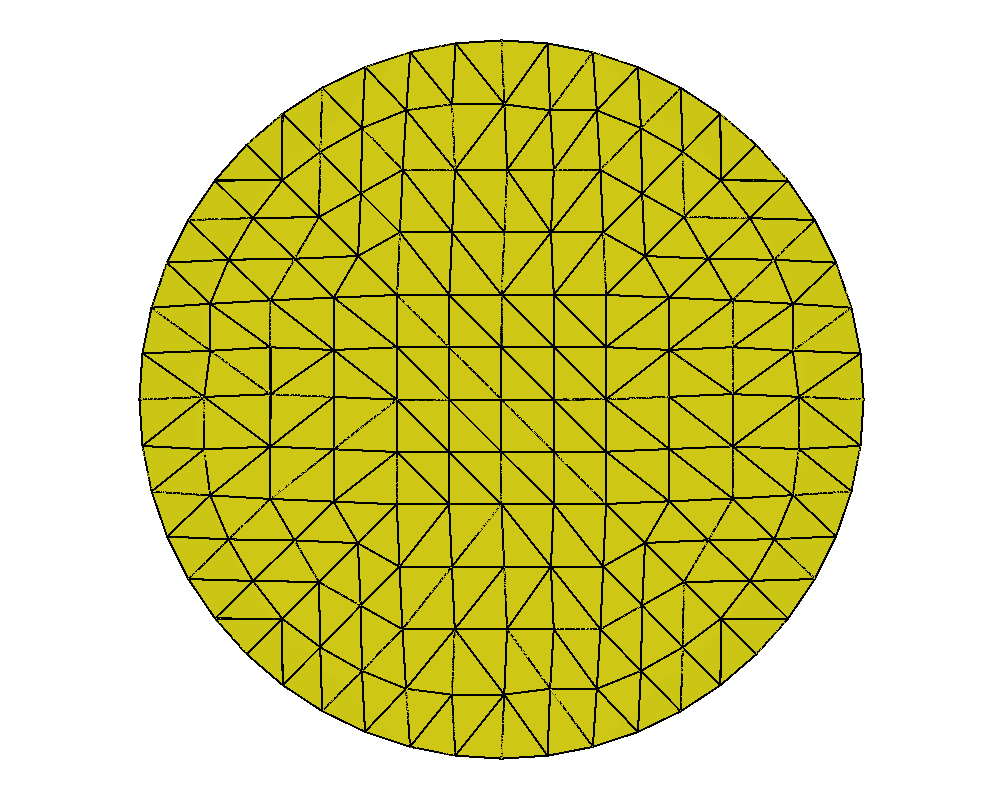} 
			\caption{\footnotesize $|M_3|=169$}
		\end{subfigure}
	    \hfill
		\begin{subfigure}{0.24\textwidth} 
			\includegraphics[width=1\linewidth]{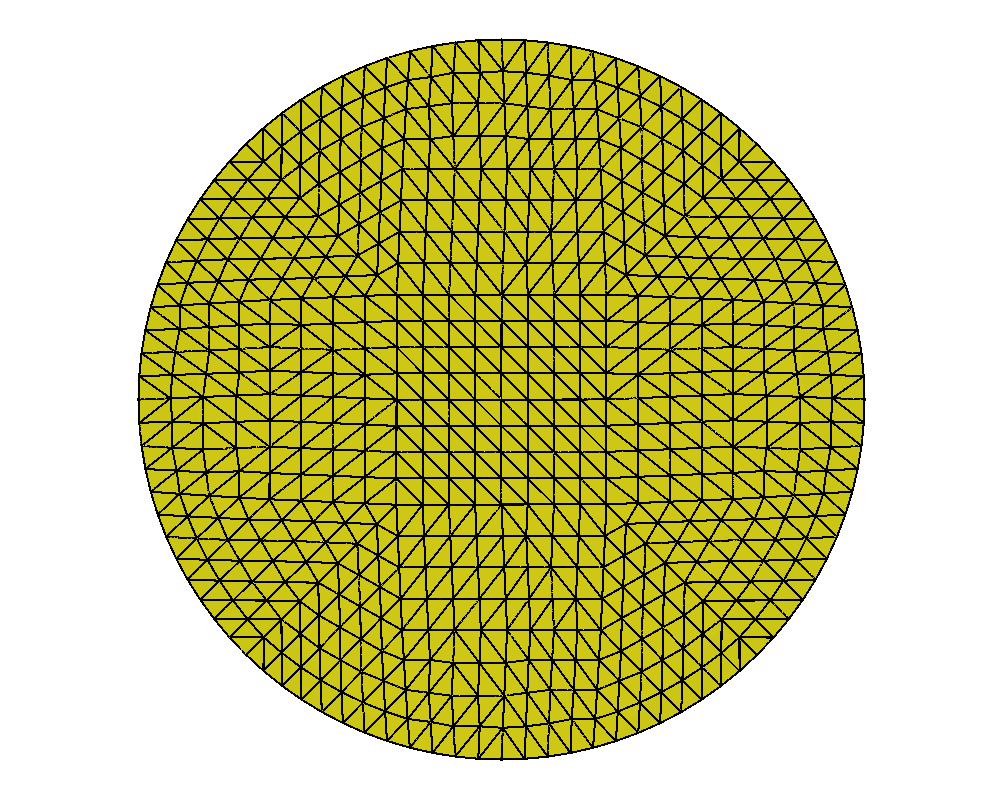} 
			\caption{\footnotesize $|M_4|=625$}
		\end{subfigure}
		\vspace{0.1cm} 
		
		
		
		\vspace{0.1cm}
		
		\begin{subfigure}{0.24\textwidth} 
			\includegraphics[width=0.95\linewidth]{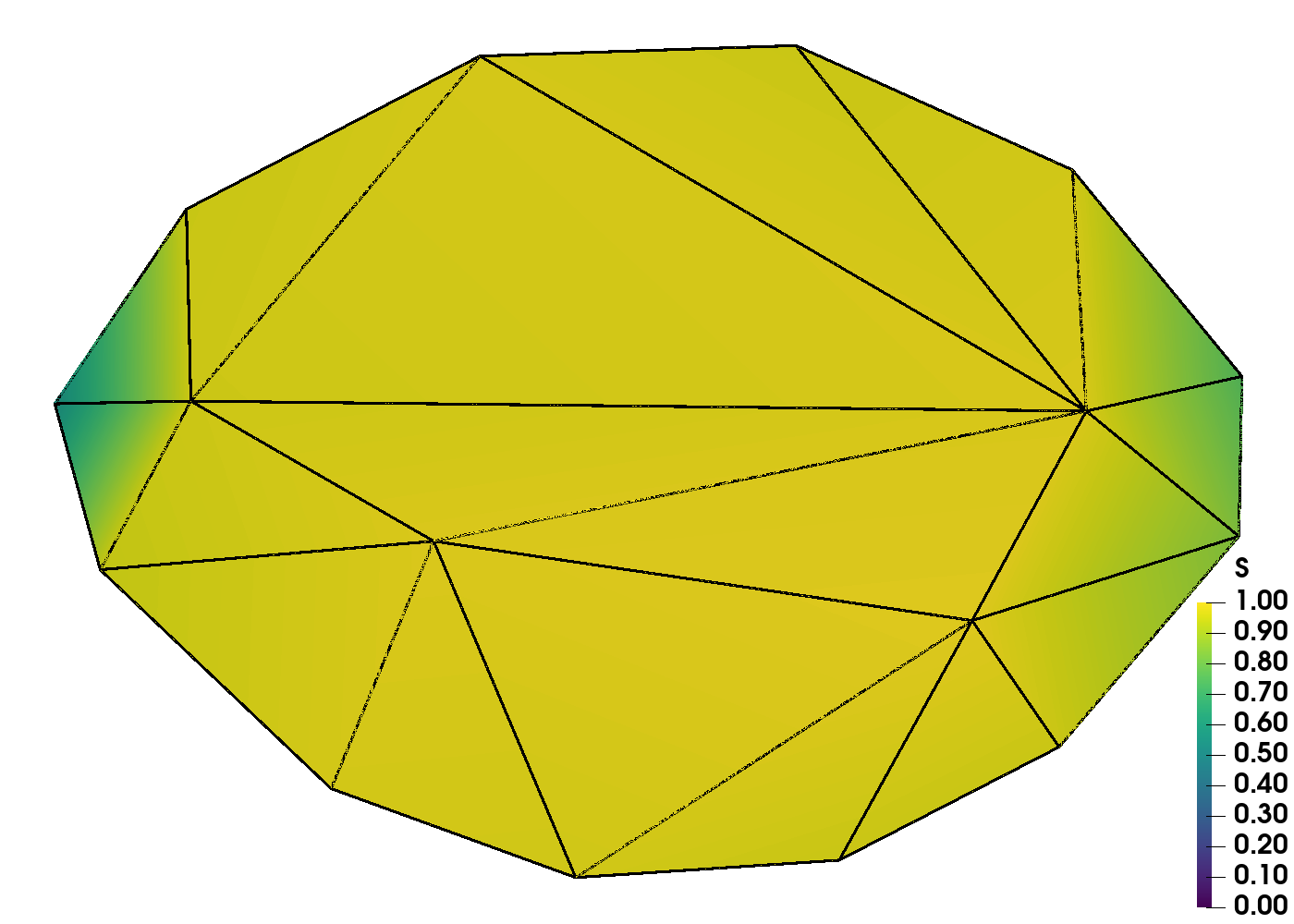} 
			\caption{\footnotesize $|M_1|=16$, $\omega=1$}
		\end{subfigure}
		\hfill
		\begin{subfigure}{0.24\textwidth} 
			\includegraphics[width=1\linewidth]{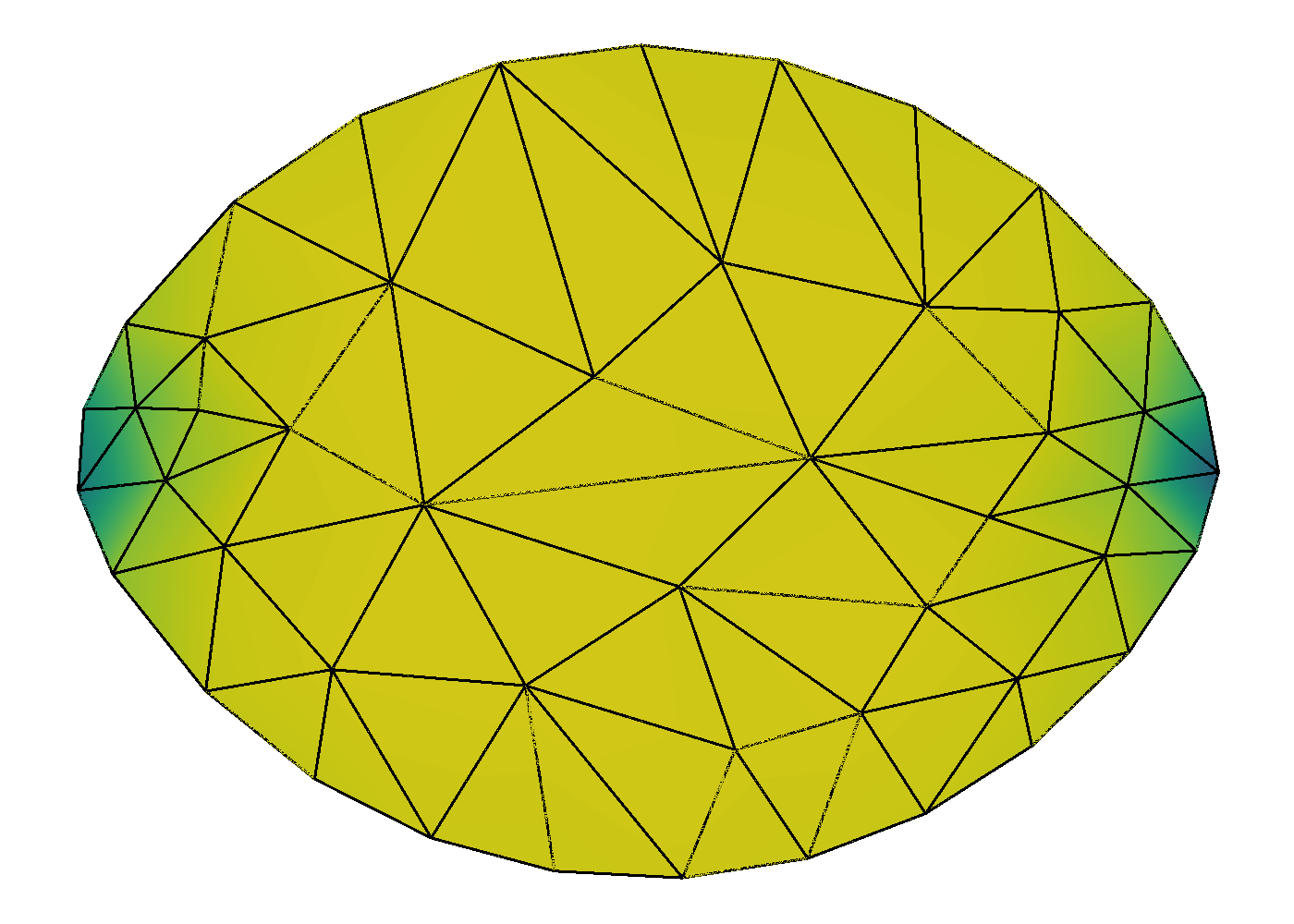} 
			\caption{\footnotesize $|M_2|=49$}
		\end{subfigure}
		\hfill
		\begin{subfigure}{0.24\textwidth} 
			\includegraphics[width=1\linewidth]{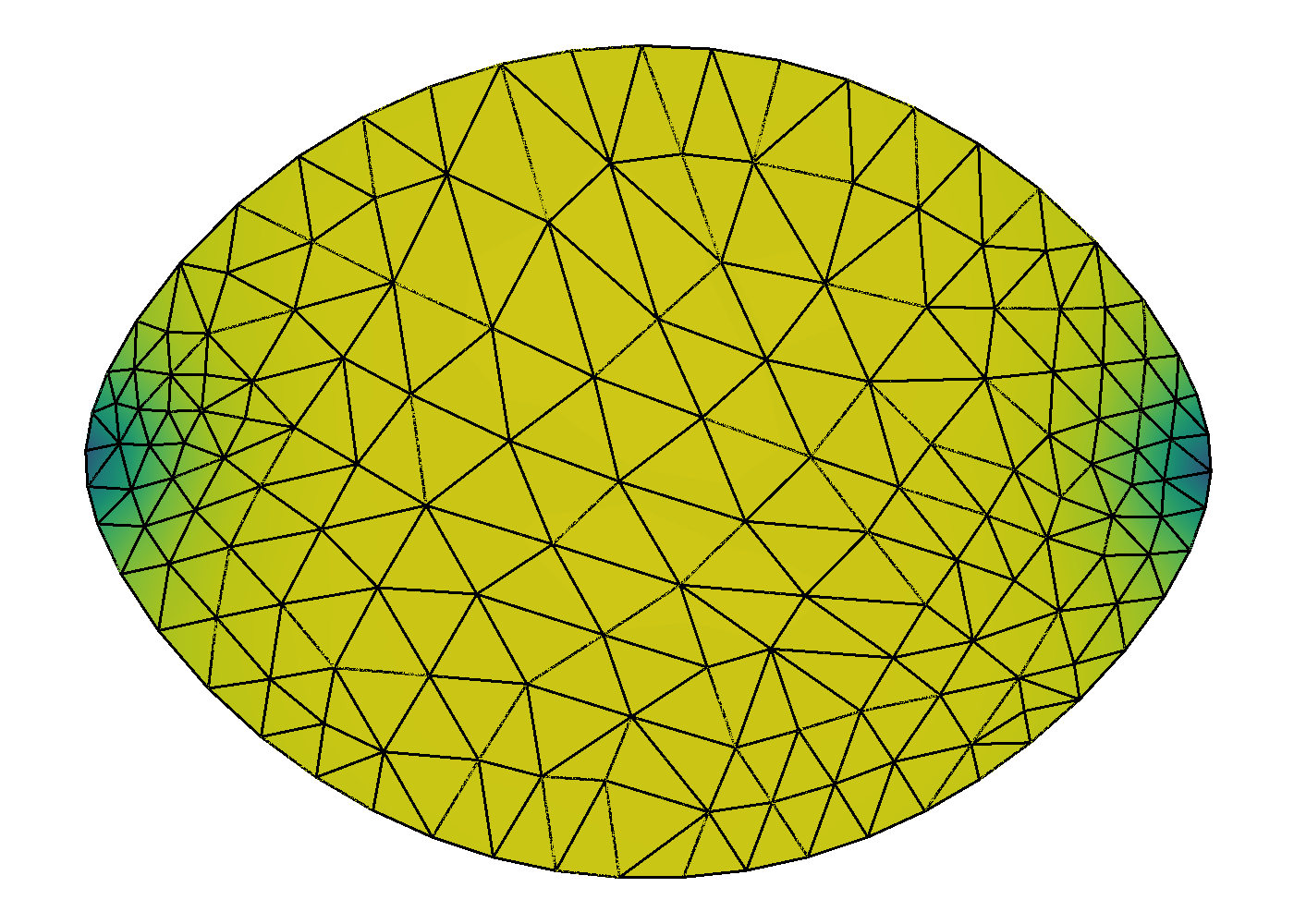} 
			\caption{\footnotesize $|M_3|=169$}
		\end{subfigure}
        \hfill
        \begin{subfigure}{0.24\textwidth} 
			\includegraphics[width=1\linewidth]{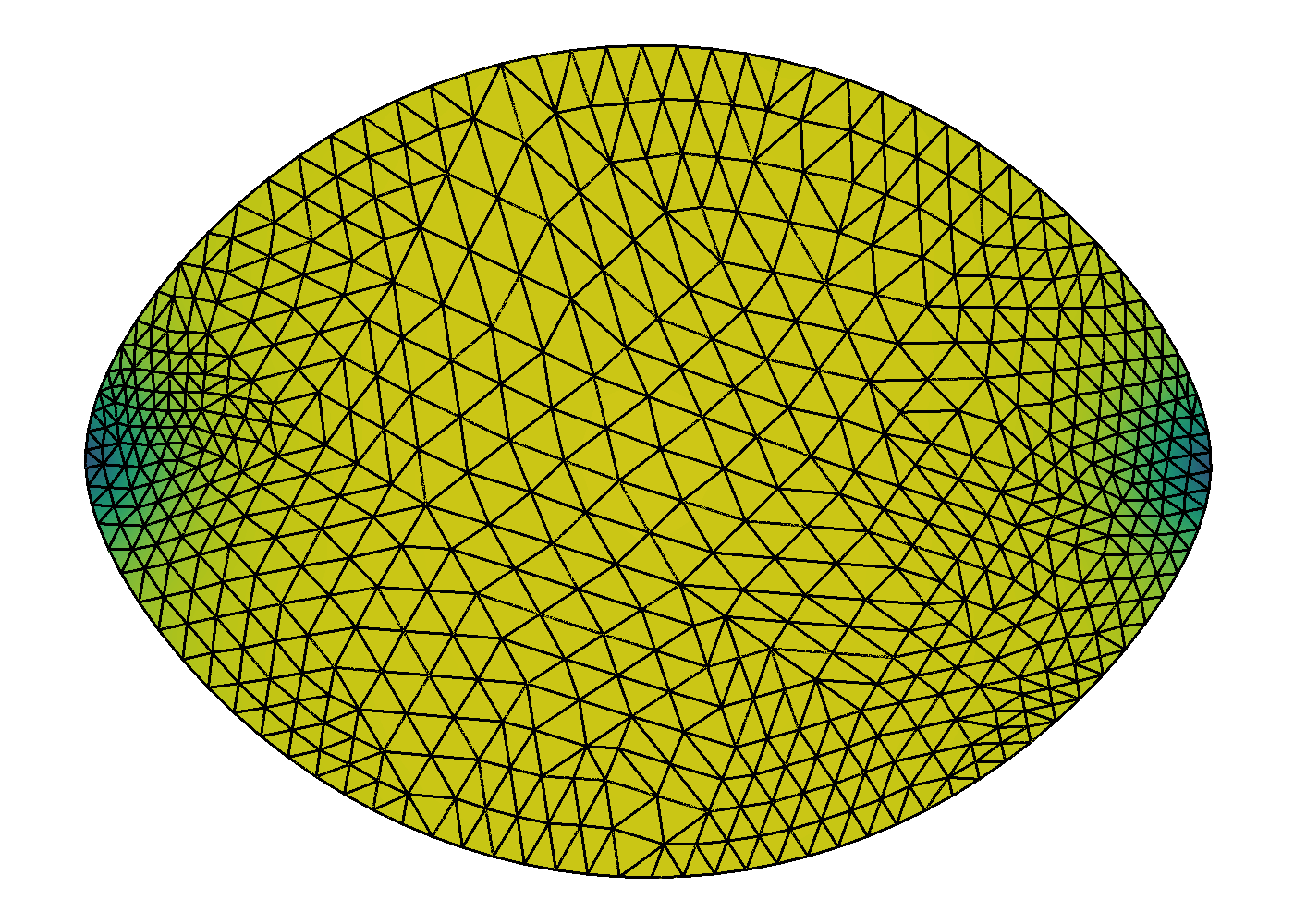}
			\caption{\footnotesize $|M_4|=625$}
		\end{subfigure}
		

		\caption{\justifying Applying QN with NI for the full $(\vec{X},Q)$ $2D$ problem: Grids for each NI level for $\omega=0.01$ (top row) and $\omega=1$ (bottom row). 
  }
		\label{fig:msmfgrids}
	\end{figure}

Similarly to the previous subproblems, we show the NI grid progression on the levels for $\omega= 0.01$ and $\omega = 1$. Figure \ref{fig:msmfgrids} illustrates the finite-element mesh evenly distributed throughout the grids in the absence of defects ($\omega=0.01$). For $\omega=1$, the vertices of the mesh tend to collect around the two defects as the grids get finer. For this problem, we also perform equiangulation after each grid refinement to maintain decent mesh quality. 
	
\begin{figure}[!ht]
		\begin{subfigure}{0.49\textwidth}
			\centering
			\includegraphics[width=\linewidth]{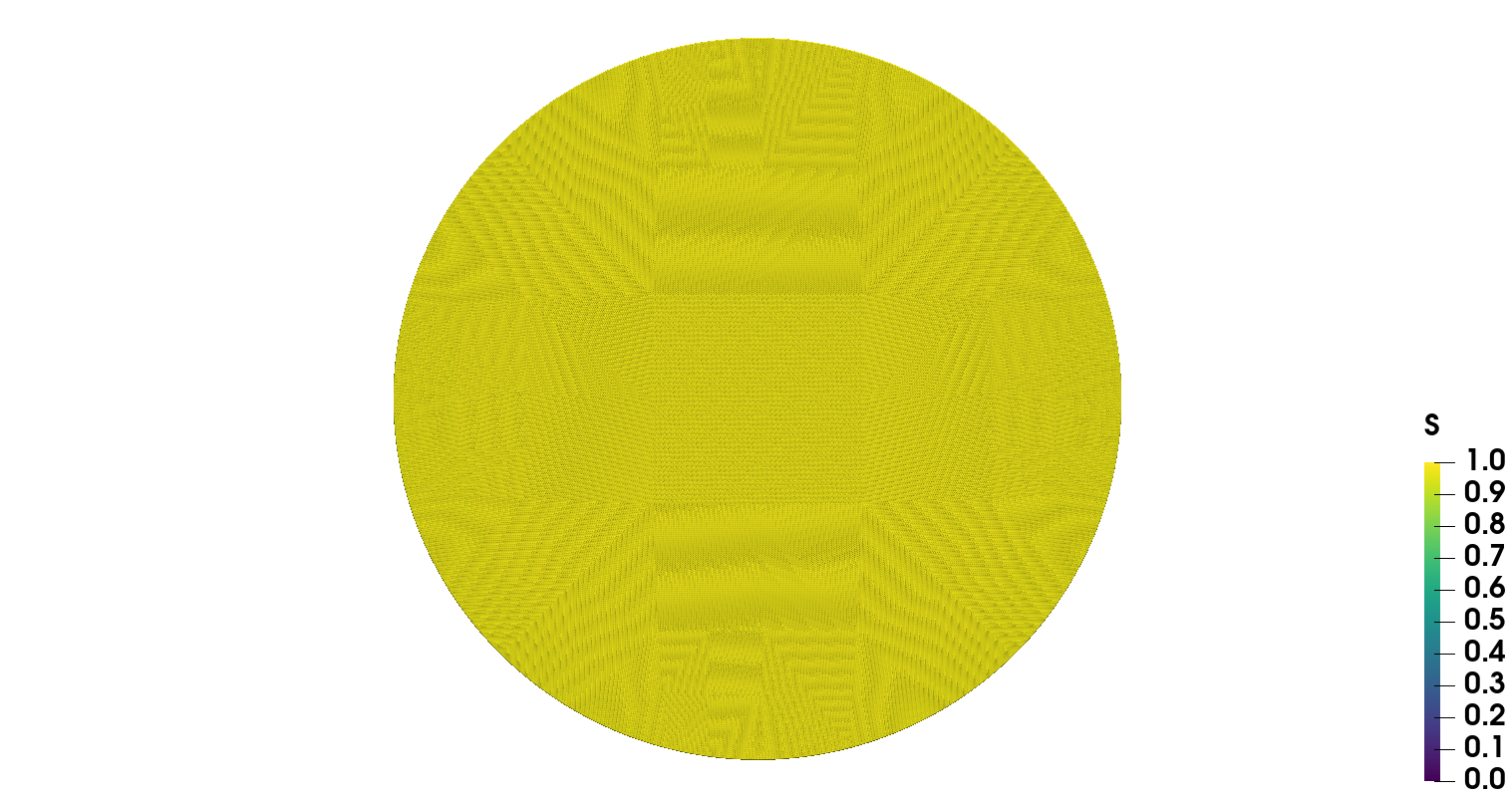}
			\caption*{(a) $\omega = 0.01$}
		\end{subfigure}%
		\begin{subfigure}{0.49\textwidth}
			\centering
			\includegraphics[width=\linewidth]{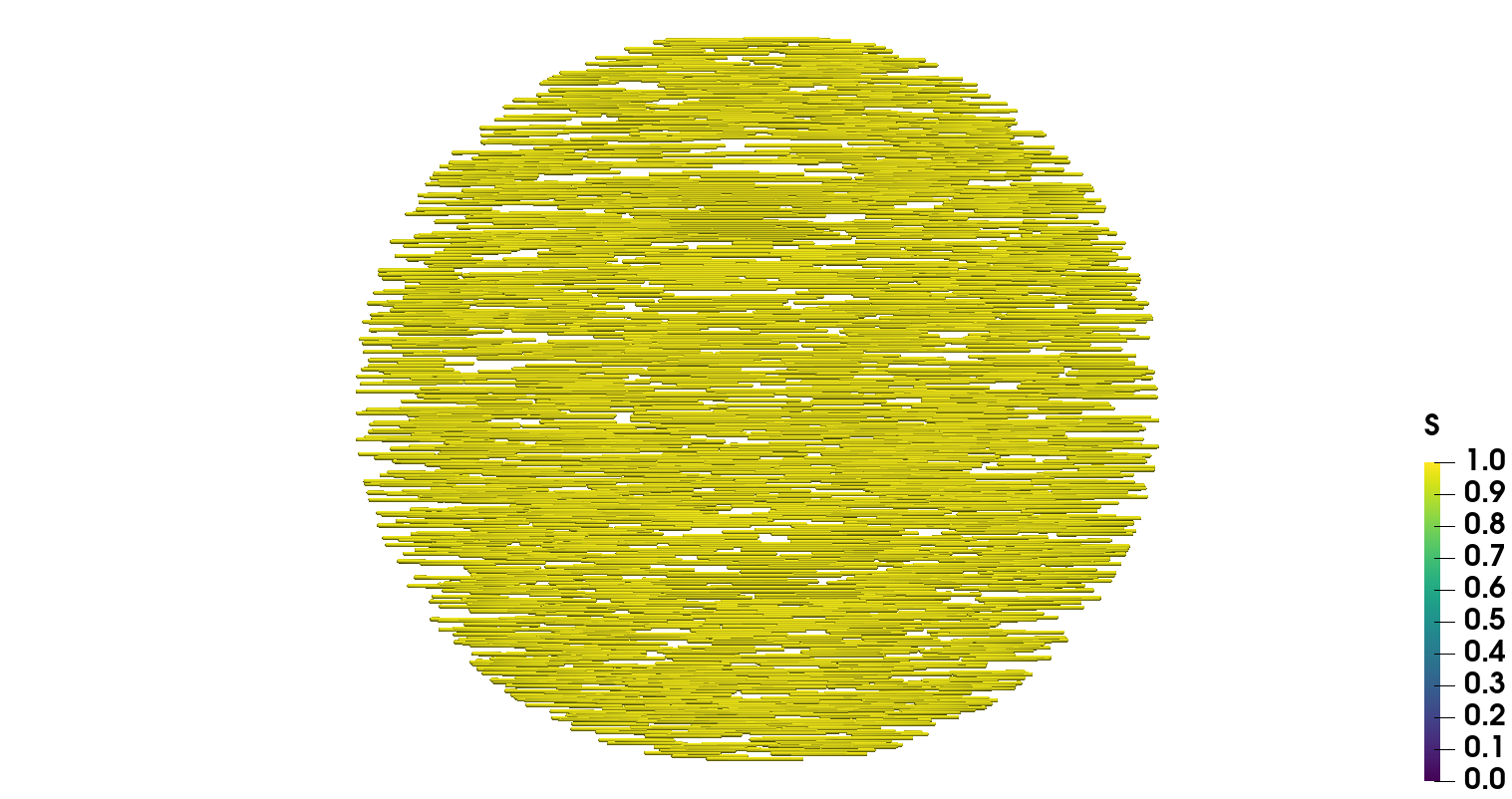}
			\caption*{(b)}
		\end{subfigure}\\
        \begin{subfigure}{0.49\textwidth}
			\centering
			\includegraphics[width=\linewidth]{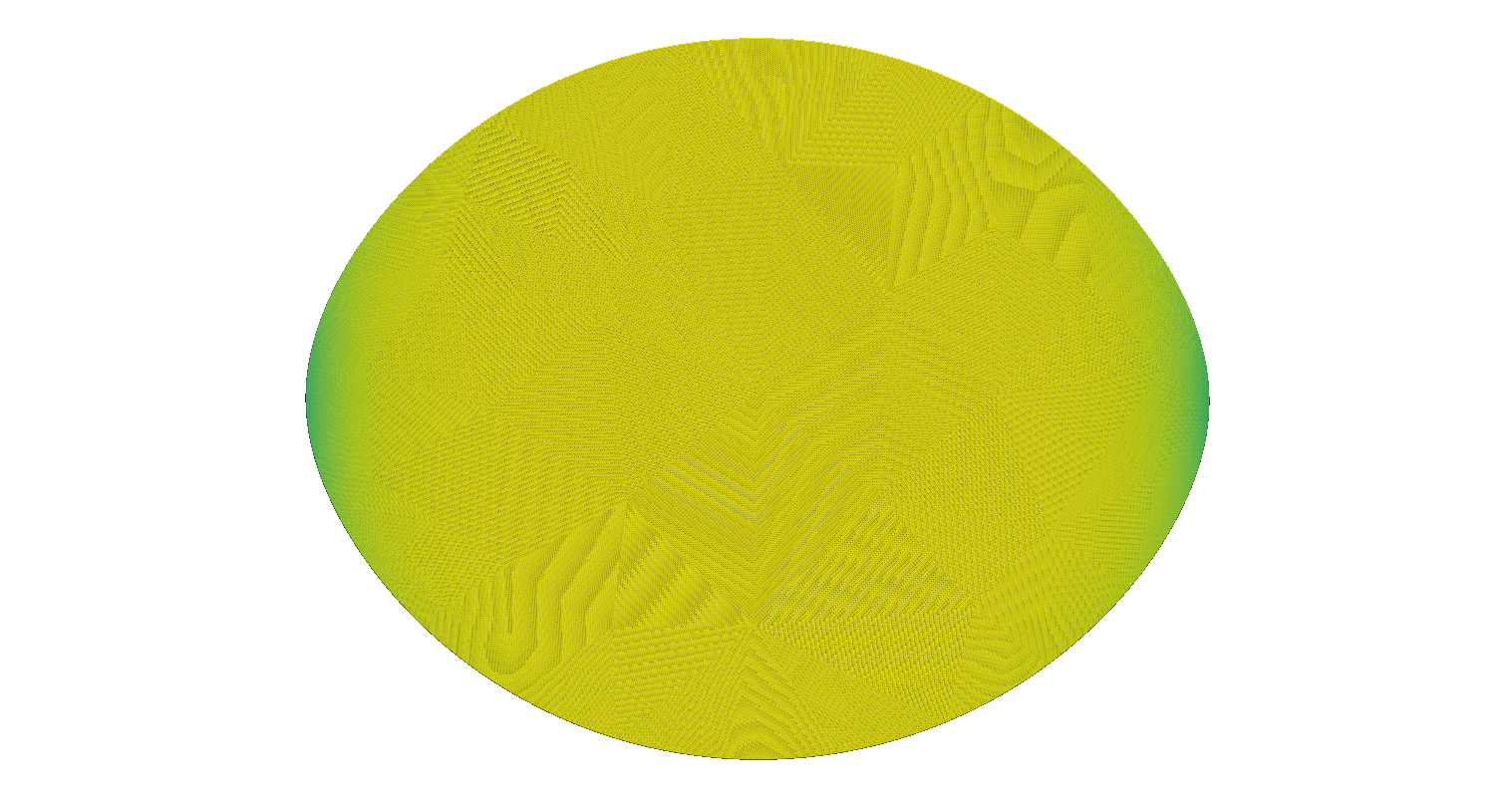}
			\caption*{(c) $\omega = 0.4$}
		\end{subfigure}%
		\begin{subfigure}{0.49\textwidth}
			\centering
			\includegraphics[width=\linewidth]{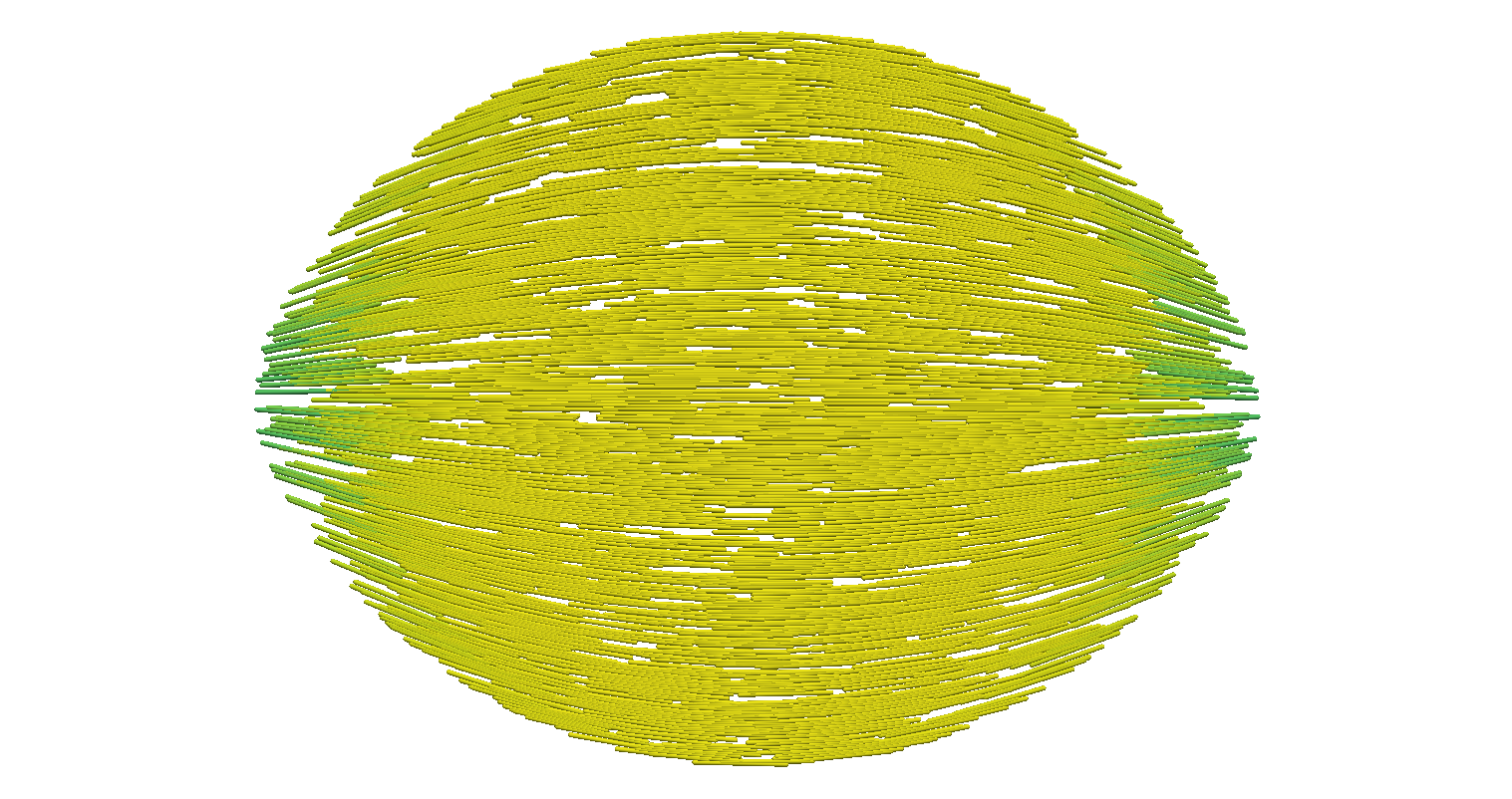}
			\caption*{(d)}
		\end{subfigure}\\
		\begin{subfigure}{0.49\textwidth}
			\centering
			\includegraphics[width=\linewidth]{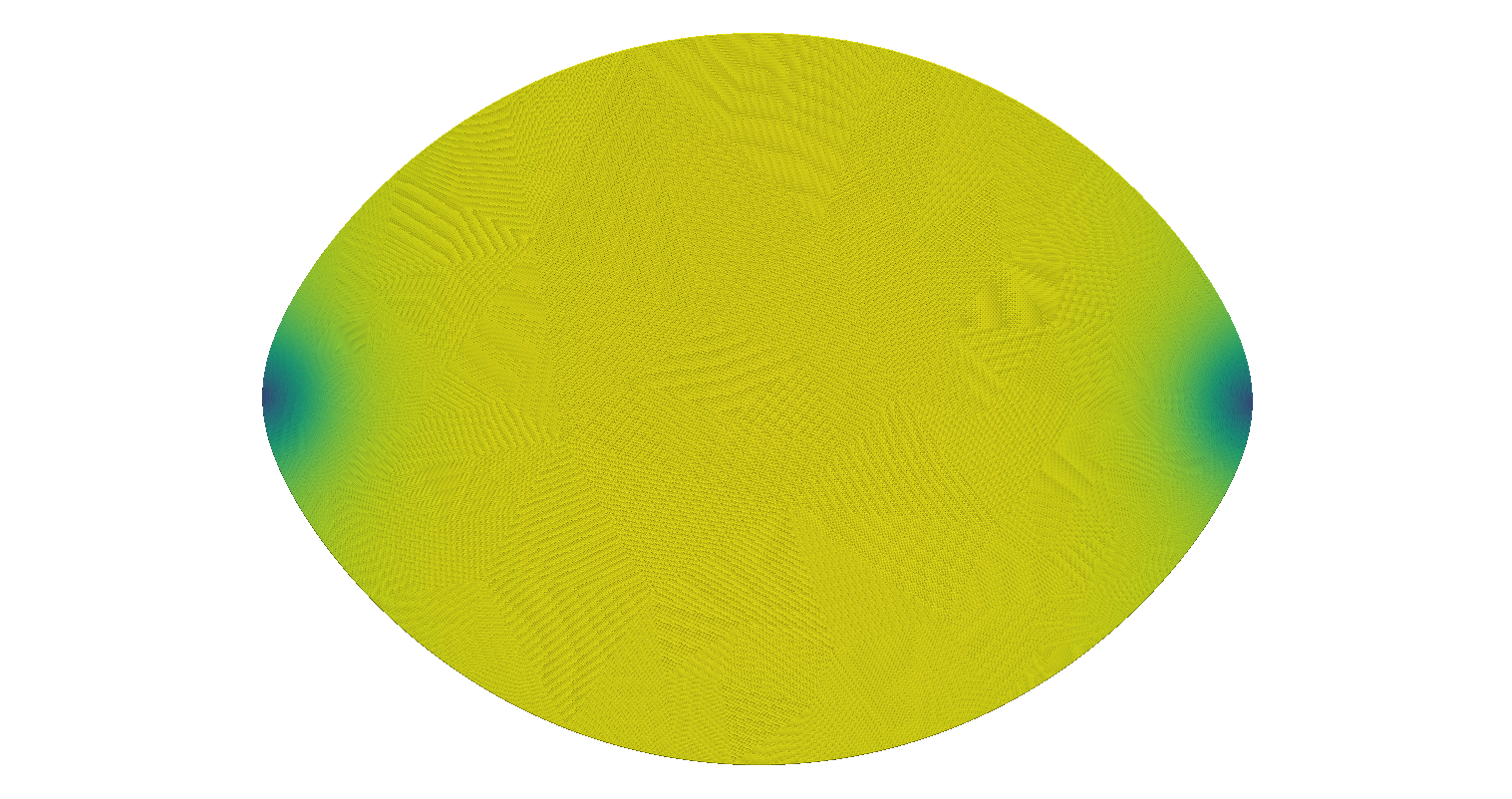}
			\caption*{(e) $\omega = 1$}
		\end{subfigure}%
		\begin{subfigure}{0.49\textwidth}
			\centering
			\includegraphics[width=\linewidth]{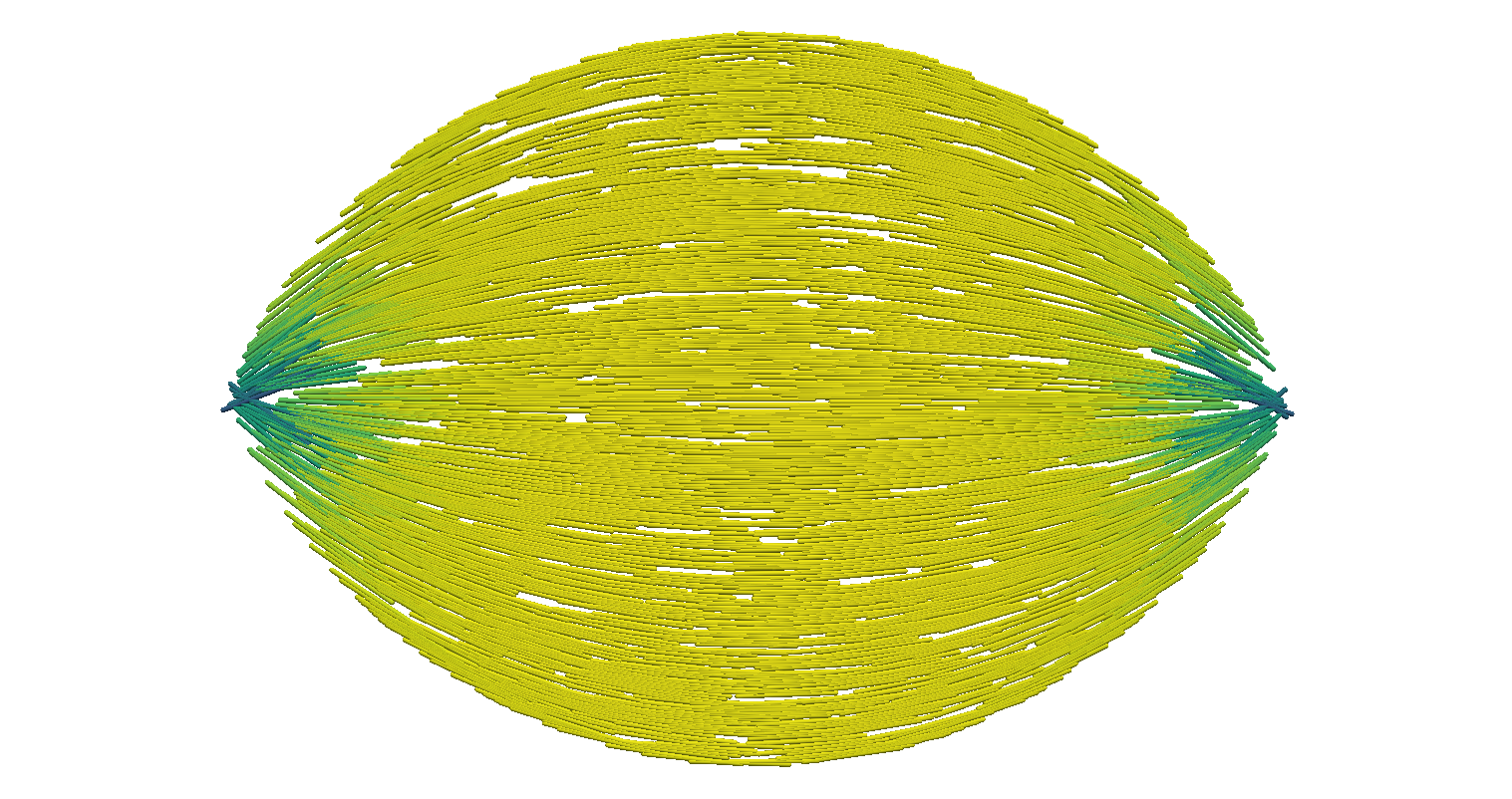}
			\caption*{(f)}
		\end{subfigure}
		\caption{\justifying Applying QN with NI to the full $(\vec{X},Q)$ $2D$ problem. Results shown on grid $M_9$. The color bar indicates the value of $S$, i.e., the order of the director field in the domain. Left plots, (a), (c), (e), depict the order's distribution. Right plots, (b), (d), (f), show the directors with stronger anchoring as the shape changes with $\omega\rightarrow 1$. Areas in green indicate less order, showing the appearance of the defects, as expected.}
		\label{fig:QNNImsmf}
	\end{figure}

\begin{figure}[!ht]
		\centering
		\includegraphics[width = \columnwidth]{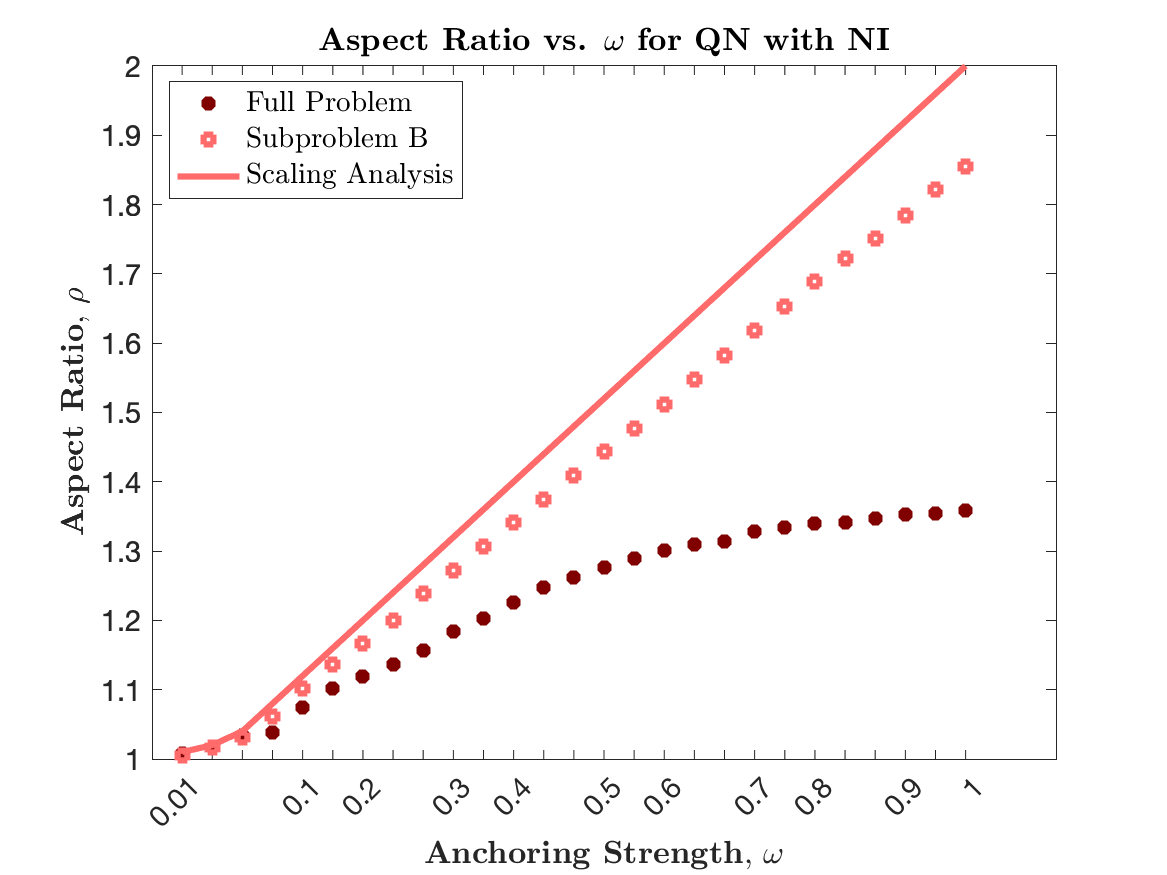}
		\caption{\justifying Applying QN with NI to the full $(\vec{X},Q)$ $2D$ problem and Subproblem B: With the presence of molecular alignment (dark red) for the full problem, we see a logarithmic effect in aspect ratio compared to the linear trend (light red) for fixed alignment in Subproblem B. For fixed alignment, this is in agreement with the scaling analysis of \cite{Prinsen2003} (solid line).}
		\label{fig:QNNImsmfRho}
	\end{figure}

Next, Figure \ref{fig:QNNImsmf} exhibits the combined effects from shape change and orientational variance as $\omega$ tends to $1$ using QN with NI. The left graphics of the figure show the distribution of the scalar order parameter $S$, and the corresponding director field is given on the right. As expected, we see non-spherical shapes with increasing anisotropic surface tension. Consequently, we see the director field anchoring to the boundary as the prefactor $\Gamma$ increases in the anchoring integral \eqref{eq:surfaceAnchoringIntegral}. Mimicking results in \cite{Bates2010}, defects appear as the shape elongates. As discussed for Subproblem A (Fixed Shape) \eqref{eq:fsmf}, our procedure does find rotationally invariant solutions. Finally, corroborating results found in \cite{Prinsen2003} and \cite{joshi2024}, Figure \ref{fig:QNNImsmfRho} shows the relationship between aspect ratio of the shape to anchoring strength, $\omega$. We see a strong logarithmic trend in aspect ratio in the presence of changing molecular alignment (full problem) in contrast to a linear trend for fixed alignment (Subproblem B) as predicting by the scaling analysis done in \cite{Prinsen2003}.
	
To compare the performance of GD and QN with and without NI, iteration counts, final energy, and runtimes are included in Table \ref{tab:PGDNIQNNImsmfResults}. Runtime for standalone GD and GD with NI includes time for equiangulation within each grid as numerical experiments indicated it is crucial for convergence. Standalone QN and QN with NI, however, do not require equiangulation within each grid. 
NI yields significant improvement for both GD and QN. We note that sensitivity to spatial deformation depends on grid size, as QN and GD both met the preferred convergence criterion of a relative energy change on the standalone $M_9$ grid simulations across all values $\omega$, but converged to a different energy, approximately $0.1-1\%$ higher than the energies computed with nested iteration.  Thus, NI yields better variational estimates of the true continuum solution. With NI, the iteration counts for the full shape optimization problem dramatically decrease as we refine the grid to $M_9$, guiding both methods to converge to the same final energy $F^k$ (again lower than the energy found from the standalone methods). In addition, while GD with NI still requires several iterations for large $\omega$ on the finest grid, QN with NI converged in only a couple of iterations on the finest grid. Due to hardware memory limitations, we were unable to perform tests on even finer grids.  However, based on the trend, we predict that the number of iterations will continue to decrease until only 1 iteration is performed per fine grid level.  Nevertheless, as shown in Table \ref{tab:PGDNIQNNImsmfResults}, GD with NI significantly improves the timing compared to standalone GD by a factor of $53$ to $219$ times. Moreover, QN with NI improves the standalone QN timings by a factor of $3$ to $52$. The latter speedup factor is lower since standalone QN runs up to $28$ times faster than standalone GD.

 \begin{table}[!ht]
        		\centering
        		\resizebox{\columnwidth}{!}{%
        			\begin{tabular}{cc|cccccc}
        				\hline 
        				\textbf{GD with NI}& & $\omega = 0.01$   & $0.2$  & $0.4$ & $0.6$ & $0.8$ & $1$     \\ \hline
        				NI Grid& $|M_i|$&\multicolumn{6}{c}{Iterations}\\
\hline
$M_1$ & $16$        & $40$          & $517$         & $568$     & $512$        & $336$         & $328$       \\ 
$M_2$ & $49$        & $16$          & $454$         & $341$     & $430$        & $996$         & $240$      \\ 
$M_3$ & $169$       & $8$           & $15$          & $22$      & $95$         & $128$         & $152$   \\ 
$M_4$ & $625$       & $7$           & $16$          & $15$      & $23$         & $24$          & $40$   \\ 
$M_5$ & $2,401$     & $6$           & $8$           & $14$      & $13$         & $15$          & $16$   \\ 
$M_6$ & $9,409$     & $5$           & $7$           & $8$       & $8$          & $8$           & $12$ \\
$M_7$ & $37,249$    & $5$           & $7$           & $7$       & $7$          & $7$           & $7$   \\ 
$M_8$ & $148,225$   & $5$           & $6$           & $7$       & $7$          & $7$           & $7$   \\
$M_9$ & $591,361$   & $5\ (1,151)$   & $5\ (535)$  & $5\ (478)$ & $6\ (2,036)$ & $7\ (1,956)$  & $7\ (2,070)$ \\
        				\hline 
$F^k$ &  & $15.41\ (15.83)$ & $17.69\ (18.78)$ & $19.16\ (21.70)$ & $20.10\ (24.13)$ & $20.71\ (26.63)$ & $21.13\ (29.01)$ \\ 
        				\hline 
\bfseries Runtime [sec]  &  & $ \vec{213.59\ (38,530.90)}$ & $\vec{239.56\ (15,037.50)}$ & $\vec{241.70\ (12,898.20)}$ & $\vec{250.78\ (54,938.78)}$ & $\vec{259.96\ (52,780.08)}$ & $\vec{260.56\ (55,856.21)}$ \\ 
        				\hline
        				\textbf{QN with NI}& & $\omega = 0.01$   & $0.2$  & $0.4$ & $0.6$ & $0.8$ & $1$     \\ \hline
        				NI Grid& $|M_i|$ &\multicolumn{6}{c}{Iterations}\\
        				\hline
$M_1$ &  $16$     & $16$          & $178$         & $161$     & $176$    & $115$      & $111$       \\ 
$M_2$ &  $49$     & $10$          & $25$         & $97$     & $55$    & $163$      & $138$      \\ 
$M_3$ &  $169$     & $4$           & $17$          & $24$     & $68$    & $18$      & $36$   \\ 
$M_4$ &  $625$     & $4$           & $10$           & $9$      & $12$    & $26$        & $18$   \\ 
$M_5$ &  $2,401$     & $3$           & $7$           & $8$       & $10$    & $9$        & $8$   \\ 
$M_6$ &  $9,409$     & $2$           & $7$           & $6$       & $5$     & $6$        & $6$ \\
$M_7$ &  $37,249$     & $2$           & $4$           & $7$       & $6$     & $3$        & $4$   \\ 
$M_8$ &  $148,225$     & $2$           & $2$           & $2$       & $5$     & $3$        & $2$   \\                     
$M_9$ &  $591,361$     & $2\ (36)$ & $2\ (10)$ & $2\ (79)$  &  $2\ (96)$ & $2\ (152)$ &  $2\ (107)$  \\ 
        				\hline 
$F^k$ & & $15.41\ (15.83)$ & $17.69\ (18.93)$ & $19.16\ (21.12)$ & $20.10\ (23.19)$ & $20.71\ (24.56)$ & $21.12\ (28.4271)$ \\ 
        				\hline
\bfseries Runtime [sec] &  & $ \vec{132.89\ (2,053.69)}$ & $ \vec{166.13\ (534.58)}$ & $ \vec{165.99\ (4,013.73)}$ &$ \vec{201.77\ (5,095.24)}$ &$ \vec{169.76\ (8,839.64)}$  & $ \vec{159.02}\ (\vec{5,373.47})$ \\ 
        				\hline 
        			\end{tabular}%
        		}
        		\caption{\footnotesize \justifying Full $(\vec{X},Q)$ $2D$ problem: Iteration count for GD with NI (top) and QN with NI (bottom) on each grid level. Iteration count for methods without NI is given on level $M_9$ in parenthesis. The final energy, $F^k$, and runtime in seconds for the full simulation with NI and without in parenthesis are also given. 
          }
        		\label{tab:PGDNIQNNImsmfResults}
        	   \end{table}

\begin{figure}[!ht]
    \centering
    \begin{subfigure}[t]{0.32\textwidth}
        \centering
        \includegraphics[width=0.9\linewidth]{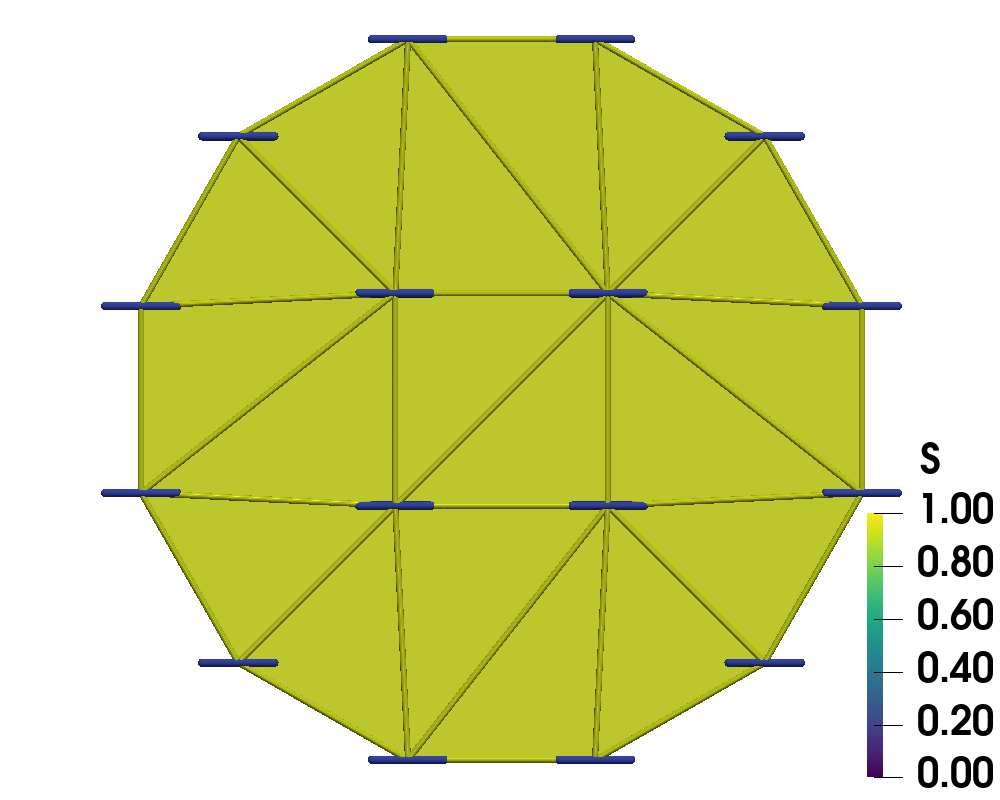}
        \caption{$\vec{n} = (1,0,0)$.}
        \label{fig:CircleShapeHorizontalField}
    \end{subfigure}
    \hfill
    \begin{subfigure}[t]{0.32\textwidth}
        \centering
        \includegraphics[width=0.9\linewidth]{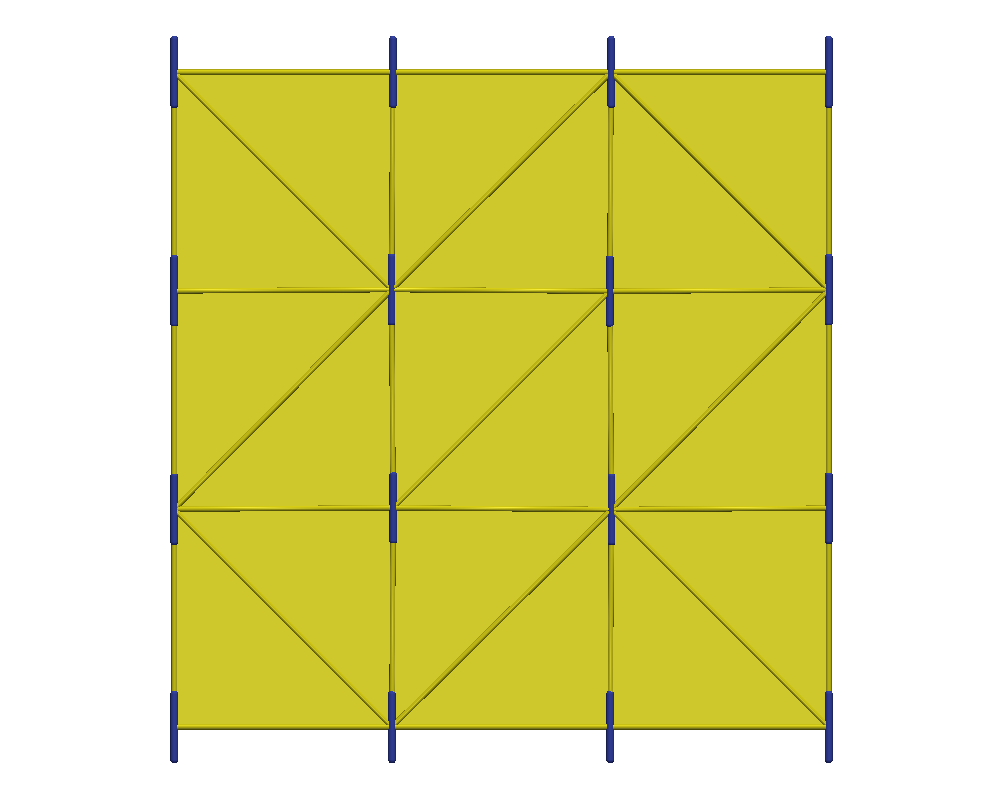}
        \caption{$\vec{n} = (0,1,0)$.}
        \label{fig:RectangleShapeVerticalField}
    \end{subfigure}
    \hfill
    \begin{subfigure}[t]{0.32\textwidth}
        \centering
        \includegraphics[width=0.9\linewidth]{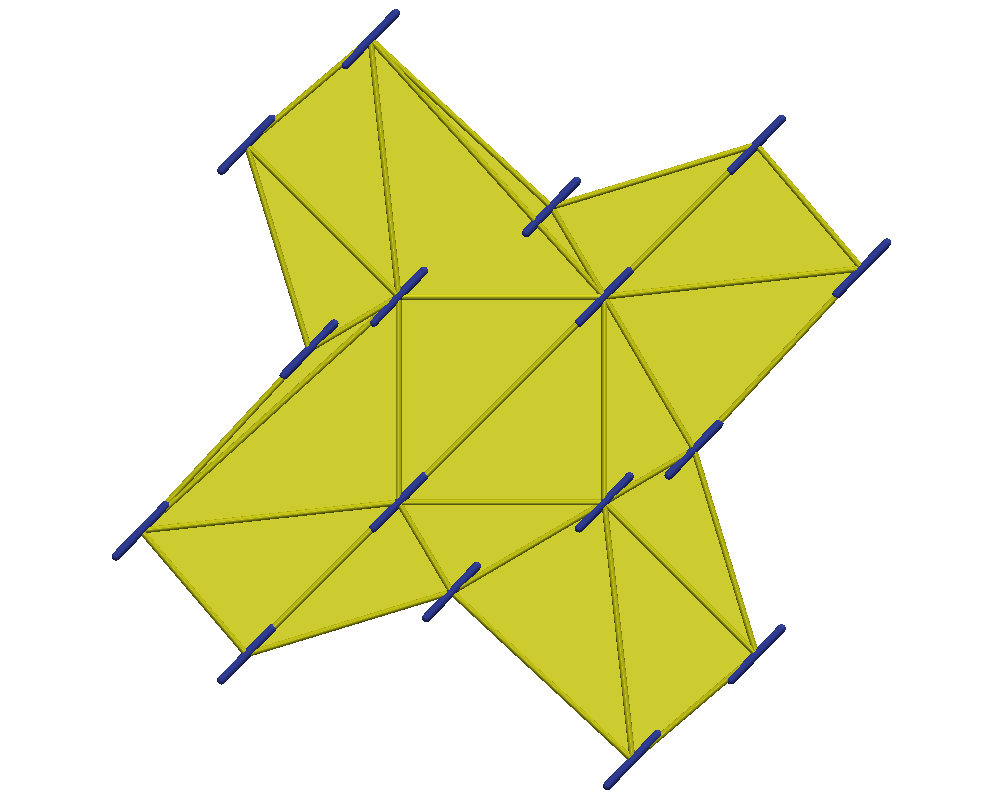}
        \caption{$\vec{n} = \left(\frac{1}{\sqrt{2}}, \frac{1}{\sqrt{2}}, 0\right)$.}
        \label{fig:RadialShapeSlantField}
    \end{subfigure}

    \caption{Initial guesses for the QN with NI scheme on the full $(\vec{X}, Q)$ $2D$ problem.}
    \label{fig:2dIG_msmf}
\end{figure}

Finally, to demonstrate the robustness of our proposed method with respect to the initial guess, we compare results for the QN with NI approach using three distinct starting configurations $(\vec{X},Q)$ all for the test case of $\omega=1$. The first is a circular domain with a regular triangulation and an initial orientation of $\vec{n} = (1,0,0)$, as was used in the previous experiments (see Figure \ref{fig:CircleShapeHorizontalField}). The second is a rectangular domain and an initial orientation of $\vec{n} = (0,1,0)$ (see Figure \ref{fig:RectangleShapeVerticalField}). The final configuration is a star-shaped domain composed of an irregular triangulation with cusps and a rotated orientation of $\vec{n} = (\frac{1}{\sqrt{2}}, \frac{1}{\sqrt{2}}, 0)$ (see Figure(\ref{fig:RadialShapeSlantField})). 
As summarized in Table~\ref{tab:Various2DInitialGuesses}, all three initial guesses converge, with the aid of nested iteration, to the same solutions with energy $F^k = 21.12$, all within comparable computation times.  While not shown, the final configurations are also identical to those shown in Figure~\ref{fig:QNNImsmf}(e) and \ref{fig:QNNImsmf}(f).

\begin{table}[!ht]
\centering
\resizebox{0.65\columnwidth}{!}{%
\begin{tabular}{cc|ccc}
\hline
\multicolumn{2}{c}{Initial Guess} & Figure \ref{fig:2dIG_msmf}(a)&Figure \ref{fig:2dIG_msmf}(b)&Figure \ref{fig:2dIG_msmf}(c)\\ 
\hline
NI Grid& $|M_i|$ &\multicolumn{3}{c}{Iterations}\\ 
\hline
$M_1$& $16$& $111$&$127$&$101$      \\ 
$M_2$& $49$& $138$&$159$&$129$     \\ 
$M_3$& $169$& $36$&$47$&$25$   \\ 
$M_4$& $625$& $18$&$11$&$14$   \\ 
$M_5$& $2,401$& $8$&$8$&$8$  \\ 
$M_6$& $9,409$& $6$&$6$&$5$  \\
$M_7$& $37,249$ &$4$ &$4$&$4$\\
$M_8$& $148,225$ &$2$ &$2$&$3$\\ 
$M_9$& $591,361$  &$2$ &$2$&$2$\\ 
\hline
$F^k$ &  & $21.12$&$21.12$&$21.12$ \\ 
\hline
\textbf{Runtime [sec]}& &$\vec{158.72}$ &$\vec{158.49}$&$\vec{168.56}$   \\ \hline
\end{tabular}
}
\caption{Full $(\vec{X},Q)$ $2D$ problem: Iteration count for QN with NI on each grid level using  $\omega = 1$ and different initial guesses corresponding to Figures \ref{fig:2dIG_msmf}(a)-(c). The final energy, $F^k$, and runtime in seconds for the full simulation are also given.}
\label{tab:Various2DInitialGuesses}
\end{table}


\subsection{Three-Dimensional Nematic Tactoids}\label{section:3dTactoids}
We conclude the numerical results section by simulating the spatial and orientational configuration of a challenging three-dimensional problem involving the formation of a nematic tactoid. 
We use the same material constants as above \cite{Mottram2014}, with $\xi=1\times10^{-7}~m$ and, as before, we
report the distribution of the scalar order parameter, $S$, the corresponding director field, the converged energy, $F^k$, the number of iterations for each grid level of the QN with NI scheme, and the runtime in seconds. For the numerical experiments below, we fix the surface tension, $\tau = 10$, and vary the surface anchoring, $\omega$ from $0.01$ to $0.18$.
The preferred convergence criterion is similarly the relative change in the energy $F^k$, this time with the tolerance set at $10^{-4}.$

Recall that $|M_i|$ denotes the number of vertices on grid $M_i$. Here in three dimensions, we consider 6 different grids of increasing size, ranging from $|M_1| = 27$ to $|M_6| = 275,793$.  Note that there are \textit{eight} degrees of freedom attached to each grid point, and the actual number of vertices on each level might vary depending on the aspect ratio of the current configuration. 

For the first value of $\omega$ ($\omega=0.01$) we consider, the initial guess is defined on a sphere with volume $1$ with equally distributed vertices representing the degrees of freedom, $\vec{X}\in M_1$. The initial director field, $\vec{n}=(1,0,0)$, is aligned parallel to the $x-$axis.
However, as we see numerically, the shape’s aspect ratio is linearly $\omega-$dependent, much stronger than the logarithmic $\omega-$dependence from the two-dimensional tactoid shapes seen earlier. Thus, the initial guess of a sphere is not adequate for convergence as $\omega$ is increased. To mitigate this, we incorporate continuation with respect to $\omega$ on the coarsest grid, $M_1$.  From there, nested iteration is implemented as before.  For example, the final equilibrium state on grid $M_1$ using $\omega=0.01$ is used as the initial guess for the simulation of $\omega=0.02$ on grid $M_1$.  The final state for $\omega=0.02$ on grid $M_1$ is then used for the $\omega=0.04$ simulation on the coarsest grid and so on.  

Table \ref{tab:QNNI3D} depicts the number of QN iterations using NI to build up through the hierarchy of grids, showing the iteration count on each level of refinement. For small $\omega$, while the iterations eventually decrease to one iteration on the finest grid, we see that the coarsest grid of only $27$ points is not enough to accurately represent the true shape and order, corroborated by the higher iteration count on grid $M_2$.  Therefore, for $\omega=0.08$ to $\omega = 0.18$, shown in the bottom of Table \ref{tab:QNNI3D}, we start the nested iteration process using grid $M_2$ as the coarsest grid, and iterate to grid $M_6$ to keep $5$ levels of refinement as before. Concurrently, we use continuation on $\omega$ across grid $M_2$. Again, the results show the iteration count decreasing down the grids as $\omega$ increases, as expected. This iteration trend also implies that modeling these complex configurations requires more points to accurately represent the shape and order.

     \begin{table}[!ht]
        		\centering
        		\resizebox{0.65\columnwidth}{!}{%
        		\begin{tabular}{cc|cccc} \hline 
& & $\omega = 0.01$   & $0.02$  & $0.04$ & $0.06$     \\ \hline
NI Grid& $|M_i|$ &\multicolumn{4}{c}{Iterations}\\ \hline
$M_1$  &  $27$  & $3$ & $1$ & $1$  & $1$ \\ 
$M_2$  &  $125$ & $9$ & $9$ & $8$ & $11$      \\ 
$M_3$  &  $729$ & $4$ & $4$ & $5$  & $9$      \\ 
$M_4$  &  $4,913$  & $4$ & $4$ & $5$ & $5$  \\ 
$M_5$  &  $35,937$ & $1$ & $1$ & $1$ & $1$      \\ \hline 
$F^k$ & & $27.62$ & $27.33$ & $26.74$ & $26.10$ \\      \hline 
\bfseries Runtime [sec]& &$\vec{1,907.29}$&$\vec{1,927.65}$&$\vec{2,111.71}$&$\vec{2,213.35}$  \\ \hline 
        			\end{tabular}%
        		}
          
          \vspace{12pt}
        		\resizebox{0.9\columnwidth}{!}{%
        			\begin{tabular}{cc|cccccc}
        				\hline 
&                       & $\omega=0.08$ & $0.10$    & $0.12$    & $0.14$     & $0.16$    & $0.18$     \\ \hline
NI Grid& $|M_i|$ &\multicolumn{6}{c}{Iterations}\\ \hline
$M_1$ &  $27$ & &          &      & & &         \\ 
$M_2$ &  $125$          & $6$          & $1$       & $15$       & $1$        & $6$     & $6$    \\ 
$M_3$ &  $729$          & $5$          & $12$      & $12$       & $11$       & $6$    & $8$   \\ 
$M_4$ &  $4,913$        & $5$          & $5$       & $5$       & $6$        & $6$     & $7$   \\ 
$M_5$ &  $36,025$       & $3$          & $3$       & $3$       & $3$        & $5$     & $5$    \\
$M_6$ &  $275,793$      & $1$          & $1$       & $1$       & $1$        & $1$     & $3$    \\  \hline 
$F^k$ &                 & $25.44$      & $24.76$   & $24.01$   & $23.28$    & $22.53$   & $21.70$   \\ \hline
\bfseries Runtime [sec]& & $\vec{14,363.60}$ & $\vec{14,340.70}$ & $\vec{14,283.6}$  & $\vec{14,462.24}$ & $\vec{19,891.51}$ &$\vec{36,430.32}$ \\ 
        				\hline 
        			\end{tabular}%
        		}
        		\caption{\footnotesize \justifying Full $(\vec{x},Q)$ $3D$ problem with $\omega = 0.01 \rightarrow 0.06$ (top) and $\omega = 0.08 \rightarrow 0.18$ (bottom): Iteration count for QN with NI on each grid level. The final energy, $F^k$, and runtime in seconds for the full simulation with NI are also given.}
        		\label{tab:QNNI3D}
        	   \end{table}
               
    \begin{figure}[!ht]
		\begin{subfigure}{0.49\textwidth}
			\centering
			\includegraphics[width=0.7\linewidth]{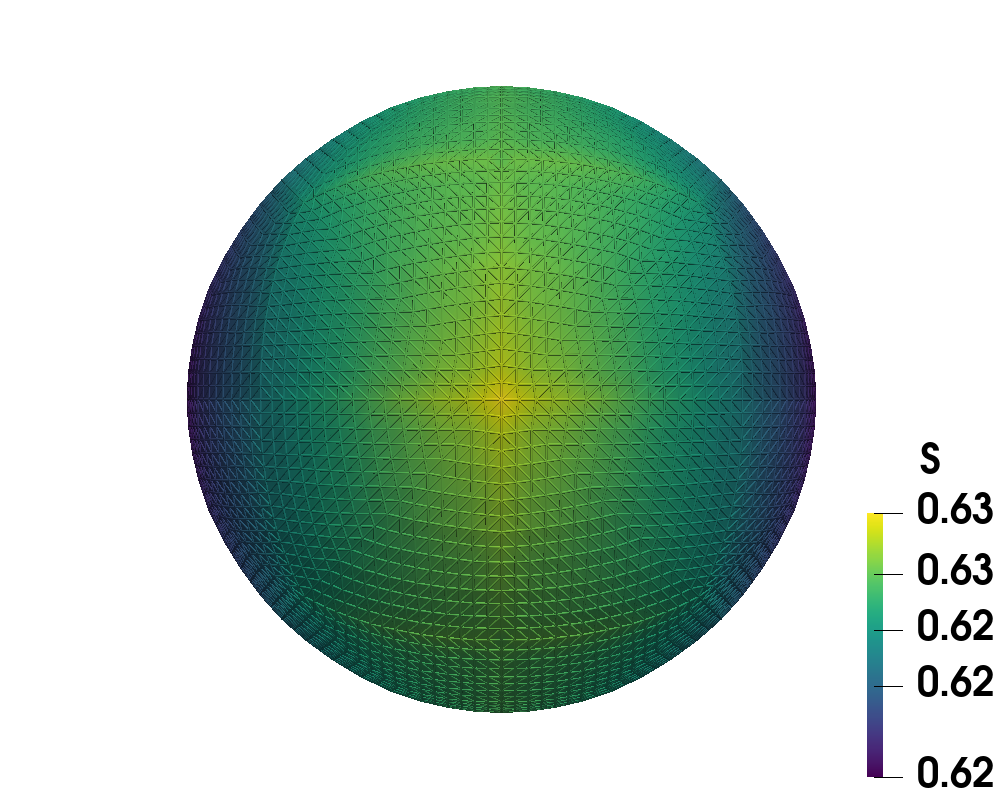}
			\caption*{(a) $\omega = 0.01$}
		\end{subfigure}%
		\begin{subfigure}{0.49\textwidth}
			\centering
			\includegraphics[width=0.7\linewidth]{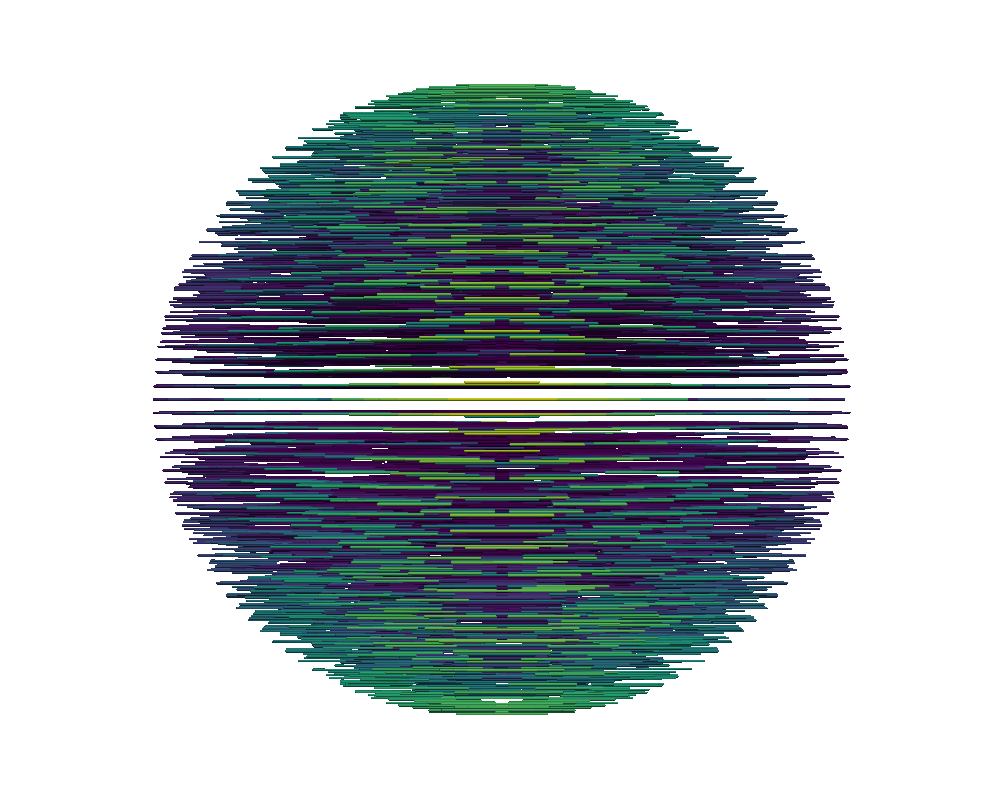}
			\caption*{(b)}
		\end{subfigure}\\
		\begin{subfigure}{0.49\textwidth}
			\centering
			\includegraphics[width=0.7\linewidth]{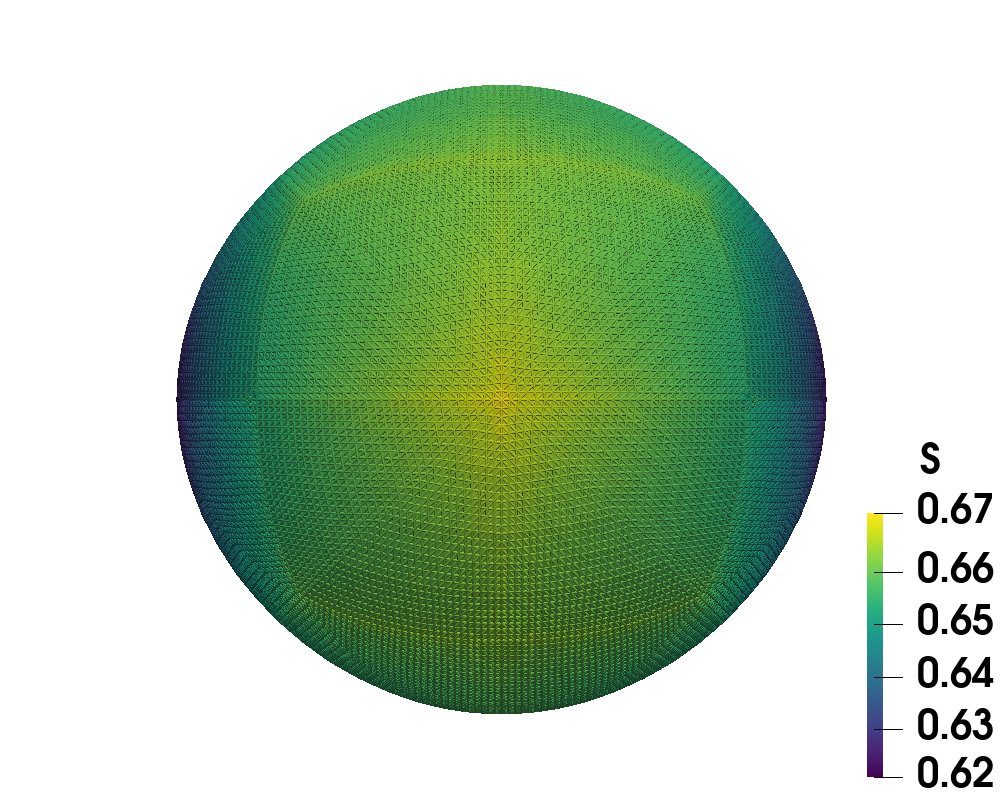}
			\caption*{(c) $\omega = 0.1$}
		\end{subfigure}%
		\begin{subfigure}{0.49\textwidth}
			\centering
			\includegraphics[width=0.7\linewidth]{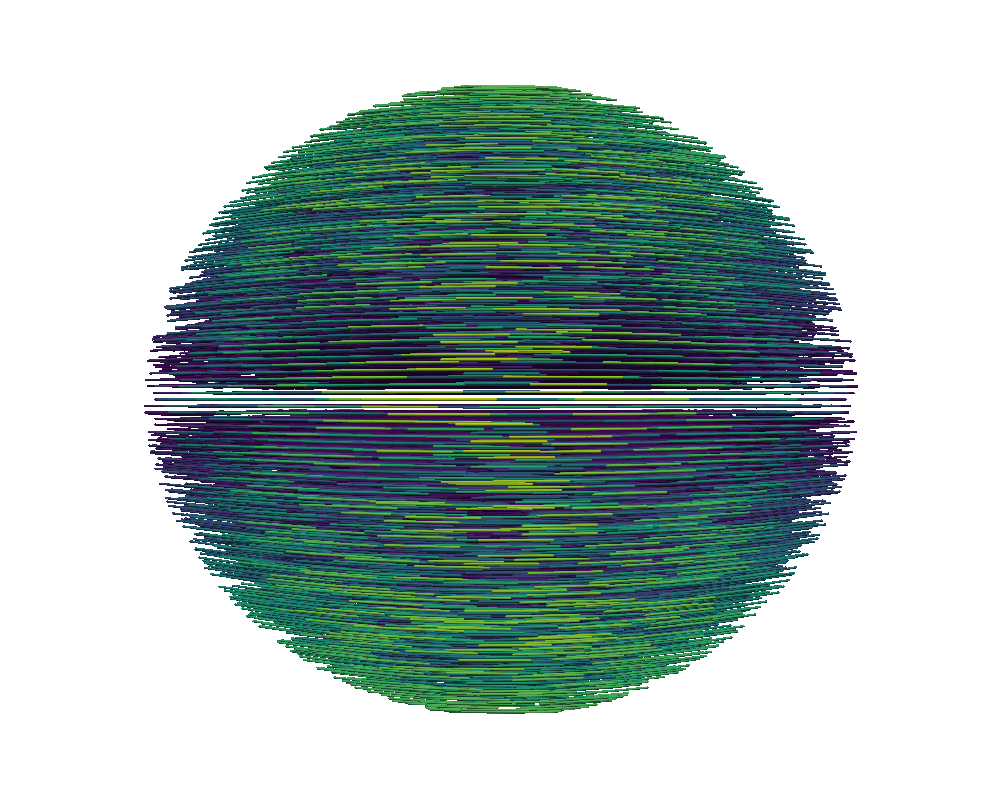}
			\caption*{(d)}
		\end{subfigure}\\
  		\begin{subfigure}{0.49\textwidth}
			\centering
			\includegraphics[width=0.7\linewidth]{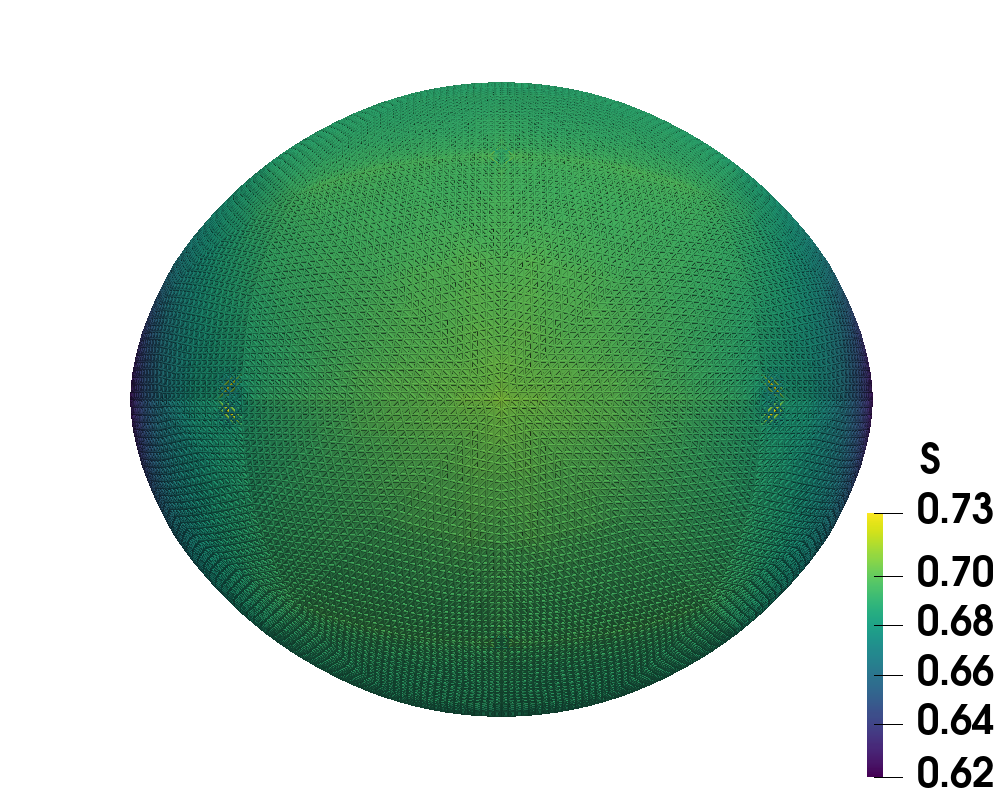}
			\caption*{(e) $\omega = 0.18$}
		\end{subfigure}%
		\begin{subfigure}{0.49\textwidth}
			\centering
			\includegraphics[width=0.7\linewidth]{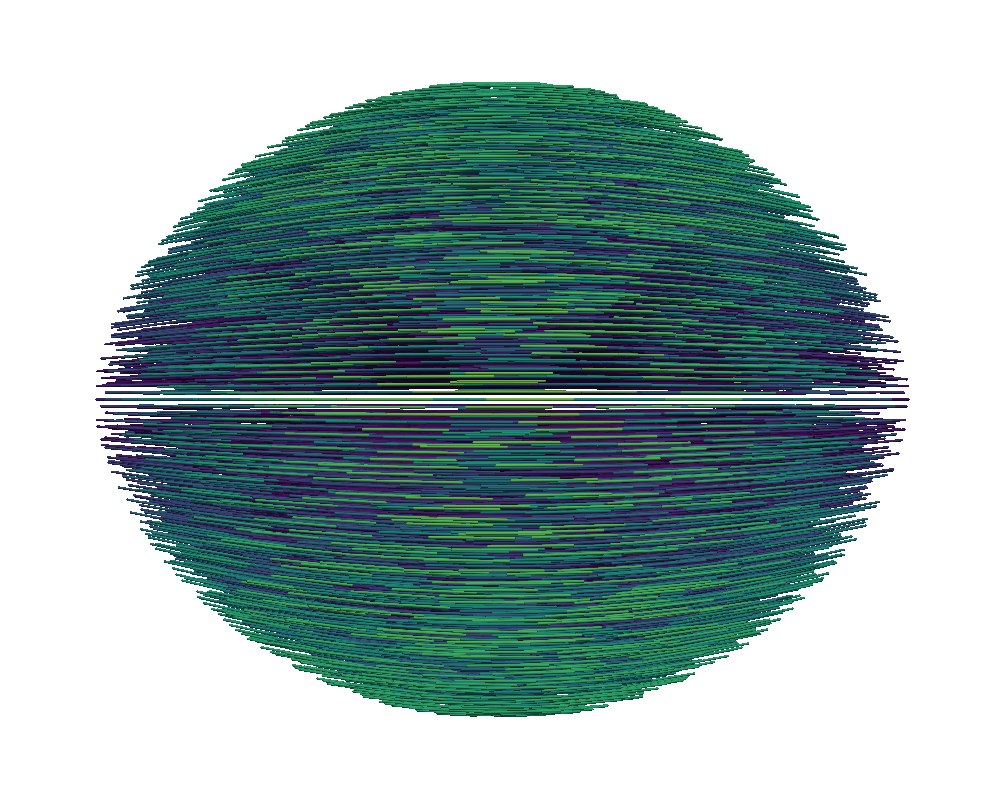}
			\caption*{(f)}
		\end{subfigure}
		\caption{\justifying Applying QN with NI to the  full $(\vec{X},Q)$ $3D$ problem. The color bar indicates the value of $S$, i.e., the order of the director field in the domain. Left plots, (a), (c), (e), depict the order's distribution. Right plots, (b), (d), (f), show the directors with stronger anchoring as the shape changes with $\omega$. Areas in dark green indicate less order, showing the appearance of the defects, as expected.}
		\label{fig:3dQNNImsmfShapes}
	\end{figure}

Figure \ref{fig:3dQNNImsmfShapes} illustrates the nematic tactoids for increasing $\omega$ with the left column showing the order parameter, $S$, indicating the two defects on opposite sides. The right column depicts the strong tangential alignment of the directors to the boundary as $\omega$ increases. 
As $\omega$ increases, the regions of disorder have even less variance indicating that the two defects have localized at the opposite ends of the elongated tactoid. Similar to the two-dimensional simulations, stronger tangential anchoring brought on by higher $\frac{\tau\omega}{2}$ demonstrates significant orientation deformation and sizable spatial horizontal deformation. 
The sharper increase in aspect ratio between the initial spherical state and highly elongated final configurations motivated our use of continuation on $\omega$ along the coarser grids.
    
\begin{figure}[!ht]
        \begin{subfigure}{0.49\textwidth}
			\centering
			\includegraphics[width=0.7\linewidth]{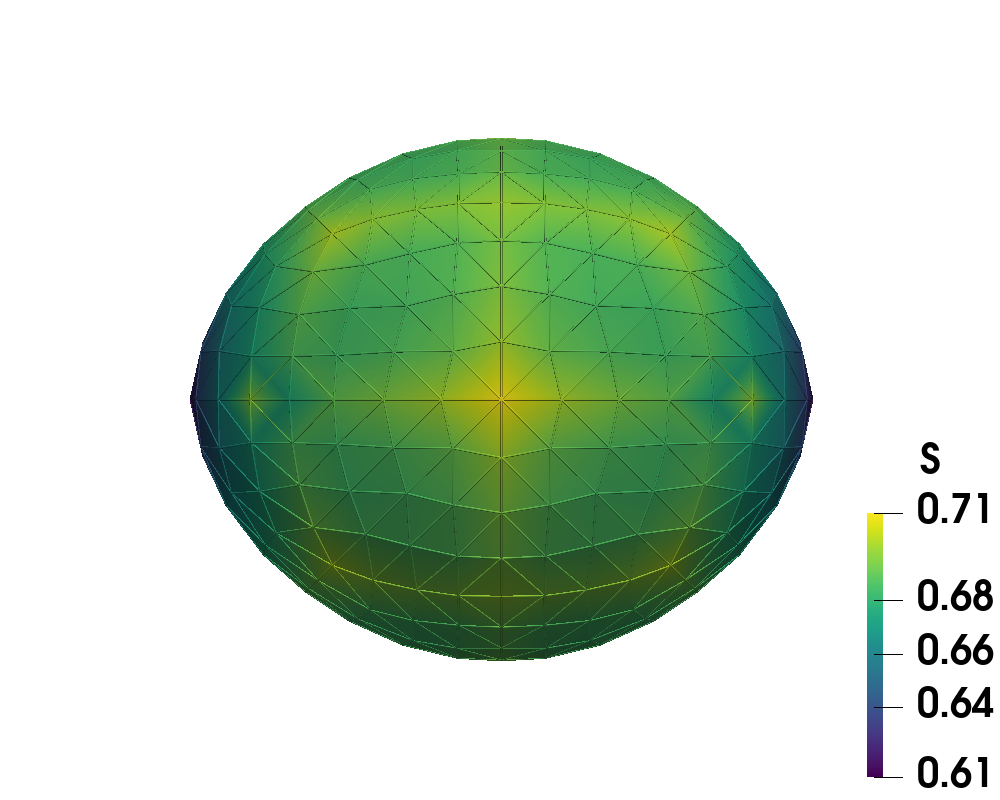}
			\caption*{(a) $\omega = 0.2$}
		\end{subfigure}%
		\begin{subfigure}{0.49\textwidth}
			\centering
			\includegraphics[width=0.7\linewidth]{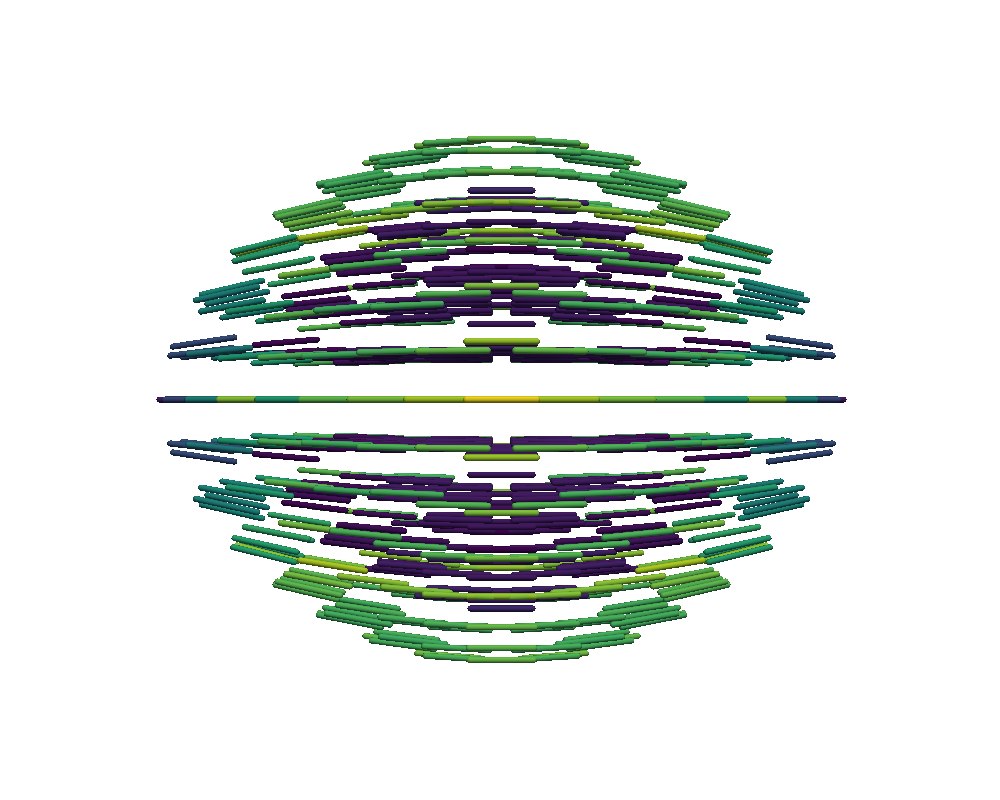}
			\caption*{(b)}
		\end{subfigure}\\
  		\begin{subfigure}{0.49\textwidth}
			\centering
			\includegraphics[width=0.7\linewidth]{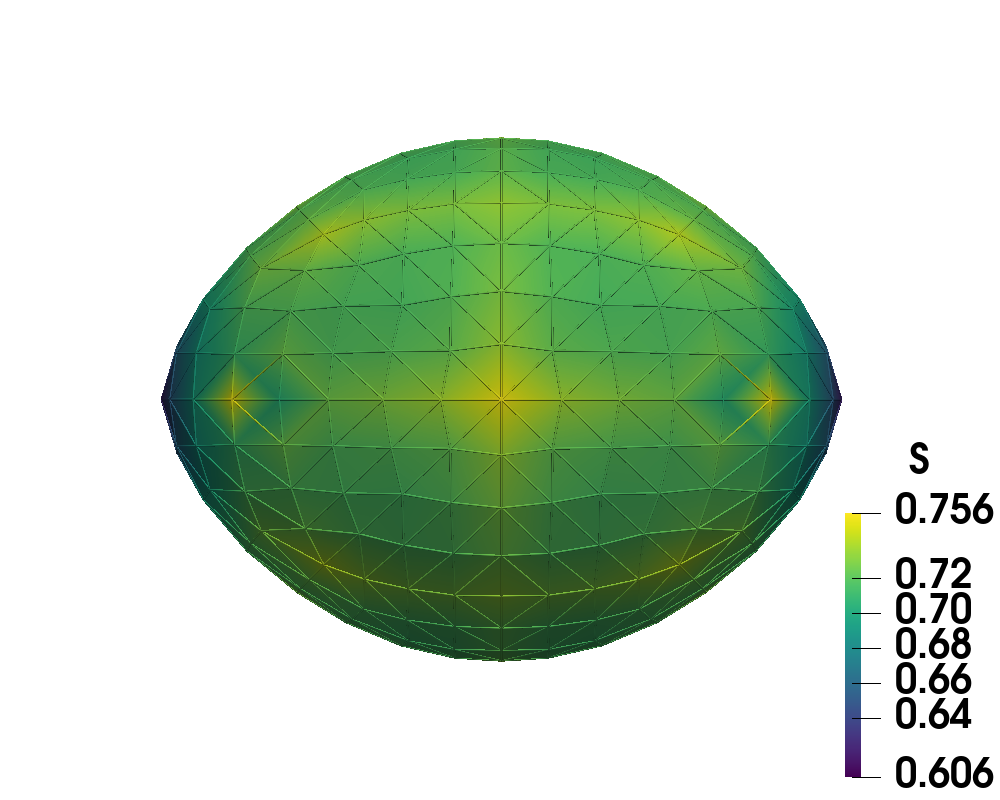}
			\caption*{(c) $\omega = 0.3$}
		\end{subfigure}%
		\begin{subfigure}{0.49\textwidth}
			\centering
			\includegraphics[width=0.7\linewidth]{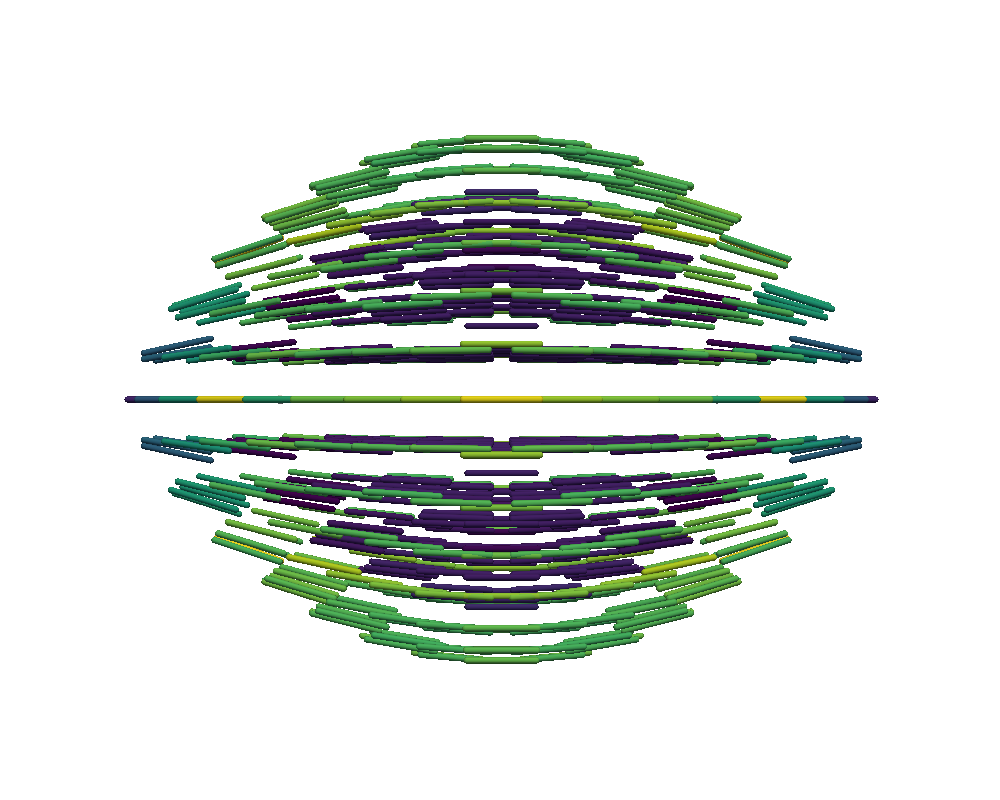}
			\caption*{(d)}
		\end{subfigure}
		\caption{\justifying  Applying QN with NI to the full $(\vec{X},Q)$ $3D$ problem for $\omega=0.2$ and $\omega=0.3$. The color bar indicates the value of $S$, i.e., the order of the director field in the domain. Left plots, (a), (c), depict the order's distribution. Right plots, (b), (d), show the directors with stronger anchoring as the shape changes. Areas in green indicate less order, showing the appearance of the defects, as expected.}
		\label{fig:3dQNNImsmfomega0203}
	\end{figure}

    Finally, we can push the bounds of admissible $\omega$ further to $\omega = 0.2$ and $\omega = 0.3$ by using continuation on all nested iteration levels, not just the coarsest. More specifically, for the first level $M_1$, we use the converged solution on $M_1$ from the previous $\omega$ value as an initial guess, then for the second level, we refine that converged solution to $M_2$ for the current $\omega$, average it with the converged solution on $M_2$ from the previous $\omega$ and use that as the initial guess for $M_2$ for the current $\omega$. We continue this process until the desired level is reached. Here, the desired level is $M_3$. Figure \ref{fig:3dQNNImsmfomega0203} illustrates the tactoid shapes with the order measured by $S$ on the left showing the two defects on the opposite ends of the shape, and the director field on the right strongly anchored tangentially to the shape's boundary. We present this discussion as a precursor to our future work of combining continuation and nested iteration more rigorously.

	\section{Conclusion and Future Work}\label{conclude} 
	The present work describes an ``all-at-once'' quasi-Newton approach to modeling a challenging class of nematic liquid crystal tactoid shape optimization problems, where the equilibrium configuration of the model is found by minimizing a free energy functional with respect to the orientational order and shape of the domain. 
	The approach is effective for this class of problems as it does not require maintenance of mesh quality during the minimization process, has an accurate line search procedure that dynamically updates all unknown variables simultaneously, and allows one to efficiently simulate the solutions on a large scale by uniformly increasing the resolution via nested iteration. 

	Exploring the space of shapes as a function of surface tension and anisotropic elastic constants, we find the nematic tactoids forming under conditions similar to those observed elsewhere \cite{Prinsen2003,Bates2010}. Our main goal here was to improve existing numerical algorithms used to find the equilibrium configurations while ensuring physical validity. We found that through the use of nested iteration, where we gradually refine the initial guess towards a solution with high resolution, we are able to preserve accuracy and robustness with low computational costs.
	
	Future work involves solving the linearized steps iteratively using multigrid methods to reduce the computational cost further while maintaining the same level of accuracy and efficiency. We also plan to further investigate the use of continuation coupled with nested iteration. In addition, we plan to apply the methods discussed in this work to other liquid crystal phases, e.g. cholesterics and smectics, and compare them against experimental results. Finally, including the study of inequality constraints that simulate the formation of nematic tactoids in bounded channels will be considered. 
	
	\section*{Acknowledgments} 
	This work was supported by the National Science Foundation under Grant No. ACI-2003820. The authors of this paper would like to thank Dr. Xiaozhe Hu, Dr. Chaitanya Joshi, and Dr. Viviana Betancur for their valuable discussions. 
	
	\bibliographystyle{elsarticle-num}   
	\bibliography{AAArefs_rev}
	
	
\end{document}